\documentclass[leqno,a4paper,12pt]{article}

\usepackage{amsmath,amsthm}

\newtheorem{Thm}{\indent Theorem}[section]
\newtheorem{Prop}[Thm]{\indent Proposition}
\newtheorem{Lem}[Thm]{\indent Lemma}

\newtheorem{Cor}[Thm]{\indent Corollary}

\theoremstyle{definition}
\newtheorem{Var}[Thm]{\indent Variant}

\newtheorem{Def}[Thm]{\indent Definition}

\newtheorem{Rem}[Thm]{\indent Remark}
\newtheorem{Ex}[Thm]{\indent Example}

\newtheorem{Over}[Thm]{\indent Overview}

\newcommand{\myheading}[1]{\medskip{\bf{#1}}\quad}

\def\qed{{\hskip0pt\unskip\unskip\nobreak\hfil\penalty50
          \hskip1em\hbox{}\nobreak\hfil
          {\bf q.e.d.}%
          \parfillskip=0pt\finalhyphendemerits=0
          \par}\medskip}

\newenvironment{Proof}
               {{\it Proof.}\quad}
               {\qed}

\newenvironment{Proofof}[1]
               {{\it Proof of #1.}\quad}
               {\qed}

\newcommand{\Prime}{\kern3\fontdimen1\font$'$\kern-7\fontdimen1\font}

\long\def\forget#1{}

\long\def\beginSIDEREMARK#1\endSIDEREMARK
    {{\par\bigskip\advance\leftskip by 2cm
                  \advance\rightskip by -2cm\noindent
      {\bf Our own side remark:} #1
      \par\bigskip\noindent}}

\long\def\beginFORGET#1\endFORGET{#1}
\long\def\beginFORGET#1\endFORGET{}

\def\?{\ ???\ \immediate\write16{}%
\immediate\write16{Warning: There was still a question mark . . . }%
\immediate\write16{}}

\usepackage{exscale}
\usepackage{amssymb}

\newcommand{\BA}{{\mathbb{A}}}

\newcommand{\BC}{{\mathbb{C}}}

\newcommand{\BE}{{\mathbb{E}}}

\newcommand{\BG}{{\mathbb{G}}}
\newcommand{\BH}{{\mathbb{H}}}

\newcommand{\BL}{{\mathbb{L}}}
\newcommand{\BM}{{\mathbb{M}}}
\newcommand{\BN}{{\mathbb{N}}}

\newcommand{\BP}{{\mathbb{P}}}
\newcommand{\BQ}{{\mathbb{Q}}}
\newcommand{\BR}{{\mathbb{R}}}
\newcommand{\BS}{{\mathbb{S}}}

\newcommand{\BV}{{\mathbb{V}}}
\newcommand{\BW}{{\mathbb{W}}}
\newcommand{\BX}{{\mathbb{X}}}

\newcommand{\BZ}{{\mathbb{Z}}}

\newcommand{\Fs}{{\mathfrak{s}}}

\newcommand{\FF}{{\mathfrak{F}}}

\newcommand{\FH}{{\mathfrak{H}}}

\newcommand{\FS}{{\mathfrak{S}}}
\newcommand{\FT}{{\mathfrak{T}}}
\newcommand{\FU}{{\mathfrak{U}}}
\newcommand{\FV}{{\mathfrak{V}}}

\newcommand{\FX}{{\mathfrak{X}}}

\newcommand{\FZ}{{\mathfrak{Z}}}

\newcommand{\CA}{{\cal A}}

\newcommand{\CD}{{\cal D}}

\newcommand{\CF}{{\cal F}}

\newcommand{\CH}{{\cal H}}

\newcommand{\Spec}{\mathop{{\bf Spec}}\nolimits}

\newcommand{\Gr}{{\rm Gr}}

\newcommand{\imm}{\mathop{{\rm im}}\nolimits}

\newcommand{\Hom}{\mathop{\rm Hom}\nolimits}

\newcommand{\Ext}{\mathop{\rm Ext}\nolimits}

\newcommand{\Cent}{\mathop{\rm Cent}\nolimits}

\newcommand{\Stab}{\mathop{\rm Stab}\nolimits}

\newcommand{\Res}{\mathop{\rm Res}\nolimits}
\newcommand{\Ind}{\mathop{\rm Ind}\nolimits}

\newcommand{\Pro}{\mathop{\rm Pro}\nolimits}

\newcommand{\coker}{\mathop{\rm Coker}\nolimits}
\newcommand{\kernel}{\mathop{\rm Ker}\nolimits}
\newcommand{\Int}{\mathop{\rm int}\nolimits}

\newcommand{\loccit}{[loc.$\;$cit.]}

\newcommand{\Ob}{\mathop{\rm Ob}\nolimits}
\newcommand{\Id}{\mathop{\rm Id}\nolimits}
\newcommand{\For}{\mathop{\rm Fo}\nolimits}

\def\tei{\, | \,}
\def\halb{\frac{1}{2}}

\def\id{{\rm id}}

\newbox\mybox
\def\arrover#1{\mathrel{
       \setbox\mybox=\hbox spread 1.4em{\hfil$\scriptstyle#1$\hfil}
       \vbox{\offinterlineskip\copy\mybox
             \hbox to\wd\mybox{\rightarrowfill}}}}
\def\larrover#1{\mathrel{
       \setbox\mybox=\hbox spread 1.4em{\hfil$\scriptstyle#1$\hfil}
       \vbox{\offinterlineskip\copy\mybox
             \hbox to\wd\mybox{\leftarrowfill}}}}

\def\ontoover#1{\mathrel{
       \setbox\mybox=\hbox spread 1.4em{\hfil$\scriptstyle#1$\hfil}
       \vbox{\offinterlineskip\copy\mybox
             \hbox to\wd\mybox{\rightarrowfill\hskip-2.8mm
                               $\rightarrow$}}}}
\def\leftontoover#1{\mathrel{
       \setbox\mybox=\hbox spread 1.4em{\hfil$\scriptstyle#1$\hfil}
       \vbox{\offinterlineskip\copy\mybox
             \hbox to\wd\mybox{$\leftarrow$\hskip-2.8mm
                               \leftarrowfill}}}}
\def\longto{\longrightarrow}
\def\into{\hookrightarrow}
\def\onto{\ontoover{\ }}
\def\longonto{\ontoover{\ }}
\def\isoto{\arrover{\sim}}

\def\longinto{\lhook\joinrel\longrightarrow}
\def\leftlonginto{\longleftarrow\joinrel\rhook}

\usepackage[curve,matrix,arrow,cmtip]{xy}
\NoComputerModernTips

\def\myxymessage{\def\messagetext
   {Here an xy-pic diagram was omitted to speed up compilation . . . }
   \immediate\write16{\messagetext}
   \hbox{\bf \messagetext}}
\def\filxymatrix#1{\myxymessage}
\def\filxyarray#1{\myxymessage}

\newdir^{ (}{{}*!/-3pt/\dir^{(}}
\newdir_{ (}{{}*!/-3pt/\dir_{(}}
\newdir^{ )}{{}*!/+3pt/\dir^{)}}
\newdir_{ )}{{}*!/+3pt/\dir_{)}}

\def\rscript#1{\hbox to 0pt{$\scriptstyle#1$\hss}}

\newcommand{\topp}{\rm{top}}

\newcommand{\is}{\mathop{i_{\sigma}}\nolimits}

\newcommand{\ist}{\mathop{i_{\tilde{\sigma}}}\nolimits}
\newcommand{\isc}{\mathop{i^{\circ}_{\sigma}}\nolimits}

\newcommand{\istc}{\mathop{i^{\circ}_{\tilde{\sigma}}}\nolimits}
\newcommand{\starT}{\mathop{star_{\FT}}\nolimits}
\newcommand{\Mp}{\mathop{M^{\pi(K_1)}}\nolimits}
\newcommand{\Mpd}{\mathop{\Delta \backslash M^{\pi(K_1)}}\nolimits}
\newcommand{\Mo}{\mathop{M^K_1}\nolimits}
\newcommand{\MoS}{\mathop{M^K_{1,\FS}}\nolimits}
\newcommand{\MKS}{\mathop{M^K (\FS)}\nolimits}
\newcommand{\MKo}{\mathop{M^{K_1}}\nolimits}
\newcommand{\MKoS}{\mathop{M^{K_1} (\FS_1^0)}\nolimits}
\newcommand{\Ms}{\mathop{Z_\sigma}\nolimits}
\newcommand{\MI}{\mathop{Z_I}\nolimits}
\newcommand{\MIb}{\mathop{Z_{\bar{I}}}\nolimits}

\newcommand{\Msb}{\mathop{Z_{\bar{\sigma}}}\nolimits}
\newcommand{\Mst}{\mathop{Z_{\tilde{\sigma}}}\nolimits}
\newcommand{\Msc}{\mathop{Z^\circ_\sigma}\nolimits}

\newcommand{\Mstc}{\mathop{Z^\circ_{\tilde{\sigma}}}\nolimits}
\newcommand{\dU}{\mathop{\partial\FU}\nolimits}

\newcommand{\Pnu}{\mathop{\BP (n)_u}\nolimits}
\newcommand{\Gm}{\mathop{\BG_m}\nolimits}
\newcommand{\Gnm}{\mathop{\BG_m^n}\nolimits}

\newcommand{\mup}{\mathop{\mu_{\pi(K_1)}}\nolimits}
\newcommand{\mupT}{\mathop{\mu_{\pi(K_1)}^{\FT}}\nolimits}

\newcommand{\HQ}{\mathop{\overline{H}_Q}\nolimits}
\newcommand{\HC}{\mathop{\overline{H}_C}\nolimits}

\newcommand{\Rep}{\mathop{\bf Rep}\nolimits}
\newcommand{\Sh}{\mathop{\bf Sh}\nolimits}

\newcommand{\Loc}{\mathop{\bf Loc}\nolimits}
\newcommand{\Perv}{\mathop{\bf Perv}\nolimits}
\newcommand{\MHS}{\mathop{\bf MHS}\nolimits}
\newcommand{\MHM}{\mathop{\bf MHM}\nolimits}
\newcommand{\VMHS}{\mathop{\bf Var}\nolimits}

\newcommand{\Tot}{{\rm Tot}}

\newcommand{\can}{\mathop{\rm can}\nolimits}
\newcommand{\var}{\mathop{\rm Var}\nolimits}
\newcommand{\tcan}{\widetilde{\can}}
\newcommand{\tvar}{\widetilde{\var}}
\newcommand{\tE}{\widetilde{E}}
\newcommand{\tH}{\widetilde{H}}
\newcommand{\tBH}{\widetilde{\BH}}
\newcommand{\tN}{\widetilde{N}}
\newcommand{\tS}{\widetilde{S}}

\newcommand{\Ab}{\mathop{\CA\it b}\nolimits}

\let\oldbullet\bullet
\def\bullet{{\mathchoice{\oldbullet}%
                        {\oldbullet}%
                        {\scriptscriptstyle\oldbullet}%
                        {\oldbullet}}}

\newcommand{\argdot}{{\;\bullet\;}}

\begin{document}

%

\hfuzz=3pt
\overfullrule=10pt                   


\setlength{\abovedisplayskip}{6.0pt plus 3.0pt}
\setlength{\belowdisplayskip}{6.0pt plus 3.0pt}
\setlength{\abovedisplayshortskip}{6.0pt plus 3.0pt}
\setlength{\belowdisplayshortskip}{6.0pt plus 3.0pt}

\setlength{\baselineskip}{13.0pt}
\setlength{\lineskip}{0.0pt}
\setlength{\lineskiplimit}{0.0pt}

%
%

\title{Hodge modules on Shimura varieties and their higher direct images
in the Baily--Borel compactification
\footnotemark
\footnotetext{To appear in Ann.\ scient.\ ENS.}
}
\author{\footnotesize by\\ \\
\mbox{\hskip-2cm
\begin{minipage}{6cm} \begin{center} \begin{tabular}{c}
Jos\'e I.~Burgos  \footnote{
Partially supported by Grant DGI BFM2000-0799-C02-01.}
\footnote{
Partially supported by the European Network ``Arithmetic Algebraic Geometry''.}\\[0.2cm]
\footnotesize Departament d'Algebra i Geometria\\[-3pt]
\footnotesize Universitat de Barcelona\\[-3pt]
\footnotesize Gran Via 585\\[-3pt]
\footnotesize E-08007 Barcelona\\[-3pt]
\footnotesize Spain\\
{\footnotesize \tt burgos@mat.ub.es}
\end{tabular} \end{center} \end{minipage}
\begin{minipage}{6cm} \begin{center} \begin{tabular}{c}
J\"org Wildeshaus $^\ddag$  \\[0.2cm]
\footnotesize Institut Galil\'ee\\[-3pt]
\footnotesize Universit\'e Paris 13\\[-3pt]
\footnotesize Avenue Jean-Baptiste Cl\'ement\\[-3pt]
\footnotesize F-93430 Villetaneuse\\[-3pt]
\footnotesize France\\
{\footnotesize \tt wildesh@math.univ-paris13.fr}
\end{tabular} \end{center} \end{minipage}
\hskip-2cm}
}
\date{January 7, 2004}
\maketitle
\quad \\[-1.7cm]
\begin{abstract}
\noindent
We study the degeneration in the Baily--Borel compactification
of va\-ria\-tions of Hodge structure on Shimura varieties.
Our main result Theorem~\ref{2D}
expresses the degeneration of
variations given by algebraic representations in terms of
Hochschild, and abstract group cohomo\-lo\-gy. It
is the Hodge theoretic analogue of Pink's theorem on degeneration
of \'etale and $\ell$-adic sheaves \cite{P2}, and completes 
results by Harder and Looijenga--Rapoport \cite{Hd,LR}.
The induced formula on the level of singular cohomology
is equivalent to the theorem of Harris--Zucker on the Hodge structure
of deleted neighbourhood cohomology of strata in toroidal
compactifications \cite{HZ1}.\\

\forget{
}

\noindent
{\bf R\'esum\'e~:} Ce travail concerne la d\'eg\'en\'erescence des
variations de structure de Hodge sur les vari\'et\'es de Shimura.
Le r\'esultat principal Thm.~\ref{2D} exprime cette d\'eg\'en\'erescence en termes de
cohomologie de Hochschild, et de cohomologie abstraite des groupes.
Ce r\'esultat est l'analogue, en th\'eorie de Hodge,
du th\'eor\`eme de Pink sur la d\'eg\'en\'erescence des faisceaux \'etales et
$\ell$-adiques \cite{P2}, et compl\`ete des r\'esultats obtenus par
Harder et Looijenga--Rapoport \cite{Hd,LR}.
Il induit une formule au niveau de la cohomologie singuli\`ere, qui est
\'equivalente au th\'eor\`eme de Harris--Zucker concernant la structure de
Hodge sur la ``deleted neighbourhood cohomology'' des strates des
compactifications toro\"{\i}dales \cite{HZ1}.\\
\end{abstract}


\bigskip
\bigskip
\bigskip

\forget{
\noindent {\footnotesize Math.\ Subj.\ Class.\ (2000) numbers:
14G35 (11F75, 14D07, 14F25, 14F40, 14F43, 32G20).
}
}

\tableofcontents


\newpage
%
%

\setcounter{section}{-1}
\section{Introduction}
\label{Intro}



In this paper, we consider the \emph{Baily--Borel
compactification} of a (pure) \emph{Shi\-mu\-ra variety}
\[
j: M \longinto M^* \; .
\]
According to \cite{P1}, the \emph{boundary} $M^* - M$ has a natural
\emph{stratification} into locally closed subsets,
each of which is itself (a quotient by
the action of a finite group of) a Shimura variety. Let
\[
i: M_1 \longinto M^*
\]
be the inclusion of an individual such stratum.
Saito's theory of \emph{mixed algebraic Hodge modules} \cite{Sa}
comes equipped with the formalism of \emph{Grothen\-dieck's functors}.
In particular, there is a functor
\[
i^*j_*
\]
from the bounded derived category of Hodge modules
on $M$ to that of Hodge modules on $M_1$. We shall refer to this functor
as the \emph{degeneration of Hodge modules on $M$ along the stratum $M_1$}.\\

The objective of the present article is a formula for the effect of $i^*j_*$
on those complexes of
Hodge modules coming about via the \emph{canonical construction},
denoted $\mu$:
the Shimura variety $M$ is associated to a reductive group
$G$ over $\BQ$, and any complex of algebraic representations
$\BV^\bullet$ of $G$ gives rise
to a complex of Hodge modules
$\mu(\BV^\bullet)$ on $M$. Let $G_1$ be the group belonging to $M_1$;
it is the maximal reductive quotient of a certain
subgroup $P_1$ of $G$:
\[
\begin{array}{ccccc}
W_1 & \trianglelefteq & P_1 & \subset & G  \\
& & \downarrow & & \\
& & G_1 & &
\end{array}
\]
($W_1:=$ the unipotent radical of $P_1$.)
The topological inertia group of $M_1$ in $M$ is an extension of
a certain arithmetic group $\HC$ by a lattice in $W_1 (\BQ)$.\\

Our main result Theorem~\ref{2D} expresses $i^*j_* \circ \mu$
as a composition of \emph{Hochschild
cohomology} of $W_1$, \emph{abstract cohomology} of $\HC$,
and the ca\-no\-ni\-cal construction on $M_1$.
It completes results of Harder and of Looijenga--Rapoport;
in fact, the result on the level of local systems is proved
in \cite{Hd}, while the result ``modulo Hodge filtrations'' is basically
contained in \cite{LR}.
Our result induces a comparison statement
on the level of singular cohomo\-lo\-gy of $M_1$, which is equivalent to one of the
main results of \cite{HZ1}.
Theorem~\ref{2D} is the analogue of the main result of \cite{P2},
which identifies the degeneration of \'etale and $\ell$-adic sheaves.\\

As far as the proof of our main result is concerned,
our \emph{geometric} approach is very close to the one
employed in \cite{P2} and \cite[Sections~4 and 5]{HZ1}:
as there, we use a \emph{toroidal
compactification}, to
reduce a difficult local calculation to an easier local calculation,
together with a global calculation on the fibres of the projection
from the toroidal compactification.
By contrast, the \emph{homo\-lo\-gi\-cal}
aspects differ drastically from \cite{P2}. The reason for this
lies in the behaviour of the formalism of Grothendieck's
functors on the two sheaf categories with respect to the
$t$-structures: roughly speaking, on the $\ell$-adic side,
the functors on the level of derived categories are obtained
by right derivation of (at worst) left exact functors.
Since the same is true for group cohomology,
the formalism of equivariant $\ell$-adic sheaves can be controlled
via the standard techniques using injective resolutions
\cite[Section~1]{P2}. Due to the perverse nature of
Hodge modules, there are no exactness properties for Grothendieck's
functors associated to arbitrary morphisms. Even when half exactness
is known (e.g., right exactness for the inverse image of a closed
immersion, left exactness for the (shift by $-d$ of the) direct
image of a smooth morphism of constant relative dimension $d$),
the corresponding functor on the level of derived categories is
not a priori obtained by derivation. As a consequence, we found
ourselves unable to establish the full formalism of Grothendieck's functors
for equivariant
Hodge modules, except for some almost obvious results when the
action of the group is free (see Section~\ref{3b}).
It turns out that these are sufficient
for our purposes, once we observe that certain combinatorial aspects
of the toroidal compactification can be translated
into group cohomo\-lo\-gy. \\

Talking about group cohomology, we should mention that
to find the correct conceptual context for
the statement of Theorem~\ref{2D} turned out
to be a major challenge in itself: recall that we express $i^*j_* \circ \mu$
as a composition of Hochschild cohomology, abstract group cohomology,
and the ca\-no\-ni\-cal construction on $M_1$.
Due to the nature of the canonical construction,
it is necessary
for abstract cohomology to map algebraic representations to
algebraic representations. We found it most na\-tu\-ral
to develop a formalism of
group cohomology ``in Abelian categories'',
which on the one hand applies to a sufficiently general si\-tu\-ation,
and on the other hand is compatible with usual group cohomology ``in the
category of $\BZ$-modules''.
This is the content of  Section~\ref{3}. \\

In \cite{HZ1}, the authors study the \emph{Hodge structure} on the
boundary cohomology of the \emph{Borel--Serre compactification} $\bar{M}$.
Their main result states that the
\emph{nerve spectral sequence} associated to the natural stratification of
$\bar{M}$ is a spectral sequence of mixed Hodge structures.
Given the non-algebraic nature of $\bar{M}$ and its strata, one of the
achievements of loc.\ cit.\ is to \emph{define} the Hodge structures in
question. It turns out that the $E_1$-terms are
given by \emph{deleted neighbourhood cohomology}
of certain strata in the toroidal compactification. Its
Hodge structure is identified in \cite[Thm.~(5.6.10)]{HZ1}.
We are able to recover this latter result,
for maximal parabolic subgroups ($R = P$ in the notation of \loccit),
by applying singular cohomology
to the comparison isomorphism of Theorem~\ref{2D}. Although Theorem~\ref{2D}
is not a formal consequence of the main results of \cite{HZ1},
it is fair to say that an important part of the local information needed
in our proof is already contained in loc.\ cit.; see also \cite[4.3]{HZ2},
where some of the statements of \cite[Section~5]{HZ1}
are strengthened. Roughly speaking, the fundamental difference between
the approach of loc.\ cit.\ and ours is the following: loc.\ cit.\ uses
the explicit description of the objects in order to deduce a compa\-ri\-son
result. We derive the comparison result from the abstract
properties of the categories involved; this gives in particular an explicit
description, which turns out to be compatible with that of \cite{HZ1,HZ2}. \\

Our article is structured as follows: we assemble the notations and results
necessary for the statement of Theorem~\ref{2D}
in Sections~\ref{1} and \ref{2}.
Sections~\ref{3}--\ref{4} contain
the material needed in its proof, which is
given in Section~\ref{5}. We refer to Overview~\ref{2L} for a more detailed
description. Let us note that because of
the homological difficulties
mentioned further above, one is forced to identify $i^*j_* \circ \mu$ with
the composition of a certain number of functors, each of which is simultaneously
(1)~relatively easy to handle, and (2)~at least half exact.
This explains the central role played by
the \emph{specialization functor} in the context of Hodge modules
(see Section~\ref{3c}). It also explains
the use of {\v{C}ech coverings} in the computation of
the direct image associated to the projection from the toroidal
to the Baily--Borel compactification (see Section~\ref{3aa}). \\

We would like to thank A.~Abbes, F.~Guill\'en,
L.~Illusie, A.~Mokrane, Ngo~B.C., and
J.~Tilouine for useful
discussions. \\

\myheading{Notations and Conventions:}
All Shimura varieties are defined over the field of complex numbers $\BC$.
Throughout the whole article, we make consistent use of the language and
the main results of \cite{P1}.
Algebraic re\-pre\-sentations of an algebraic group are finite dimensional
by definition. If a group $G$ acts on $X$, then we write $\Cent_G X$ for
the kernel of the action. If $Y$ is a sub-object of $X$, then $\Stab_G Y$
denotes the subgroup of $G$ stabilizing $Y$.
\forget{
If $X$ is a variety over $\BC\,$,
then $D^b_\con (X(\BC))$ denotes the full triangulated sub-category
of complexes of sheaves of Abelian groups on $X(\BC)$ with
constructible cohomology. The sub category of complexes whose cohomology
is \emph{algebraically} constructible is denoted by $D^b_\con (X)$.
If $F$ is a coefficient field,
then we define triangulated categories of complexes of sheaves of
$F$-vector spaces
\[
D^b_\con (X , F) \subset D^b_\con (X(\BC) , F)
\]
in a similar fashion. The category
$\Perv_F X$ is defined as the heart of the perverse
$t$-structure on $D^b_\con (X , F)$.\\
}
Finally, the ring of finite ad\`eles over $\BQ$ is denoted by $\BA_f$.


\bigskip
%
%

\section{Strata in the Baily--Borel compactification}
\label{1}

This section is intended for reference. We recall and prove
what is stated in \cite[(3.7)]{P2}, but using group actions from the
left (as in \cite[6.3]{P1}). This can be seen as an adelic
version of the
description of these actions
contained in \cite[(6.1)--(6.2)]{LR}. \\

Let $(G,\FH)$ be \emph{mixed Shimura data} \cite[Def.~2.1]{P1}.
The \emph{Shimura varieties} associated to $(G,\FH)$ are
indexed by the
open compact subgroups of $G (\BA_f)$. If $K$ is one such, then the
analytic
space of $\BC$-valued points of the corresponding variety
$M^K := M^K (G,\FH)$ is given as
\[
M^K (\BC) := G (\BQ) \backslash ( \FH \times G (\BA_f) / K ) \; .
\]
We assume
that $G$ is reductive, and hence that $(G,\FH)$ is \emph{pure} in the
sense of \cite{P1}.
In order to describe the \emph{Baily--Borel compactification}
$(M^K)^*$ of $M^K$ \cite{BaBo,AMRT}, recall that for any
\emph{admissible parabolic subgroup} $Q$ of $G$ \cite[Def.~4.5]{P1},
there is associated a
canonical normal subgroup $P_1$ of $Q$ \cite[4.7]{P1}.
There is a finite collection
of \emph{rational boundary components} $(P_1 , \FX_1)$ associated to
$P_1$, and indexed by the
$P_1 (\BR)$-orbits in $\pi_0 (\FH)$ \cite[4.11]{P1}.
The $(P_1 , \FX_1)$ are themselves mixed Shimura data. Denote by $W_1$
the unipotent radical of $P_1$, and by $(G_1, \FH_1)$ the quotient
of $(P_1 , \FX_1)$ by $W_1$ \cite[Prop.~2.9]{P1}.
From the proof
of \cite[Lemma~4.8]{P1}, it follows that $W_1$
equals the unipotent radical of $Q$. \\

One defines
\[
\FH^* := \coprod_{(P_1 , \FX_1)} \FH_1 \; ,
\]
where the disjoint union is extended over all rational boundary components
$(P_1 , \FX_1)$. This set comes equipped with the \emph{Satake topology}
(see \cite[p.~257]{AMRT}, or \cite[6.2]{P1}), as well as a natural action
of the group $G(\BQ)$ (see \cite[4.16]{P1}). One defines
\[
M^K (G,\FH)^*(\BC) := G (\BQ) \backslash ( \FH^* \times G (\BA_f) / K ) \; .
\]
This object is endowed with the quotient topology. By
\cite[10.4, 10.11]{BaBo} (whose proof works equally well in the more general
context considered by Pink; see \cite[8.2]{P1}),
it can be canonically identified
with the space of $\BC$-valued points of a normal projective complex variety
$(M^K)^* = M^K (G,\FH)^*$, containing $M^K$ as a Zariski-open dense subset.
The stratification of $\FH^*$ induces a stratification of $(M^K)^*$.
Let us explicitly describe this stratification:
fix an admissible parabolic subgroup $Q$ of $G$, and let $(P_1,\FX_1)$,
$W_1$, and
\[
\pi: (P_1,\FX_1) \longonto (G_1,\FH_1) = (P_1,\FX_1)/W_1
\]
as above. Fix an open compact subgroup $K \subset G (\BA_f)$, and an element
$g \in G(\BA_f)$. Define
$K' := g \cdot K \cdot g^{-1}$, and
$K_1  :=  P_1 (\BA_f) \cap K'$.
We have the following natural morphisms (cmp.\ \cite[(3.7.1)]{P2}):
\[
\vcenter{\xymatrix@R-10pt{
M^{\pi(K_1)}(G_1,\FH_1) (\BC) \ar@{=}[d] & \\
G_1(\BQ) \backslash ( \FH_1 \times G_1 (\BA_f) / \pi(K_1) ) &
                                [(x, \pi(p_1)] \\
P_1(\BQ) \backslash ( \FH_1 \times P_1 (\BA_f) / K_1 ) \ar[u] \ar[d] &
                                [(x, p_1)] \ar@{|->}[u] \ar@{|->}[d] \\
\forget{
\Stab_{Q (\BQ)} (\FH_1) \backslash
           ( \FH_1 \times \Stab_{Q (\BQ)} (\FH_1)P_1 (\BA_f) / K_1 ) \ar[d] &
                                [(x, p_1)] \ar@{|->}[d] \\
}
G (\BQ) \backslash ( \FH^* \times G (\BA_f) / K ) &
                                [(x, p_1g)] \\
M^K(G, \FH)^* (\BC) \ar@{=}[u] &
\\}}
\]
The map $[(x,p_1)] \mapsto [(x, \pi(p_1))]$ is an
isomorphism of complex spaces.
The composition of these morphisms comes from
a unique morphism of algebraic varieties
\[
i = i_{G_1,K,g}: M^{\pi(K_1)}:= M^{\pi(K_1)}(G_1,\FH_1) \longto (M^K)^*
\]
(\cite[Prop.~15]{Se}, applied to the Baily--Borel compactification of
$M^{\pi(K_1)}$; see \cite[7.6]{P1}).
Define the following groups (cmp.\ \cite[(3.7.4)]{P2}, where the same
notations are used for the groups
corresponding to the actions from the right):
\begin{eqnarray*}
H_Q & := & \Stab_{Q (\BQ)} (\FH_1) \cap P_1 (\BA_f) \cdot K' \; , \\
H_C & := & \Cent_{Q (\BQ)} (\FH_1) \cap W_1 (\BA_f) \cdot K'  \; ;
\end{eqnarray*}
note that these are indeed groups since $Q$ normalizes $P_1$ and $W_1$.
The group $H_Q$ acts by analytic automorphisms on
$\FH_1 \times P_1 (\BA_f) / K_1$ (see Lemma~\ref{1D} below for an
explicit description of this action). Hence the group
$\Delta_1 := H_Q / P_1(\BQ)$ acts naturally on
\[
M^{\pi(K_1)} (\BC) =
P_1(\BQ) \backslash ( \FH_1 \times P_1 (\BA_f) / K_1 ) \; .
\]
This action is one by algebraic automorphisms \cite[Prop.~9.24]{P1}.
By \cite[6.3]{P1} (see also Proposition~\ref{1A} below),
it factors through a finite quotient of $\Delta_1$, which we shall denote by
$\Delta$.
The quotient by this action is precisely the image of $i$:
\[
\vcenter{\xymatrix@R-10pt{
M^{\pi(K_1)} \ar@{->>}[r] \ar@/_2pc/[rr]_i&
          M_1^K := \Delta \backslash M^{\pi(K_1)} \ar@{^{ (}->}[r]^-i & (M^K)^*
\\}}
\]
By abuse of notation, we denote by the same letter $i$ the inclusion of
$M_1^K$ into $(M^K)^*$.
We need to identify the group $\Delta$, and the nature of its action
on $M^{\pi(K_1)}$. Let us introduce
the following condition on $(G, \FH)$:
\begin{enumerate}
\item [$(+)$] The neutral connected component $Z (G)^0$ of the center
$Z (G)$ of $G$
is, up to isogeny, a direct product of a $\BQ$-split torus with a torus $T$
of compact type (i.e., $T(\BR)$ is compact) defined over $\BQ$.
\end{enumerate}
The proof of \cite[Cor.~4.10]{P1} shows that
$(G_1, \FH_1)$ satisfies $(+)$ if $(G, \FH)$ does.

\begin{Prop} \label{1A}
(a) The subgroup $P_1(\BQ) H_C$ of $H_Q$ is of finite index.\\[0.2cm]
(b) The group $H_C / W_1 (\BQ)$ centralizes $G_1$,
and $H_C$ is the kernel of the
action of $H_Q$ on $\FH_1 \times G_1 (\BA_f) / \pi(K_1)$.
The group $P_1(\BQ) H_C$ acts trivially on $M^{\pi(K_1)}$.\\[0.2cm]
(c) Assume that $(G,\FH)$ satisfies $(+)$, and that $K$ is \emph{neat}
(see e.g.\ \cite[0.6]{P1}). Then we have an equality
\[
P_1(\BQ) \cap H_C = W_1(\BQ)
\]
of subgroups of $H_Q$.\\[0.2cm]
(d) Under the hypotheses of (c), the action of the finite group
$H_Q / P_1(\BQ) H_C$ on $M^{\pi(K_1)}$ is free. In particular, we have
$\Delta = H_Q / P_1(\BQ) H_C$.
\end{Prop}

For the proof of this result, we shall need three lemmata.
By slight abuse of notation, we denote by the letter $\pi$
the projection $Q \to Q/W_1$ as well:
\[
\vcenter{\xymatrix@R-10pt{
P_1 \ar@{^{ (}->}[r] \ar@{->>}[d]_{\pi} & Q \ar@{->>}[d]^{\pi} \\
G_1 \ar@{^{ (}->}[r] & Q/W_1
\\}}
\]
Note that since $Q/W_1$ is reductive, it is possible to choose a
complement of $G_1$ in $Q/W_1$, i.e., a normal connected reductive
subgroup $G_2$ of $Q/W_1$ such that
\[
Q/W_1 = G_1 \cdot G_2 \; ,
\]
and such that the intersection $G_1 \cap G_2$ is finite.
Let us mention that in the literature, the groups $G_1$ and $G_2$ are
sometimes called $G_h$ and $G_{\ell}$, respectively.

\begin{Lem} \label{1B}
Let $\gamma \in \Cent_{Q (\BQ)} (\FH_1)$. Then $\pi(\gamma) \in Q/W_1$
centralizes $G_1$.
\end{Lem}

\begin{Proof}
For any $g_1 \in G_1(\BQ)$, the element $g_1 \pi(\gamma) g_1^{-1}$ centralizes
$\FH_1$. It follows that $\pi(\gamma)$ centralizes
the normal subgroup of $G_1$ generated by the images of the morphisms
\[
h_x: \BS \longto G_{1,\BR} \quad , \quad x \in \FH_1 \; .
\]
But by definition of $P_1$ (see \cite[4.7]{P1}), this subgroup is $G_1$
itself.
\end{Proof}

\begin{Lem} \label{1C}
Let $(G', \FH')$ be Shimura data satisfying condition $(+)$, and
$\Gamma \subset G'(\BQ)$ an arithmetic subgroup (see
e.g.\ \cite[0.5]{P1}). Then $\Gamma$ acts
properly discontinuously on $\FH'$. In particular, the stabilizers
of the action of $\Gamma$ are finite.
\end{Lem}

\begin{Proof}
We cannot quote \cite[Prop.~3.3]{P1} directly because loc.\ cit.\ uses
a more general notion of properly discontinuous actions (see \cite[0.4]{P1}).
However, the proof of \cite[Prop.~3.3]{P1} shows that the action of $\Gamma$ is
properly discontinuous in the usual sense
if condition $(+)$ is satisfied. We refer to \cite[Prop.~1.2~b)]{W1} for the
details.
\end{Proof}

Let us identify explicitly the action of $H_Q$ on
$\FH_1 \times P_1 (\BA_f) / K_1$:

\begin{Lem} \label{1D}
Let $x \in \FH_1$, $p_1 \in P_1 (\BA_f)$, and
$\gamma \in H_Q \subset \Stab_{Q (\BQ)} (\FH_1)$. Write
\[
\gamma = p_2 k \; ,
\]
with $p_2 \in P_1 (\BA_f)$, and $k \in K'$. Since $Q$ normalizes $P_1$, we have
\[
p_3 := \gamma p_1 \gamma^{-1} \in P_1 (\BA_f) \; .
\]
We then have
\[
\gamma \cdot (x,[p_1]) = (\gamma (x), [p_3 p_2] )
                       = (\gamma (x), [\gamma p_1 k^{-1}])
\]
in $\FH_1 \times P_1 (\BA_f) / K_1$.
\end{Lem}

We leave the proof of this result to the reader.\\

\begin{Proofof}{Proposition~\ref{1A}}
As for (a), observe that by Lemma~\ref{1B},
the images of both $H_C$ and $H_Q$ in $Q/P_1 (\BQ)$
are arithmetic subgroups. Part (b) results directly from
Lemmata~\ref{1B} and \ref{1D}.
Let us turn to (c). By Lemma~\ref{1B}, the image of the
group $P_1(\BQ) \cap H_C$ under $\pi$ is
an arithmetic subgroup of the center of $G_1$. It is neat because $K$ is.
Because of $(+)$, it must be trivial.
It remains to show (d). Fix $x \in \FH_1$, $p_1 \in P_1 (\BA_f)$, and
$\gamma \in H_Q$ as well as $k \in K'$
as in Lemma~\ref{1D}. Suppose that
\[
\gamma \cdot [(x,p_1)] = [(x,p_1)]
\]
in $M^{\pi(K_1)}(\BC)$. There is thus an element $\gamma'$ of $P_1(\BQ)$,
such that
\[
\gamma \cdot (x,[p_1]) = \gamma' \cdot (x,[p_1])
\]
in $\FH_1 \times P_1 (\BA_f) / K_1$.
In other words, we can find $k_1 \in K_1$ such that
\begin{enumerate}
\item[(1)] $\gamma' (x) = \gamma (x)$, i.e.,
$\gamma'' := \gamma^{-1} \gamma' \in Q(\BQ)$
stabilizes $x$.
\item[(2)] $\gamma' p_1 = \gamma p_1 k^{-1} k_1$. We thus have
$\gamma'' \in p_1 K' p_1^{-1}$, which is a neat subgroup of
$G(\BA_f)$.
\end{enumerate}
Choose a complement $G_2$ of $G_1$ in $Q/W_1$.
The groups $G_1$
and $G_2$ centralism each other. Denote by $\Pi(\gamma'')$ the image of
$\gamma''$ in $Q/G_2$. We identify this group with
$G_1/G_1 \cap G_2$. Because of (1), the element $\Pi(\gamma'')$
stabilizes a point in the space $\FH_1 / G_1 \cap G_2$
belonging to the quotient Shimura data
\[
(Q/G_2, \FH_1 / G_1 \cap G_2) := (G_1,\FH_1) / G_1 \cap G_2 \; .
\]
By Lemma~\ref{1C}, the element $\Pi(\gamma'')$ is of finite order. Because
of (2), it must be trivial. We conclude:
\begin{enumerate}
\item[(3)] $\pi(\gamma'')$ centralizes $G_1$.
Hence $p_1 k^{-1} k_1 = \gamma'' p_1 =
p_1 \gamma'' \mod W_1 (\BA_f)$, and we get:
\item[(4)] $\gamma''$ lies in $W_1 (\BA_f) \cdot K'$.
\end{enumerate}
Since $G_1(\BR)$ acts transitively on $\FH_1$, (1) and (3) imply
that $\gamma''$ acts trivially on $\FH_1$.
Because of (4), we then have $\gamma'' \in H_C$, hence
\[
\gamma = \gamma' (\gamma'')^{-1} \in P_1(\BQ) H_C \; ,
\]
as claimed.
\end{Proofof}

For future reference, we note:

\begin{Cor} \label{1E}
Assume that $(G, \FH)$ satisfies $(+)$, and that $K$ is neat. Then the kernel
of the projection $\Delta_1 \to \Delta$ is canonically isomorphic to
$\pi (H_C)$.
\end{Cor}

\begin{Proof}
This follows immediately from Proposition~\ref{1A}~(c).
\end{Proof}


\bigskip
%
%

\section{Statement of the main result}
\label{2}

Let $(M^K)^* = M^K (G, \FH)^*$ be the Baily--Borel compactification of a
Shimura variety $M^K = M^K (G, \FH)$, and
$M^K_1 = \Mpd = \Delta \backslash M^{\pi(K_1)}(G_1, \FH_1)$
a boundary stratum.
Consider the situation
\[
M^K \stackrel{j}{\longinto} (M^K)^*
\stackrel{i}{\leftlonginto} \Mo \; .
\]
Saito's formalism \cite{Sa} gives a functor $i^{\ast} j_{\ast}$ between the
bounded derived categories of \emph{algebraic mixed Hodge modules} on
$M^K$ and on $\Mo$ respectively. Our main result (Theorem~\ref{2D})
gives a formula for the restriction of $i^{\ast} j_{\ast}$ onto the image of
the natural functor $\mu_K$ associating to an
algebraic representation of $G$ a
variation of Hodge structure on $M^K$. Its proof, which will rely on the
material developed in the next six sections, will be given in Section~\ref{5}.
In the present section, we shall restrict ourselves to a concise presentation
of the ingredients necessary for the formulation of Theorem~\ref{2D}
(\ref{2A}--\ref{2C}), and we shall state the main corollaries
(\ref{2E}--\ref{2J}). Let us mention that part of these results are already
contained in work of Harder, Looijenga--Rapoport, and Harris--Zucker
\cite{Hd,LR,HZ1}
(see Remark~\ref{2K}).
\ref{2D}--\ref{2J} are the Hodge theoretic analogues of
results obtained by Pink in the $\ell$-adic context \cite{P2}. \\

Fix pure Shimura data $(G, \FH)$
satisfying the hypothesis $(+)$, and an open compact neat subgroup $K$ of
$G(\BA_f)$. Let $F$ be a subfield of $\BR$.
By definition of Shimura data, there is a tensor functor
associating to an algebraic $F$-representation $\BV$ of $G$
a variation of Hodge structure $\mu (\BV)$ on $\FH$ \cite[1.18]{P1}.
It descends to a variation $\mu_K (\BV)$ on $M^K (\BC)$.
We refer to the tensor functor $\mu_K$ as the \emph{canonical construction} of
variations of Hodge structure from re\-pre\-sentations of $G$. Since
the weight cocharacter associated to $(G, \FH)$ is central
\cite[Def.~2.1~(iii)]{P1},
$\mu_K (\BV)$ is the direct sum of its weight graded objects.
By Schmid's Nilpotent Orbit Theorem \cite[Thm.~(4.9)]{Sch}, the
image of $\mu_K$ is contained in the category
$\VMHS_F M^K$ of \emph{admissible} variations, and hence
\cite[Thm.~3.27]{Sa}, in the category $\MHM_F M^K$ of algebraic mixed Hodge
modules. Since the functor $\mu_K$ is exact, it descends to the level of
derived categories:
\[
\mu_K: D^b (\Rep_F G) \longto D^b (\MHM_F M^K) \; .
\]
In order to state the main result, fix
a rational boundary component $(P_1, \FX_1)$ of $(G, \FH)$,
and an element $g \in G(\BA_f)$.
We shall use the notation of Proposition~\ref{1A}. In particular,
we have the following diagram of algebraic groups:
\[
\vcenter{\xymatrix@R-10pt{
P_1 \ar@{^{ (}->}[r] \ar@{->>}[d]_{\pi} &
                        Q \ar@{->>}[d]^{\pi} \ar@{^{ (}->}[r] & G \\
G_1 \ar@{^{ (}->}[r] & Q/W_1 &
\\}}
\]
By Proposition~\ref{1A}~(c), we have
a Cartesian diagram of subgroups of $Q(\BQ)$,
all of which are normal in $H_Q$:
\[
\vcenter{\xymatrix@R-10pt{
& P_1(\BQ) \ar@{^{ (}->}[dr] & \\
W_1(\BQ) \ar@{^{ (}->}[ur] \ar@{_{ (}->}[dr] && H_Q \\
& H_C \ar@{_{ (}->}[ur] &
\\}}
\]
Writing $\HQ$ for $\pi(H_Q)$, and $\HC$ for $\pi(H_C)$, we thus have
\[
G_1(\BQ) \cap \HC = \{ 1 \} \; .
\]

\begin{Def} \label{2A}
(a) The category $(\Rep_F G_1, \HQ)$ consists of pairs
\[
(\BV_1, (\rho_\gamma)_{\gamma \in \HQ}) \; ,
\]
where $\BV_1 \in \Rep_F G_1$, and
$(\rho_\gamma)_{\gamma \in \HQ}$ is a family of
isomorphisms
\[
\rho_\gamma: (\Int \gamma)^* \BV_1 \isoto \BV_1
\]
in $\Rep_F G_1$ ($\, \Int \gamma := $ conjugation by $\gamma$ on $G_1$)
such that the
following hold:
\begin{enumerate}
\item[(i)] $\rho_\gamma$ is given by $v \mapsto \gamma^{-1} (v)$ if
$\gamma \in G_1(\BQ)$,
\item[(ii)] the cocycle condition holds.
\end{enumerate}
Morphisms in $(\Rep_F G_1, \HQ)$ are defined in the obvious way. \\[0.2cm]
(b) The category $(\Rep_F G_1, \HQ/\HC)$ is defined as the full
sub-category of $(\Rep_F G_1, \HQ)$ consisting of objects
\[
(\BV_1, (\rho_\gamma)_{\gamma \in \HQ})
\]
for which $\rho_\gamma$ is the identity whenever $\gamma$ lies in $\HC$.
\end{Def}

We also define variants of the above on the level of pro-categories, i.e.,
categories
$(\Pro(\Rep_F G_1), \HQ)$ and $(\Pro(\Rep_F G_1), \HQ/\HC)$.
Note that by Proposition~\ref{1A}~(b),
we have $(\Int \gamma)^* \BV_1 = \BV_1$ for any
$\BV_1 \in \Pro(\Rep_F G_1)$ and $\gamma \in \HC$. \\

The functor $\mup$ extends to give
an exact tensor functor from the category $(\Rep_F G_1, \HQ/\HC)$ to
the category of objects of $\VMHS_F \Mp$ with an action of the finite group
$\Delta = \HQ/G_1(\BQ)\HC$. From Proposition~\ref{1A}~(d), we get a canonical
equivalence of this category and $\VMHS_F \Mo$. Altogether, we have defined
a tensor functor
\[
(\Rep_F G_1, \HQ/\HC) \longto \VMHS_F \Mo \subset \MHM_F \Mo \; ,
\]
equally referred to as $\mup$. It is exact, and hence defines
\[
\mup: D^b (\Rep_F G_1, \HQ/\HC) \longto D^b (\MHM_F \Mo) \; .
\]

\begin{Def} \label{2B}
Denote by
\[
\Gamma(\HC,\argdot):
(\Rep_F G_1, \HQ) \longto (\Rep_F G_1, \HQ/\HC)
\]
the left exact functor associating to $\BV_1 = (\BV_1, (\rho_\gamma)_\gamma)$
the largest sub-object $\BV_1'$ on which
the $\rho_\gamma$ act as the identity whenever $\gamma \in \HC$.
\end{Def}

Instead of $\Gamma(\HC, (\BV_1, (\rho_\gamma)_\gamma))$, we shall often write
$(\BV_1, (\rho_\gamma)_\gamma)^{\HC}$, or simply $\BV_1^{\HC}$.
Note that $\Gamma(\HC, \argdot)$ extends to a functor
\[
(\Pro(\Rep_F G_1), \HQ) \longto (\Pro(\Rep_F G_1), \HQ/\HC) \; .
\]
This functor will be studied
in Section~\ref{3}; for the time being, let us accept that the total right
derived functor of $\Gamma(\HC,\argdot)$ exists (Theorem~\ref{3T}~(a)):
\[
R\Gamma(\HC,\argdot): D^+ (\Pro(\Rep_F G_1), \HQ) \longto
                      D^+ (\Pro(\Rep_F G_1), \HQ/\HC) \; ,
\]
and that it respects the sub-categories $D^{b}(\Rep_F G_{1}, ? )$
(Theorem~\ref{3T}~(b)):
\[
R\Gamma(\HC,\argdot): D^b (\Rep_F G_1, \HQ) \longto
                      D^b (\Rep_F G_1, \HQ/\HC) \; .
\]
The cohomology functors associated to
$R\Gamma(\HC,\argdot)$ will be referred to
by $H^p(\HC,\argdot)$, for $p \in \BZ$. Let us assemble the properties of
these functors necessary for the understanding of our main result. For their
proof, we refer to Section~\ref{3}.

\begin{Prop} \label{2Ba}
(a) The vector space underlying $H^p(\HC,\argdot)$ is
given by the cohomology of the abstract group $\HC$. More precisely,
there is a commutative diagram of functors
\[
\vcenter{\xymatrix@R-10pt{
D^b (\Rep_F G_1, \HQ) \ar[d]_{R\Gamma(\HC,\argdot)} \ar[r] &
D^+ (\Rep_F \HQ) \ar[d]^{R\Gamma(\HC,\argdot)} \\
D^b (\Rep_F G_1, \HQ/\HC) \ar[r] & D^+ (\Rep_F \HQ/\HC)
\\}}
\]
Here, the categories at the right hand side denote the derived categories
of abstract representations, and the arrow $R\Gamma(\HC,\argdot)$
between them denotes the total derived functor of
the functor associating to an representation its $\HC$-invariants.
The horizontal arrows are the natural forgetful functors. \\[0.2cm]
(b) Let $\BV_1 \in (\Rep_F G_1, \HQ)$, and $p \in \BZ$. Consider the
algebraic representations $\Res^{\HQ}_{G_1} \BV_1$ and
$\Res^{\HQ/\HC}_{G_1} H^p(\HC,\BV_1)$ of $G_1$. Then any irreducible
factor of $\Res^{\HQ/\HC}_{G_1} H^p(\HC,\BV_1)$ is an irreducible
factor of $\Res^{\HQ}_{G_1} \BV_1$.
\end{Prop}

Observe that the weight cocharacter associated to the Shimura data
$(G_1, \FH_1)$ maps
to the center of $G_1$, and hence to the center of $Q/W_1$. It follows that
any object of $(\Rep_F G_1, \HQ)$ or of $(\Rep_F G_1, \HQ/\HC)$
is the direct sum of its
weight-graded objects. Proposition~\ref{2Ba} implies in particular:

\begin{Cor} \label{2Bb}
The functors
\[
H^p(\HC,\argdot): (\Rep_F G_1, \HQ) \longto (\Rep_F G_1, \HQ/\HC)
\]
respect the sub-categories of
pure objects. Hence they preserve the weight decompositions
in both categories.
\end{Cor}

\begin{Def} \label{2C}
Denote by
\[
\Gamma(W_1,\argdot):
\Rep_F Q \longto (\Rep_F G_1, \HQ)
\]
the left exact functor associating to a representation $\BX$ of $Q$ its
invariants under $W_1$.
\end{Def}

Instead of $\Gamma(W_1, \BX)$, we shall often write
$\BX^{W_1}$. The total right derived functor of
$\Gamma(W_1,\argdot)$ in the sense of
\cite[1.2]{D} exists (see e.g.\ \cite[Thm.~2.2, 2.3]{W3}):
\[
R\Gamma(W_1,\argdot): D^b (\Rep_F Q) \longto D^b (\Rep_F G_1, \HQ)  \; .
\]
In fact, this is the composition of the functor
\[
D^b (\Rep_F Q) \longto D^b (\Rep_F Q / W_1)
\]
calculating algebraic (Hochschild) cohomology of $W_1$,
and the forgetful functor from $D^b (\Rep_F Q / W_1)$
to $D^b (\Rep_F G_1, \HQ)$.
The cohomology functors associated to $R\Gamma(W_1,\argdot)$
will be referred to
by $H^q(W_1,\argdot)$, for $q \in \BZ$.\\

Denote by $c$ the codimension of $\Mo$ in $(M^K)^*$.
Our main result reads as follows:

\begin{Thm} \label{2D}
There is a natural commutative diagram
\[
\vcenter{\xymatrix@R-10pt{
D^b (\Rep_F G) \ar[rr]^-{\mu_K} \ar[dd]_{\Res^G_Q} &&
                                  D^b (\MHM_F M^K) \ar[dddddd]^{i^*j_*[-c]} \\
         &&  \\
D^b (\Rep_F Q) \ar[dd]_{R\Gamma(W_1,\argdot)} && \\
         &&  \\
D^b (\Rep_F G_1, \HQ) \ar[dd]_{R\Gamma(\HC,\argdot)} && \\
         &&  \\
D^b (\Rep_F G_1, \HQ/\HC) \ar[rr]^-\mup && D^b (\MHM_F \Mo)
\\}}
\]
In particular, the functor $i^*j_* \circ \mu_K$
takes values in the sub-category of objects of $D^b (\MHM_F \Mo)$ that can
be represented by complexes of direct sums of
pure variations of Hodge structure.
\end{Thm}

The proof of this result will be given in Section~\ref{5}.
Theorem~\ref{2D} is the Hodge theoretic analogue of the main result
of \cite{P2}. It expresses
$i^*j_* \circ \mu_K$ as a composition of two derived functors.

\begin{Cor} \label{2E}
For any $\BV^\bullet \in D^b (\Rep_F G)$, there is a canonical and
functorial spectral sequence
\[
E^{p,q} = \mup \circ H^p (\HC, H^q (W_1, \Res^G_Q \BV^\bullet)) \Longrightarrow
                       \CH^{p+q-c} i^*j_* \circ \mu_K (\BV^\bullet)
\]
in the category of admissible variations on $\Mo$.
\end{Cor}

The central ingredient for the analysis of this spectral sequence is the
following well known fact (see e.g.\ \cite[proof of 1.6.2, Satz 1]{Hd}):

\begin{Prop} \label{2F}
For any $\BX^\bullet \in D^b (\Rep_F Q)$, there is a canonical and functorial
isomorphism in $D^b (\Rep_F G_1, \HQ)$
\[
R\Gamma(W_1, \BX^\bullet) \isoto
                  \bigoplus_{q \in \BZ} H^q (W_1, \BX^\bullet)[-q] \; .
\]
\end{Prop}

\begin{Proof}
In fact, this isomorphism exists already on the level of
the derived category $D^b (\Rep_F Q/W_1)$.
Since $Q/W_1$ is reductive,
the category $\Rep_F Q/W_1$ is semi-simple. Its derived category
is therefore canonically
equivalent to the category of graded objects in $\Rep_F Q/W_1$.
\end{Proof}

Consequently, we have:

\begin{Thm} \label{2G}
The spectral sequence of Corollary~\ref{2E} degenerates and splits canonically.
Therefore, there is a canonical and
functorial isomorphism in $\VMHS_F \Mo$
for any $\BV^\bullet \in D^b (\Rep_F G)$
\[
\CH^{n} i^*j_* \circ \mu_K (\BV^\bullet) \isoto
\bigoplus_{p+q = n+c} \mup \circ H^p (\HC, H^q (W_1, \Res^G_Q \BV^\bullet)) \; ,
\]
for any $n \in \BZ$.
\end{Thm}

This is the Hodge theoretic analogue of \cite[Thm.~5.3.1]{P2}. The
isomorphism of local systems underlying the
isomorphism of Theorem~\ref{2G}
is known; see e.g.\ \cite[proof of 1.6.2, Satz 1]{Hd}.

\begin{Cor} \label{2H}
Let $\BV^\bullet \in D^b (\Rep_F G)$, and $n \in \BZ$. \\[0.2cm]
(a) The admissible variation
$\CH^{n} i^*j_* \circ \mu_K (\BV^\bullet)$ on $\Mo$
is the direct sum of its weight
graded objects. In particular, it is semi-simple. \\[0.2cm]
(b) For any $k \in \BZ$, we have a canonical and functorial isomorphism
in $\VMHS_F \Mo$
\[
\Gr^W_k \CH^{n} i^*j_* \circ \mu_K (\BV^\bullet) \isoto
\bigoplus_{p+q = n+c} \mup \circ \,
                  H^p (\HC, \Gr^W_k H^q (W_1, \Res^G_Q \BV^\bullet)) \; .
\]
\noindent (c) For any $k \in \BZ$, any Hodge type occurring in
$\Gr^W_k \CH^{n} i^*j_* \circ \mu_K (\BV^\bullet)$ occurs
already in one of the $\Gr^W_k H^q (W_1, \Res^G_Q \BV^\bullet)$,
$q \in \BZ$, $q \le n+c$.
\end{Cor}

\begin{Rem} \label{2I}
Observe that the weights and the Hodge types of the objects
$H^q (W_1, \Res^G_Q \BV^\bullet)$ are computed in the category $\Rep_F G_1$.
They thus depend on the restriction $\Res^G_{P_1} \BV^\bullet$ of $\BV^\bullet$
to the subgroup $P_1$ of $G$. We remind the reader that in general, the functor
$\Res^G_{P_1}$ preserves neither the weights nor the Hodge types
of a representation. However, it \emph{does} preserve the Hodge
\emph{filtration} \cite[Prop.~4.12]{P1}.
\end{Rem}

\begin{Proofof}{Corollary~\ref{2H}}
(a) follows from the existence of the weight decomposition in the
category $(\Rep_F G_1, \HQ/\HC)$.
The only point that needs to be explained in (b) is the equality
\[
\Gr^W_k H^p (\HC, H^q (W_1, \Res^G_Q \BV^\bullet)) =
H^p (\HC, \Gr^W_k H^q (W_1, \Res^G_Q \BV^\bullet)) \; .
\]
But this follows from Corollary~\ref{2Bb}. (c) results from (b), and from
Proposition~\ref{2Ba}~(b).
\end{Proofof}

Passage to singular cohomology yields the following:

\begin{Cor} \label{2J}
Denote by $a_1$ the structure morphism of $M_1^K$,
and by $\MHS_F$ the category $\MHM_F (\Spec \BC)$, that is,
the category of mixed graded-polarizable
$F$-Hodge structures.\\[0.2cm]
(a) There is a natural commutative diagram
\[
\vcenter{\xymatrix@R-10pt{
D^b (\Rep_F G) \ar[rr]^-{\mu_K} \ar[dd]_{\Res^G_Q} &&
                                  D^b (\MHM_F M^K) \ar[dddddd]^{i^*j_*[-c]} \\
         &&  \\
D^b (\Rep_F Q) \ar[dd]_{R\Gamma(W_1,\argdot)} && \\
         &&  \\
D^b (\Rep_F G_1, \HQ) \ar[dd]_{R\Gamma(\HC,\argdot)} && \\
         &&  \\
D^b (\Rep_F G_1, \HQ/\HC) \ar[dd]_\mup && D^b (\MHM_F \Mo) \ar[dd]^{{a_1}_*} \\
         &&  \\
D^b (\MHM_F \Mo) \ar[rr]^-{{a_1}_*} && D^b (\MHS_F)
\\}}
\]
(b) For any $\BV^\bullet \in D^b (\Rep_F G)$, there is a canonical and
functorial spectral sequence
\[
E^{p,q} = H^p (\Mo, \HC, H^q (W_1, \Res^G_Q \BV^\bullet)) \Longrightarrow
          H^{p+q} (\Mo, i^* j_* \circ \mu_K (\BV^\bullet) )
\]
in the category of Hodge structures. Here, $H^p (\Mo, \HC, \argdot)$
denotes the cohomology functors associated to the composition
\[
{a_1}_* \circ \mup \circ R\Gamma(\HC,\argdot): D^b (\Rep_F G_1, \HQ) \longto
D^b (\MHS_F) \; .
\]
(c) The spectral sequence of (b) degenerates and splits canonically.
\end{Cor}

\begin{Proof}
(a) follows from Theorem~\ref{2D},
(b) is clear, and (c) follows from
Proposition~\ref{2F}.
\end{Proof}

\begin{Rem} \label{2K}
(a) Part of the information contained in Corollary~\ref{2H} is known.
The splitting of the weight filtration on
the variation $\CH^{n} i^*j_* \circ \mu_K (\BV^\bullet)$ is proved in
\cite[Prop.~(6.4)]{LR}. According to Corollary~\ref{2H}~(b), the
local system underlying $\Gr^W_k \CH^{n} i^*j_* \circ \mu_K (\BV^\bullet)$
is the direct sum of the local systems underlying
the $\mup \circ \, H^p (\HC, \Gr^W_k H^q (W_1, \Res^G_Q \BV^\bullet))$,
for $p+q = n+c$.
This result can also be obtained by combining
Proposition~\ref{2F} and \cite[Cor.~(6.6)]{LR}.\\[0.2cm]
(b) Note that we can do better than Corollary~\ref{2H}~(c). Let $\BV$
be in $\Rep_F G$. Kostant's theorem
\cite[3.2.3]{Vo} allows us to identify the $Q/W_1$-representations
$H^q (W_1, \Res^G_Q \BV)$, for $\BV \in \Rep_F G$. Combining this with
Corollary~\ref{2H}~(b) and Proposition~\ref{2Ba}~(b), we get an explicit
a priori list of possible constituents of
$\Gr^W_k \CH^{n} i^*j_* \circ \mu_K (\BV)$.\\[0.2cm]
(c) The cohomology
$H^n (\Mo, i^* j_* \circ \mu_K (\BV^\bullet) )$
coincides with what is called \emph{deleted neighbourhood cohomology}
(of certain strata in a toroidal compactification of $M^K$) in
\cite{HZ1}, as can be seen from proper base change for the morphism from
a toroidal to the Baily--Borel compactification.
Thus, Corollary~\ref{2J} is equivalent to
\cite[Theorem~(5.6.10)]{HZ1}
for maximal parabolic subgroups ($R = P$ in the notation of \loccit).
\end{Rem}

\begin{Over} \label{2L}
Let us finish this section by an overview of the
proof of Theorem~\ref{2D}, and of the remaining sections of this paper.
We start by developing the basics of abstract group cohomology in
Abelian categories (Section~\ref{3}). We shall see in particular
(Proposition~\ref{3N}) that as in the case of usual group cohomology,
it can be calculated using free resolutions of the trivial module $\BZ$.
In Section~\ref{3b}, we state basic results on equivariant Hodge modules.\\
For simplicity, let us assume that the finite group $\Delta$ is trivial,
and hence, that $M^K_1$ is an actual Shimura variety.
The proof of Theorem~\ref{2D} involves the analysis of the degeneration
in a toroidal compactification
\[
j_\FS: M^K \longinto \MKS
\]
of $M^K$, and
the computation of its direct image under the proper morphism $p$ from
$\MKS$ to $(M^K)^*$. \\
The pre-image $i_\FS: \MoS \into \MKS$
of $M^K_1$ is itself stratified.
We recall the description of this stratification
in Section~\ref{4}, following the
presentation of \cite[(3.9)]{P2}. We recall in particular (Proposition~\ref{4D})
that the formal
completion of $\MKS$ along $\MoS$ is isomorphic to
the quotient by the action of $\Delta_1$
of the formal completion of a certain torus embedding $\MKoS$
along a union $Z$ of strata.
The action of $\Delta_1$ is free, and so is the induced action on the
set $\FT$ indexing the strata of $Z$.
In particular, each individual stratum of $\MKS$ can be seen as
a stratum of the torus embedding.
We are able to identify the composition
\[
p_* \is_* i^*_{\sigma} \, {j_\FS}_* \circ \mu_K
\]
for any stratum $\is: \Ms \into Z$ of $Z$, by
appealing to known results, which we recall in
Section~\ref{3a}, on degeneration along strata,
and on direct images to the base of torus compactifications: the value of
$p_* \is_* i^*_{\sigma} \, {j_\FS}_*$ on $\mu_K (\BV^\bullet)$ is equal to
$\mup \circ R\Gamma(W_1,\Res^G_{P_1} \BV^\bullet)$. In particular, it
does not depend on the stratum $\Ms$.
Since our isomorphisms are well behaved
under the action of $\Delta_1$, we get an object $\BX^\bullet$ in the
category $\Delta_1 \text{-} D^b ((\MHM_F M^K_1)^{\FT})$, i.e., a class
in $D^b ((\MHM_F M^K_1)^{\FT})$ with an action of $\Delta_1$.
It is constant in the sense that its components, indexed by $\FT$, are
all isomorphic. \\
The problem is now to put the information together, in order to compute
\[
p_* \, i_\FS^* \, {j_\FS}_* \circ \mu_K \; .
\]
The formal setting for this is provided by the theory of
\v{C}ech complexes for Hodge modules, the basics of which are contained
in Section~\ref{3aa}. In particular (Corollary~\ref{5R}), we show how
to compute the direct image $p_*$ with the help of stratifications.
The answer we get is perfectly adapted to the formalism
of group cohomology:
assume that $\BX^\bullet$ is concentrated in a single
degree. It is given by the same
object for all strata of $\FT$. Since the action of $\Delta_1$ is
free on $\FT$, Proposition~\ref{3N} and Corollary~\ref{5R} tell us that
the functor $p_* \, i_\FS^* \, {j_\FS}_*$
maps $\mu_K (\BV^\bullet)$ to $R\Gamma(\Delta_1,\BX^\bullet)$.\\
Unfortunately, $\BX^\bullet$ is always concentrated in more than one degree,
unless $\BV^\bullet = 0$.
The formalism of \v{C}ech complexes does not in general
allow to recover $p_* \, i_\FS^* \, {j_\FS}_*$ just from the collection
of the $p_* \is_* i^*_{\sigma} \, {j_\FS}_*$,
viewed as functors on derived
categories; we need to relate complexes on different strata \emph{before}
passing to the derived category. In other words, the information provided
by the object
\[
\BX^\bullet \in \Delta_1 \text{-} D^b ((\MHM_F M^K_1)^{\FT})
\]
is too weak; what is needed is an object in
\[
D^b (\Delta_1 \text{-} (\MHM_F M^K_1)^{\FT}) \; .
\]
This explains the presence of Section~\ref{3c}, which
provides the missing global information on the degeneration along $Z$.
We work on the normal cone of $\MKoS$
along $Z$, and identify the value of Saito's \emph{specialization functor}
$Sp_Z$ on ${j_\FS}_* \circ \mu_K (\BV^\bullet)$. The most difficult ingredient
is control of the \emph{monodromy weight filtration}. We
recall a number of results from the literature:
first, the explicit description, due to
Galligo--Granger--Maisonobe, of the category of perverse sheaves \emph{of
normal crossing type} on a product of unit disks; then,
Saito's identification of
the specialization functor in this description;
next, the fundamental theorems,
due to Schmid and Cattani--Kaplan, on \emph{nilpotent orbits};
finally, Kashiwara's permanence result on nilpotent orbits under
the \emph{nearby cycle functor}. It then suffices
to combine all these results in order to deduce the desired statement on
the monodromy weight filtration of the composition
$Sp_Z \, {j_\FS}_* \circ \mu_K (\BV^\bullet)$
(Theorem~\ref{3cY}, Corollary~\ref{3cZ}). \\
Section~\ref{5} puts everything together, and concludes the proof
of Theorem~\ref{2D}.
\end{Over}


\bigskip
%
%

\section{On the formalism of group cohomology}
\label{3}

Let $\mathcal{A}$ be an Abelian category, and $H$ an
abstract group. We shall denote by $H \text{-}\mathcal{A}$ the category
of objects of $\mathcal{A}$ provided with a left $H$-action and by
$\Pro(\mathcal{A})$ the pro-category associated to
$\mathcal{A}$ (see \cite[0.5]{D}).
Hence $H \text{-}\Pro(\mathcal{A})$ is the category of
pro-objects
of $\mathcal{A}$ provided with a left $H$-action. All these
categories are also Abelian. In this section all functors will be
additive. If $\gamma \in H$ and
$A\in \Ob(H \text{-}\mathcal{A}) $, we denote by the same letter
$\gamma$ the corresponding automorphism of $A$. We denote by $e$
the unit element of $H$.

\begin{Def} \label{3A}
  The \emph{fixed point functor} associated to $H$
  is the functor
  \begin{displaymath}
   \Gamma (H,\argdot) = (\argdot)^{H}:
   H \text{-}\Pro(\mathcal{A}) \longrightarrow \Pro(\mathcal{A})
  \end{displaymath}
  given by
  \begin{displaymath}
   \Gamma (H,A) = (A)^{H} := \bigcap_{\gamma \in H}\kernel(e-\gamma) \; .
  \end{displaymath}
\end{Def}

In general, the image of the category $H \text{-}\mathcal{A}$ under
the functor $\Gamma (H,\argdot)$ is not contained in $\mathcal{A}$,
unless certain conditions on $H$ or $\mathcal{A}$ are satisfied. Examples
for such conditions are:
$H$ is finitely generated, or $\mathcal{A}$ contains arbitrary products,
or $\mathcal{A}$ is Artinian.\\

The main aim of this section
is to show the existence of the right derived functor
\begin{displaymath}
  R\Gamma (H,\argdot) :D^{+}(H \text{-}\Pro(\mathcal{A}))
  \longrightarrow
  D^{+}(\Pro(\mathcal{A}))
\end{displaymath}
(Theorem~\ref{3L}).
We shall also show that, under some finiteness conditions on $H$,
the above derived functor can be lifted  to define
functors
\begin{displaymath}
  R\Gamma(H,\argdot) :D^{+}(H \text{-}\mathcal{A})\longrightarrow
  D^{+}(\mathcal{A})
\end{displaymath}
or
\begin{displaymath}
  R\Gamma(H,\argdot) :D^{b}(H \text{-}\mathcal{A})\longrightarrow
  D^{b}(\mathcal{A})
\end{displaymath}
(Theorem~\ref{3O}). As an application of these abstract principles,
we establish in Theorem~\ref{3T} the existence of the
functor $R\Gamma (\HC,\argdot)$, which occurs in the statement of our
main result, Theorem~\ref{2D}. We end the section by giving a proof of
Proposition~\ref{2Ba}. \\

The strategy for the construction of the derived functor is an
abstract version of a well known theme. The main interest of this
approach is that one does not need to suppose
the existence of sufficiently many
injective objects. We shall only treat the case of
covariant left exact functors, the other cases being completely
analogous.

\begin{Def} \label{3B}
  Let $\mathcal{A}$ be an Abelian category and let $\Id$ be the
  identity functor. A \emph{resolution functor} is an exact functor
  $C:\mathcal{A} \to \mathcal{A}$ provided with a morphism
  of functors $\Id \to C$ such that, for every $A\in
  \Ob(\mathcal{A})$ the map $A \to C(A)$ is a monomorphism.
\end{Def}

\begin{Def} \label{3C}
  Let $\mathcal{A}$ be an Abelian category, and $C$ a resolution
  functor. For any object $A$ of $\mathcal{A}$, the \emph{$C$-resolution} of
  $A$, denoted $C^{\ast}(A)$, is defined inductively as follows:
  \begin{alignat*}{2}
    K_{C}^{0}(A)&:=A \; ,&&\\
    C^{i}(A)&:=C(K_{C}^{i}(A)) \; ,&i&\ge 0 \; ,\\
    K_{C}^{i+1}(A)&:=\coker(K_{C}^{i}(A)\longinto
    C^{i}(A))\; , \quad&i&\ge 0 \; .
  \end{alignat*}
  The differential $d:C^{i} \to C^{i+1}$ is defined as the
  composition
  \begin{displaymath}
    C^{i}(A)\longonto K_{C}^{i+1}(A)\longinto C^{i+1}(A) \; .
  \end{displaymath}
\end{Def}

By definition, the sequence
\begin{displaymath}
  0\longrightarrow A\longrightarrow C^{0}(A)\longrightarrow C^{1}(A)
  \longrightarrow \ldots
\end{displaymath}
is exact.

\begin{Prop} \label{3D}
  Let $\mathcal{A}$ and $\mathcal{B}$ be Abelian categories. Let
  $F:\mathcal{A} \to \mathcal{B}$ be a left exact covariant
  additive functor. Let $C:\mathcal{A} \to \mathcal{A}$ be
  a resolution functor such that the composition $F\circ C$ is
  exact. Then the functor
  \begin{displaymath}
    F\circ C^{*}:K^{+}(\mathcal{A})\longrightarrow K^{+}(\mathcal{B})
  \end{displaymath}
descends to the level of derived categories. The resulting functor
  \begin{displaymath}
    D^{+}(\mathcal{A})\longrightarrow D^{+}(\mathcal{B}) \; ,
  \end{displaymath}
equally denoted by $F\circ C^{*}$,
is the total right derived functor $RF$ of $F$ in the sense of \cite[1.2]{D}.
\end{Prop}

\begin{Proof}
We need to show that the functor
  \begin{displaymath}
    F\circ C^{*}:K^{+}(\mathcal{A})\longrightarrow K^{+}(\mathcal{B})
    \longto D^{+}(\mathcal{B})
  \end{displaymath}
transforms quasi-isomorphisms into isomorphisms; from the construction of
the total derived functor \cite[1.2]{D}, it is clear that
this will imply the desired equality
$RF = F\circ C^{*}$. Using the cone of such a quasi-isomorphism, we are reduced
to showing that $F \circ C^{*}(K^\bullet)$ is acyclic for any acyclic complex
$K^\bullet$ in $C^+(\mathcal{A})$.
For this, it is enough
to show that for each $i\ge 0$, the functor $F\circ C^{i}$ is exact,
because in this case  $F\circ C^{\ast}(K^{\bullet})$ is the simple complex
associated to a double
complex with exact rows, hence acyclic.

\forget{To this end, we shall show: (1) any object $C(A)$ in the image of $C$ is
$F\circ C^{*}$-acyclic in the sense that $F\circ C^{*}$
transforms this object into a complex whose cohomology
is concentrated in degree $0$;
(2) any acyclic complex in
$C^+(\mathcal{A})$, whose components are $F\circ C^{*}$-acyclic, is transformed
into an acyclic complex under $F$.
}
The following lemma follows by induction from the exactness of the
functor $C$.

\begin{Lem} \label{3E}
The functors $K_{C}^i$ and $C^i$ are exact for all $i\ge 0$.
\forget{

    There are canonical isomorphisms
    \begin{displaymath}
      C^{i}(C(A))=C(C^{i}(A)) \; ,\quad K_{C}^{i}(C(A))=C(K_{C}^{i}(A))
    \end{displaymath}
for any object $A$ of $\mathcal{A}$.
}
\end{Lem}

In the situation of Proposition~\ref{3D},
we thus see that $F\circ C^i=(F\circ C)\circ K_{C}^i $ is the composition of two exact
functors, hence exact.
\forget{
Hence, by exactness of $F\circ C$,
  \begin{displaymath}
    H^{i}(F\circ C^{\ast}(C(A)))=H^{i}(F\circ C(C^{\ast}(A)))=0 \; , \quad
    i\ge 1 \; .
  \end{displaymath}
This proves claim (1). For (2), first use the fact that
$F\circ C^{\ast}$ is a $\delta$-functor in order to see that the quotient of
two $F\circ C^{*}$-acyclic objects is again $F\circ C^{*}$-acyclic.
Then observe that by left exactness of $F$, for any object $A$,
  \begin{displaymath}
    H^{0}(F\circ C^{\ast}(A))=F(A) \; .
  \end{displaymath}
Finally, one shows by induction on $i$ that for any complex $K^\bullet$
satisfying the hypotheses of (2), one has
\[
H^i (F(K^\bullet)) = 0 \; , \quad \text{for all} \; i \in \BZ \; .
\]
}
\end{Proof}

Recall \cite[p.~23]{D} that the derived functor $RF$ in the
sense of \cite[1.2]{D} satisfies the universal property of \cite[II.2.1.2]{V1}.

\begin{Ex} \label{3F}
  We recall how group homology can be defined using (the dual of) the above
  method. Let $H$ be a group, We denote by $\Ab$ the
  category of Abelian groups. Then $H \text{-}\Ab $ is the
  category of left $\mathbb{Z}H$-modules. If $A$ is an object of
  $H \text{-}\Ab$, then the group of co-invariants is
  \begin{displaymath}
    A_{H}=\mathbb{Z}\underset{\mathbb{Z}H}{\otimes} A
  \end{displaymath}
  where $\mathbb{Z}$ has the trivial $\mathbb{Z}H$ action. This
  defines a right exact functor $H \text{-}\Ab $ to $\Ab$ that we
  want to derive. Let $\For:H \text{-}\Ab \to \Ab $
  be the forgetful functor. We define the functor
  $C^H:H \text{-}\Ab \to H \text{-}\Ab $ by
  \begin{displaymath}
    C^H(A) := \mathbb{Z}H \underset{\mathbb{Z}}{\otimes}\For(A)=
    \Ind^{H}_{\{1\}}\Res^{H}_{\{1\}}(A) \; .
  \end{displaymath}
  The functor $C^H$ is exact and is equipped with a natural equivariant
  epimorphism
  $\epsilon : C^H(A) \to A$ given by
  \begin{displaymath}
    \epsilon (\sum n_{i}g_{i}\otimes a_{i})=
    \sum n_{i} g_{i} a_{i} \; .
  \end{displaymath}
  Applying the dual of the
  above construction we get a
  resolution
  \begin{displaymath}
    \ldots \longrightarrow C^H_{2}(A)\longrightarrow
    C^H_{1}(A)\longrightarrow C^H_{0}(A)
    \longrightarrow A\longrightarrow 0 \; .
  \end{displaymath}
  Since the composition
  $(\argdot)_{H}\circ C^H$ is the forgetful functor,
  the total left derived functor
  of the co-invariant functor is given by
  $(C^H_{\ast}(A))_{H}$.
  Note that $C^H_{\ast}(\mathbb{Z})$ is a free resolution of
  $\mathbb{Z}$.
\end{Ex}

Next we use the general theory to define group cohomology in an
arbitrary Abelian category.
Let $\mathcal{A}$, $H \text{-} \mathcal{A}$,
$\Pro(\mathcal{A})$ and $H \text{-}\Pro(\mathcal{A})$
be as in the beginning of the section.  We denote by $\For$ the
forgetful functor from $H \text{-}\mathcal{A}$ to $\mathcal{A}$,
as well as the forgetful functor between the pro-categories.

\begin{Def} \label{3G}
The \emph{resolution functor associated to} $H$ is
  the functor
  $C_{H}$ defined as follows.  Given an object $A$ of
$H \text{-}\Pro(\mathcal{A})$, the underlying object of
  $C_{H}(A)$ is $\prod_{h \in H} A\in
  \Ob(\Pro(\mathcal{A}) )$. Let
$p_{h}:C_{H}(A) \to A$ be the projection over the factor $h$. The
  action of an element $\gamma \in H$ over $C_{H}(A)$ is
  defined by the family of morphisms
  $$\gamma_{h}:C_{H}(A)\longrightarrow A \; ,\quad h \in H \; ,$$
  where $\gamma_{h}=\gamma \circ p_{\gamma^{-1} h}$.
\end{Def}

By definition, we have
\begin{displaymath}
  C_{H}(A)=\Hom_{\mathbb{Z}}(\mathbb{Z}H,A) \; ,
\end{displaymath}
with the diagonal action of $H$.\\

The following result is immediate from the definition of the action of
$H$ on $C_H(A)$. It is the
basic ingredient to define morphisms to $C_{H}(A)$.

\begin{Lem} \label{3H}
  Let $A$ and $B$ be objects of
$H \text{-}\Pro(\mathcal{A})$. Then there are canonical
  bijections between (a)~the set of morphisms $f:B \to C_{H}(A)$
  in the ca\-te\-gory $H \text{-}\Pro(\mathcal{A})$,
  (b)~the set of families of morphisms
$f_{h}:\For(B) \to \For(A)$ in $\Pro(\mathcal{A})$, $h\in H$,
  such that
  \begin{displaymath}
    f_{h}\circ \gamma =\gamma \circ f_{\gamma^{-1}h} \; , \quad
    \gamma \in H \; ,
  \end{displaymath}
and (c)~the set of morphisms
$f_{e}:\For(B) \to \For(A)$ in $\Pro(\mathcal{A})$.
In other words, the functor $C_H$ represents the functor on
$H \text{-}\Pro(\mathcal{A})$ given by
\[
B \longmapsto Hom_{\Pro(\mathcal{A})} (\For(B) , \For(A)) \; .
\]
\end{Lem}

\begin{Def} \label{3I}
  Denote by $\iota :A \to C_{H}(A)$ the
  canonical equivariant monomorphism
  determined by the family of morphisms $\iota _{h}=\Id$. We
  denote by $\tau :\For(A) \to \For(C_{H} (A))$ the
  morphism determined by the family of morphisms $\tau_{h}=h$.
\end{Def}

The following result follows easily from the definitions.

\begin{Prop} \label{3K}
  (a) The functor $C_{H}$ together with the morphism of functors
  $\iota $ is a resolution functor. \\[0.2cm]
  (b) The morphism $\tau $ induces an isomorphism of functors
    between $\For$ and $(\argdot)^{H}\circ C_{H}$.
\end{Prop}

Since $C_{H}$ is a resolution functor and
$(\argdot)^{H}\circ C_H$ is exact (because the functor $\For$ is exact),
we obtain the
following result:

\begin{Thm} \label{3L}
  The functor $(\argdot)^{H}:H \text{-}\Pro(\mathcal{A})
  \to \Pro(\mathcal{A})$ is right derivable, and the total
  right derived functor
  \begin{displaymath}
    R\Gamma(H,\argdot)
    :D^{+}(H\text{-}\Pro(\mathcal{A}))\longrightarrow
    D^{+}(\Pro(\mathcal{A}))
  \end{displaymath}
  is the functor induced by the exact functor $(\argdot)^{H}\circ
  C_{H}^{\ast}$.
\end{Thm}

The cohomology functors associated to
$R\Gamma(H,\argdot)$ will be denoted
by $H^p(H,\argdot)$, for $p \in \BZ$.
We shall see that the functor $R\Gamma(H,\argdot)$ can be computed using
\emph{any} right resolution of the trivial $H$-module $\mathbb{Z}$ by
free $\mathbb{Z}H$-modules. For any Abelian category
$\mathcal{A}$, we define a bifunctor
\begin{displaymath}
  \Hom:\Ab\times \Pro(\mathcal{A})\longrightarrow
  \Pro(\mathcal{A}) \; :
\end{displaymath}
Let $M \in \Ob(\Ab)$ and $A \in \Ob(\Pro(\mathcal{A}))$, and consider
the following contravariant functor $\CF$ on $\Pro(\mathcal{A})$:
by definition, $\CF (B)$ is the set of all group homomorphisms
\[
\alpha: M \longto Hom_{\Pro(\mathcal{A})} (B,A) \; .
\]
To see that this functor is representable, we first treat the case of
a free Abelian group $M$. Choose a basis
$\{x_{i}\}_{i\in I}$ of $M$. For any object $A$
of $\Pro(\mathcal{A})$, we see that
\begin{displaymath}
  \prod_{i\in I} A
\end{displaymath}
represents the functor $\CF$.
If $M$ is any Abelian group, we choose any two step free resolution
$$F_{2}\longrightarrow F_{1}\longrightarrow M\longrightarrow 0 \; .$$
We then have
\begin{displaymath}
  \Hom(M,A)=\kernel(\Hom(F_{1},A)\longrightarrow \Hom(F_{2},A)) \; .
\end{displaymath}

We can take into account the action of $H$:

\begin{Def} \label{3M}
  Let $A$ be an object of
  $H \text{-}\Pro(\mathcal{A})$, and let $M$ be a
  $\mathbb{Z}H$-module. The \emph{diagonal action} of $H$ over
  $\Hom(M,A)$ is defined as follows: for $B$ in $\Pro(\mathcal{A})$
  and $\alpha: M \to Hom_{\Pro(\mathcal{A})} (B,A)$ in $\Hom(M,A) (B)$,
  define $\gamma \alpha: M \to Hom_{\Pro(\mathcal{A})} (B,A)$ as
\[
m \longmapsto \gamma \circ \alpha (\gamma^{-1}(m)) \; .
\]
\end{Def}

When $M$ is a free $\mathbb{Z}H$-module, we can give an explicit description
of this action:
we choose a basis $\{x_{i}\}_{i\in I}$ of $M$ as
$\mathbb{Z}H$-module. Then
\begin{displaymath}
  \Hom(M,A)=\prod_{i}\prod_{h \in H}A \; .
\end{displaymath}
We write $p_{i,h}$ for the projection over the factor $(i,h)$.
Then the action of an element $\gamma \in H$ is determined by
the family of morphisms
\begin{displaymath}
  \gamma_{i,h}:\Hom(M,A)\longrightarrow A \; ,
\end{displaymath}
with $\gamma_{i,h}=\gamma \circ p_{i,\gamma^{-1}h}$.

\begin{Prop} \label{3N}
   (a) There are canonical equivalences of functors between $C_{H}$ and
    $\Hom(C^H(\mathbb{Z}),\argdot)$, and between $C_{H}^{\ast}$ and
    $\Hom(C^H_{\ast}(\mathbb{Z}),\argdot)$.\\[0.2cm]
    (b) Let $F_{\ast} \to \mathbb{Z}$ be any resolution of
    the trivial $\mathbb{Z}H$-module $\mathbb{Z}$ by free
    $\mathbb{Z}H$-modules. Then the functor $R\Gamma (H,\argdot)$ is
    induced by $(\Hom(F_{\ast},\argdot))^{H}$.
\end{Prop}

\begin{Proof}
  The fact that $C_{H}=\Hom(C^H(\mathbb{Z}),\argdot)$
  is a direct consequence of the
  definitions. Since the sequence
  \begin{displaymath}
    0\longrightarrow K_{1}(\mathbb{Z})\longrightarrow
    C^H_{0}(\mathbb{Z})\longrightarrow \mathbb{Z}\longrightarrow 0
  \end{displaymath}
  splits as a sequence of Abelian groups we obtain that the sequence
  \begin{displaymath}
    0\longrightarrow \Hom(K_{1}(\mathbb{Z}),A)\longrightarrow
    \Hom(C^H_{0}(\mathbb{Z}),A)\longrightarrow
    \Hom(\mathbb{Z},A)\longrightarrow 0
  \end{displaymath}
  is exact in $\Pro(\mathcal{A})$. Since $\Hom(\mathbb{Z},A)=A$
  we can prove by induction that $K_{i}(\mathbb{Z})$ is projective as
  Abelian group, that $\Hom(K_{i}(\mathbb{Z}),A)=K^{i}(A)$ and that
  $\Hom(C^H_{i}(\mathbb{Z}),A)=C_{H}^{i}(A)$. This proves part (a).

  For (b), use the fact that any
  $\mathbb{Z}H$-free resolution $F_\ast$
  of $\mathbb{Z}$ is homotopically equivalent to
  $C^H_{\ast}(\mathbb{Z})$. Therefore the complex
  $\Hom(F_{\ast},A)^{H}$ is homotopically equivalent to
  $(C^{\ast}_{H}(A))^{H}$.
\end{Proof}

Next we shall put some finiteness conditions on the group $H$.
Recall that a group is of type $FL$ if the trivial
$\mathbb{Z}H$-module $\mathbb{Z}$ admits a finite resolution
\begin{displaymath}
  0\longrightarrow F_{n}\longrightarrow \ldots \longrightarrow
  F_{1}\longrightarrow F_{0}\longrightarrow  \mathbb{Z}\longrightarrow 0
\end{displaymath}
by finitely generated free $\mathbb{Z}H$-modules.
A group is called $FP_{\infty}$ if $\mathbb{Z}$ admits a resolution by
finitely generated free $\mathbb{Z}H$-modules.

\begin{Thm} \label{3O}
    (a) If the group $H$ is of type $FP_{\infty}$, then
    there exists a canonical functor, also denoted by
    $$R\Gamma(H,\argdot)
    :D^{+}(H \text{-}\mathcal{A} ) \longrightarrow
    D^{+}(\mathcal{A}) \; ,$$
    and a natural commutative diagram
\[
\vcenter{\xymatrix@R-10pt{
        D^{+}(H \text{-}\mathcal{A} )
                                   \ar[rr]^-{R\Gamma(H,\argdot)} \ar[d] &&
        D^{+}(\mathcal{A}) \ar[d] \\
        D^{+}(H \text{-}\Pro(\mathcal{A}) )
                                   \ar[rr]^-{R\Gamma(H,\argdot)} &&
        D^{+}(\Pro(\mathcal{A}))
\\}}
\]
\noindent (b) If the group $H$ is of type $FL$, then the functor in (a)
    respects the bounded derived categories. We thus get
    a canonical functor, still denoted by
    $$R\Gamma(H,\argdot)
    :D^{b}(H \text{-}\mathcal{A} ) \longrightarrow
    D^{b}(\mathcal{A}) \; ,$$
    and a natural commutative diagram
\[
\vcenter{\xymatrix@R-10pt{
        D^{b}(H \text{-}\mathcal{A} )
                              \ar[rr]^-{R\Gamma(H,\argdot)} \ar[d] &&
        D^{b}(\mathcal{A}) \ar[d] \\
        D^{+}(H \text{-}\Pro(\mathcal{A}) ) \ar[rr]^-{R\Gamma(H,\argdot)} &&
        D^{+}(\Pro(\mathcal{A}))
\\}}
\]
\end{Thm}

\begin{Proof}
  If the group is of type $FP_\infty$, then there exists a resolution
  \begin{displaymath}
    \ldots \overset{f_{n+1}}{\longrightarrow} F_{n}
           \overset{f_{n}}{\longrightarrow}
    \ldots \overset{f_{2}}{\longrightarrow} F_{1}
           \overset{f_{1}}{\longrightarrow} F_{0}
           \longrightarrow \mathbb{Z} \longrightarrow 0 \; ,
  \end{displaymath}
  where every $F_{i}$ is a finitely generated free $\mathbb{Z}H$-module.
For each $i$, we choose a basis $(x_{i,j})_{j\in J_{i}}$
  of $F_{i}$. Then the morphism $f_{i}$ is determined by
  \begin{displaymath}
    f_{i}(x_{i,j})=\sum_{\substack{k\in J_{i-1}\\
        h \in H}}n^{k,h}_{j}h(x_{i-1,k}) \; .
  \end{displaymath}
  For any object $A$ of $H \text{-}\mathcal{A} $, we write
  $$S^{i}(A)=\prod_{j\in J_{i}}A, $$ and let
  $d^{i-1}:S^{i-1}(A) \to S^{i}(A)$ be the morphism
  determined by the family of morphisms
  \begin{displaymath}
    d^{i-1}_{j}:S^{i-1}(A)\longrightarrow A,\quad j\in J_{i} \; ,
  \end{displaymath}
  given by
  $$d^{i-1}_{j}=\sum_{\substack{k\in J_{i-1}\\ h \in H}}
  n^{k,h}_{j} h \circ p_{i-1,k} \; ,$$
  where $p_{i-1,k}$ is the projection of $S^{i-1}(A)$ onto the $k$th
  factor. Then there is a natural isomorphism of complexes
  $S^{\ast}(A)=\Hom(F_{\ast},A)^{H}$. Since the complex
  $S^{\ast}(A)$ determines an element of
  $D^{+}(\mathcal{A})$, we have proved (a). The proof
  of (b) is analogous. Note that our construction is canonical, since
  it does not depend on the choice of the resolution $F_\ast$, any two such
  choices being homotopically equivalent.
\end{Proof}

\begin{Rem} \label{3P}
(a) A natural question to ask is whether
under the above finiteness conditions on
$H$, the functors $R\Gamma(H,\argdot)$ of Theorem~\ref{3O}~(a)
and (b) are the actual
derived functors of the functor
\[
\Gamma(H,\argdot): H \text{-}\mathcal{A}  \longrightarrow \mathcal{A} \; .
\]
The answer in general is negative, and counterexamples occur right in the
context of arithmetic groups: let $\mathcal{A}$ be the category of finite
dimensional
vector spaces over $\BC$. Choose a connected, simply connected algebraic
group $P$ over $\BQ$, which is simple over $\BQ$, and of $\BQ$-rank at least
two. Let $H$ be an arithmetic subgroup of $P(\BQ)$.
Then $H \text{-}\mathcal{A}$,
the category of abstract representations in finite dimensional $\BC$-vector
spaces, is semi-simple: indeed, for two objects $\BV$ and $\BW$ of
$H \text{-}\mathcal{A}$, we have
\[
\Ext^1_{H \text{-}\mathcal{A}} (\BV,\BW) = H^1 (H, \BV^* \otimes_\BC \BW) \; ,
\]
and the latter group is zero by \cite[Cor.~2 of Thm.~2]{Rg}.
Therefore, any additive functor on $H \text{-}\mathcal{A}$ is automatically
exact. In particular, the derived functor of $\Gamma(H,\argdot)$ takes
the value $\Gamma(H,\BV) [0]$ on any object $\BV$ of $H \text{-}\mathcal{A}$.
On the other hand, if $H$ is neat in $P (\BQ)$,
then it is of type $FL$ by \cite[11.1~(c)]{BS}, hence Theorem~\ref{3O}~(b)
is applicable; but there exist examples of such $H$, and objects $\BV$ of
$H \text{-}\mathcal{A}$, for which
\[
\{ p \ge 1 \tei H^p (H,\BV) \ne 0 \}
\]
is not empty (e.g., \cite[Prop.~11.3~(b)]{BS}).\\[0.2cm]
(b) If the category $\mathcal{A}$ is Artinian, then the natural functor
\[
D^{+}(\mathcal{A}) \longto D^{+}(\Pro(\mathcal{A}))
\]
is a full embedding. Its image consists of the
complexes whose cohomology objects
lie in $\mathcal{A}$, the sub-category of Artinian
objects of $\Pro(\mathcal{A})$.
Consequently, the conclusions of (a) and (b)
of Theorem~\ref{3O} are equivalent to the following:
(a')~for any object $A$ of $H \text{-}\mathcal{A}$,
the group cohomology objects
$H^p(H,A)$ lie all in $\mathcal{A}$;
(b')~for any object $A$ of $H \text{-}\mathcal{A}$,
the group cohomology objects
$H^p(H,A)$ lie all in $\mathcal{A}$, and are trivial for large $p$.
If (a') or (b') is satisfied, then the respective
functor $R\Gamma(H,\argdot)$ is uniquely
determined by the commutative diagram in \ref{3O}.
\end{Rem}

We note the following consequence of Theorem~\ref{3O}:

\begin{Cor} \label{3Q}
Assume that the group $H$ is of type $FP_{\infty}$,
and that the Abelian category
$\mathcal{A}$ is semi-simple. Let $A$ be an object of $H \text{-}\mathcal{A}$,
and $p \in \BZ$.
Then any irreducible factor of $H^p (H,A) \in \Ob(\mathcal{A})$
is an irreducible factor of
$\For (A)$.
\end{Cor}

\begin{Proof}
By Theorem~\ref{3O}~(a), the object $H^p (H,A)$ is the cohomology object
of a complex, all of whose components are
finite products of copies of $\For (A)$.
\end{Proof}

Next we state the compatibility of group cohomology with respect to exact
functors.  Let $\mu :\mathcal{A} \to \mathcal{B}$ be an exact
functor between Abelian categories. We denote by the same symbol $\mu$ the
induced functor between the categories $H \text{-}\mathcal{A}$ (resp.
$\Pro(\mathcal{A})$, $H\text{-}\Pro(\mathcal{A})$) and
$H\text{-}\mathcal{B}$ (resp.
$\Pro(\mathcal{B})$, $H\text{-}\Pro(\mathcal{B})$).
The proof of the following
result is immediate and left to the reader.

\begin{Prop} \label{3R}
  Let $\mu :\mathcal{A} \to \mathcal{B}$ be an exact functor
  between Abelian categories. Then there is a
  natural commutative diagram
\[
\vcenter{\xymatrix@R-10pt{
D^{+}(H \text{-}\Pro(\mathcal{A}))
                                \ar[rr]^-{R\Gamma(H,\argdot)} \ar[d]^{\mu} &&
                    D^{+}(\Pro(\mathcal{A}))\ar[d]^{\mu} \\
                    D^{+}(H \text{-}\Pro(\mathcal{B}))
                                \ar[rr]^-{R\Gamma(H,\argdot)} &&
                    D^{+}(\Pro(\mathcal{B}))
\\}}
\]

  If $H$ is of type $FL$ or $FP_{\infty}$, then there are
  natural commutative diagrams
\[
\vcenter{\xymatrix@R-10pt{
D^{?}(H \text{-}\mathcal{A})
                         \ar[rr]^-{R\Gamma(H,\argdot)} \ar[d]^{\mu} &&
              D^{?}(\mathcal{A})\ar[d]^{\mu} \\
              D^{?}(H \text{-}\mathcal{B})
                         \ar[rr]^-{R\Gamma(H,\argdot)} &&
              D^{?}(\mathcal{B})
\\}}
\]
  with $?=b$ if $H$ is of type $FL$, and $?=+$ if $H$ is of type
  $FP_{\infty}$.
\end{Prop}

\begin{Var} \label{3S}
(a) There are obvious variants of \ref{3G}--\ref{3R} for Abelian ca\-te\-gories
$\mathcal{A}$ which are closed under arbitrary products. More precisely, in this
case, the use of the pro-category $\Pro(\mathcal{A})$ is unnecessary, and the
constructions and statements of \ref{3G}--\ref{3R} remain valid when the
symbol $\Pro(\mathcal{A})$ is replaced by $\mathcal{A}$.\\[0.2cm]
(b) Consider the case
when the group $H$ is normal in a larger group $L$. Then we may study the
fixed point functor
\begin{displaymath}
\Gamma(H,\argdot) = (\argdot)^{H}:
               L \text{-} \Pro(\mathcal{A}) \longrightarrow
               L/H \text{-} \Pro(\mathcal{A})
\end{displaymath}
defined in the same way as in \ref{3A}. The resolution functor is the functor
$C_L$ of \ref{3G}. The analogue of Proposition~\ref{3K}~(b) reads as follows:
the functor $(\argdot)^H \circ C_L$ is isomorphic to $C_{L/H} \circ \For$.
In particular, it is exact. Therefore, the analogue of Theorem~\ref{3L} holds:
the above fixed point functor is right derivable, and
\[
R\Gamma(H,\argdot) = (\argdot)^{H} \circ C_L^\ast \; .
\]
Furthermore, by Proposition~\ref{3N}~(b) (applied to the
group $H$ and the free $\BZ H$-resolution $C_\ast^L (\BZ)$ of
Example~\ref{3F}),
we see that the diagram
\[
\vcenter{\xymatrix@R-10pt{
D^{+}(L \text{-}\Pro(\mathcal{A}))
                      \ar[rr]^-{R\Gamma(H,\argdot)} \ar[d]_{\Res^L_H} &&
            D^{+}(L/H \text{-}\Pro(\mathcal{A})) \ar[d]^{\Res^{L/H}_{\{e\}}} \\
            D^{+}(H \text{-}\Pro(\mathcal{A}))
                      \ar[rr]^-{R\Gamma(H,\argdot)} &&
            D^{+}(\Pro(\mathcal{A}))
\\}}
\]
is commutative.
Finally, the analogues of Theorem~\ref{3O}~(a), (b) hold if the group $L$
is of type $FP_\infty$, resp.\ of type $FL$.\\[0.2cm]
(c) The construction of (b) continues to work, and the
statements made in (b) continue to hold in a somewhat larger generality.
Namely, let $\mathcal{A}$ be an Abelian category, on which the action of an
abstract group $L$ is given. This means that there are given
contravariant functors $\gamma^*$ on $\mathcal{A}$, for $\gamma \in L$,
such that
$(\gamma_1 \cdot \gamma_2)^* = \gamma_2^* \circ \gamma_1^*$ for all
$\gamma_1,\gamma_2 \in L$, and such that $e^* = \Id$.
We denote by $L \text{-}\mathcal{A}$ the category of pairs
\[
(A, (\rho_\gamma)_{\gamma \in L}) \; ,
\]
where $A \in \Ob(\mathcal{A})$, and $(\rho_\gamma)_{\gamma \in L}$ is a
family of
isomorphisms
\[
\rho_\gamma: \gamma^* A \isoto A
\]
in $\mathcal{A}$ such that the
cocycle condition holds. In the same way, define
the category $L \text{-}\Pro(\mathcal{A})$.
We assume that the action of a given normal
subgroup $H$ of $L$ on $\mathcal{A}$ is trivial: $\gamma^* = \Id$
for all $\gamma \in H$.
Therefore, the action of $L$ on $\mathcal{A}$ is induced by an action of the
quotient $L/H$.
The fixed point functor
\begin{displaymath}
\Gamma(H,\argdot) = (\argdot)^{H}:
               L \text{-}\Pro(\mathcal{A}) \longrightarrow
               L/H \text{-}\Pro(\mathcal{A})
\end{displaymath}
is defined by the same formula as in \ref{3A}. The resolution
functor is the functor $C_L$ of (b).
The action of an element $\gamma \in L$,
\begin{displaymath}
  \rho_{\gamma}: \gamma^{\ast} C_L(A) \isoto C_L A
\end{displaymath}
is determined by the family of morphisms
$\gamma_{h}= \rho_{\gamma} \circ \gamma^*p_{\gamma h}: \gamma^* C_L(A) \to A$.
\end{Var}

It is clear that there is a variant of Proposition~\ref{3N}~(b) in
the setting of Variant~\ref{3S}~(c).
We quote the precise result for further reference:

\begin{Prop}\label{3Ta}
  Let $\mathcal{A}$ be an Abelian category with an action of a group
  $L$. Let $H$ be a normal subgroup of $L$, which acts trivially
  on $\mathcal{A}$. Let
  $F_{\ast}$ be a free $\mathbb{Z}L$-resolution of the trivial
  $L$-module $\mathbb{Z}$. Then the functor
\[
R\Gamma (H,\argdot): D^+ (L \text{-}\Pro(\mathcal{A})) \longrightarrow
               D^+ (L/H \text{-}\Pro(\mathcal{A}))
\]
is represented by the functor $(\Hom(F_{\ast},\argdot))^{H}$.
\end{Prop}

It is possible to further enlarge the degree of generality
by imposing conditions
on the action of a second normal subgroup $H'$ of the group
$L$ in \ref{3S}~(c).
This applies in particular to the situation considered in Definition~\ref{2A},
where $H' = G_1(\BQ)$, and $(\Rep_{F} G_{1},\HQ)$ is the full sub-category
of $\HQ \text{-}\Rep_{F} G_{1}$ of objects satisfying condition
\ref{2A}~(a)~(i).
We want to derive the functor $\Gamma(\HC,\argdot)$ of
Definition~\ref{2B}.
In this case, we use the resolution functor
\begin{displaymath}
  C_{\Delta_1}: (\Pro(\Rep_{F} G_{1}),\HQ) \longrightarrow
                (\Pro(\Rep_{F} G_{1}),\HQ)
\end{displaymath}
given by
\begin{displaymath}
  C_{\Delta_1}(\BV_1) = \prod_{h \in \Delta_1} \BV_1 \; .
\end{displaymath}
Recall that by definition, we have $\Delta_1 = \HQ/G_{1} (\BQ)$,
and $\Delta = \HQ/G_{1} (\BQ) \HC$.
The composition $\Gamma(\HC,\argdot) \circ C_{\Delta_1}$ maps $\BV_1$
to $\prod_{\overline{h} \in \Delta} \BV_1$, and hence is exact.
Applying freely the results obtained so far, we get:

\begin{Thm} \label{3T}
(a) The functor
\begin{displaymath}
  \Gamma (\HC,\argdot):(\Pro(\Rep_F G_{1}),\HQ)
  \longrightarrow (\Pro(\Rep_F G_{1}),\HQ/\HC)
\end{displaymath}
is right derivable:
\begin{displaymath}
  R\Gamma (\HC,\argdot): D^{+}(\Pro(\Rep_F G_{1}),\HQ) \longrightarrow
                         D^{+}(\Pro(\Rep_F G_{1}),\HQ/\HC) \; .
\end{displaymath}
\noindent (b) The functor $R\Gamma (\HC,\argdot)$ respects the
sub-categories $D^{b}((\Rep_F G_{1}), ? )$: there is a commutative diagram
\[
\vcenter{\xymatrix@R-10pt{
D^{b}(\Rep_F G_{1},\HQ) \ar[rr]^-{R\Gamma(\HC,\argdot)} \ar[d] &&
              D^{b}(\Rep_F G_{1},\HQ/\HC) \ar[d] \\
      D^{+}(\Pro(\Rep_F G_{1}),\HQ) \ar[rr]^-{R\Gamma(\HC,\argdot)} &&
      D^{+}(\Pro(\Rep_F G_{1}),\HQ/\HC)
\\}}
\]
\end{Thm}

\begin{Proof}
Part (a) follows from the general formalism developed above.  For
(b), we intend to apply the criterion of Theorem~\ref{3O}.  Because
of the form of our resolution functor $C_{\Delta_1}$, we have to
impose the finiteness condition on the group $\Delta_1$ (see
Variants~\ref{3S}~(b) and (c)).  More precisely, we need to know
that $\Delta_1$ is of type $FL$.  By definition, this group is an
arithmetic subgroup of $Q/P_1 (\BQ)$, which is neat since $K$ is. By
\cite[11.1~(c)]{BS}, such a group is indeed of
type $FL$.
\end{Proof}

We still need to prove what was left open in Section~\ref{2}:\\

\begin{Proofof}{Proposition~\ref{2Ba}}
Part (a) is a special case of Proposition~\ref{3R}, and (b) follows from
semi-simplicity of the category $\Rep_F G_{1}$, and from Corollary~\ref{3Q}.
\end{Proofof}


\bigskip
%
%

\section{Equivariant algebraic Hodge modules}
\label{3b}

The aim of this short section is to develop
some ele\-mentary theory of equivariant algebraic Hodge modules. \\

Because of the local nature of Hodge modules, the category
$\MHM_F X$ can be defined for reduced schemes $X$, which are only \emph{locally}
of finite type over $\BC$. However, the formalism of Grothendieck's
functors \cite[Section~4]{Sa} is constructed on the bounded derived categories
of Hodge modules on reduced schemes which are (globally) of finite type over
$\BC$. It does not obviously extend to the derived categories
of Hodge modules on schemes of a more general type.

\begin{Def} \label{5C}
Let $X$ be a reduced scheme which is locally of finite type over $\BC$,
and $H$ an abstract group acting on $X$ by algebraic automorphisms.
The category $H \text{-} \MHM_F X$ consists of pairs
 \[
(\BM, (\rho_\gamma)_{\gamma \in H}) \; ,
\]
where $\BM \in \MHM_F X$, and $(\rho_\gamma)_{\gamma \in H}$
is a family of isomorphisms
\[
\rho_\gamma: \gamma^* \BM \isoto \BM
\]
in $\MHM_F X$ such that the cocycle condition holds.
\end{Def}

Note that this is a special case of what was done in Variant~\ref{3S}~(c).
We shall repeatedly use the following principle:

\begin{Prop} \label{5D}
In the situation of Definition~\ref{5C}, suppose that the action of
$H$ on $X$ is free and proper in the sense of \cite[(1.7)]{P2}, with
quotient $H \backslash X$. Denote by $\Pi$ the
morphism from $X$ to $H \backslash X$. Then the inverse image
\[
\Pi^*: \MHM_F (H \backslash X) \longto H \text{-} \MHM_F X
\]
is an equivalence of categories, which possesses a canonical pseudo-inverse.
\end{Prop}

\begin{Proof}
The pseudo-inverse is given
by the direct image $\Pi_*$, followed by the $H$-invariants
$\Gamma(H,\argdot)$. Since the direct
image is not in general defined for morphisms which are only locally of finite
type, this definition needs to be explained:
choose an $H$-equivariant open covering $\FV$ of $X$, such that each open
subset $V$ in $\FV$ satisfies
\[
\Pi^{-1} (\Pi (V)) = \coprod_{h \in H} h(V) \; .
\]
This is possible because of our assumption on the action of $H$.
It is then clear how to define the restriction of
$\Gamma(H,\argdot) \circ \Pi_*$ to any open subset in the quotient
$H \backslash \FV$. The resulting collection of Hodge modules glues to
give a Hodge module on $H \backslash X$.
\end{Proof}

\begin{Cor} \label{5E}
In the situation of Proposition~\ref{5D}, the inverse image
\[
\Pi^*: D^b (\MHM_F (H \backslash X)) \longto D^b (H \text{-} \MHM_F X)
\]
is an equivalence of categories, which possesses a canonical pseudo-inverse.
\end{Cor}

\begin{Rem} \label{5F}
Using the formalism developed in Section~\ref{3}, we can give a more
conceptual meaning of the canonical pseudo-inverse
\[
D^b (H \text{-} \MHM_F X) \longto D^b (\MHM_F (H \backslash X))
\]
of Corollary~\ref{5E}.
As in the proof of \cite[Thm.~4.3]{Sa}, it is possible, using a covering as
in the proof of Proposition~\ref{5D}, to define the direct image
\[
\begin{array}{cccc}
\Pi_*: D^b (H \text{-} \MHM_F X)
& \! \longto  \! &
       D^b (H\text{-}\Pro(\MHM_F (H \backslash X)))    & \\
& \! \subset  \! &
       D^+ (H\text{-}\Pro(\MHM_F (H \backslash X)))  & .
\end{array}
\]
Its image consists of $\Gamma(H,\argdot)$-acyclic complexes.
The composition of $\Pi_*$ and the functor
\[
R\Gamma (H,\argdot): D^+ (H\text{-}\Pro(\MHM_F (H \backslash X)))
\longto D^+ (\Pro(\MHM_F (H \backslash X)))
\]
of Theorem~\ref{3L} factors through $D^b (\MHM_F (H \backslash X))$.
Our quasi-inverse is the functor
\[
D^b (H \text{-} \MHM_F X) \longto
D^b (\MHM_F (H \backslash X))
\]
induced by the composition $R\Gamma (H,\argdot) \circ \Pi_*$.
\end{Rem}


\bigskip
%
%

\section{\v{C}ech complexes for Hodge modules}
\label{3aa}

For later purposes, we need
to develop the basics of the formalism of \v{C}ech complexes
associated to closed coverings
in the context of Hodge modules.\\

Fix a reduced scheme $Z$, which is separated and of finite type
over $\BC$. Let $\FZ = \{ Z_\sigma \}_{\sigma \in
  \Sigma }$ be a finite covering of $Z$
by reduced closed sub-schemes, not necessarily different from
each other. We
denote by $\FZ_{\bullet}$ the free simplicial set generated by the set
of indices $\Sigma $. That is,
$\FZ_{p}=\Sigma ^{p+1}$ is the set of
$p+1$-tuples $(\sigma _{0}, \dots , \sigma _{p})$. If
\[
\tau :\{0,\ldots,q\}\longrightarrow \{0,\ldots,p\}
\]
is an increasing map and $I=(\sigma _{0}, \dots , \sigma _{p})\in
\FZ_{p}$, then
\[
\FZ_{\bullet}(\tau )(I)=
(\sigma _{\tau (0)}, \dots , \sigma _{\tau (q)}) \in \FZ_q \; .
\]

\begin{Def} \label{5Ga}
Define the Abelian category
$(\MHM_F Z)^{\FZ}_\infty$ as the category of mixed Hodge modules over the
simplicial scheme $Z\times \FZ_{\bullet}$.
\end{Def}

Explicitly, an element of $(\MHM_F Z)^{\FZ}_\infty$ is a family $(\BM_{I})_I$ of
objects of $\MHM_F Z$ indexed by $\FZ_{\bullet}$, and for every
increasing map
\[
\tau :\{0,\ldots,q\} \longto \{0,\ldots,p\}
\]
and each $I\in \FZ_{p}$, a morphism
\[
\tau_{I}:\BM_{\FZ_{\bullet}(\tau )(I)}\longrightarrow
\BM_{I} \; ,
\]
equal to the identity if $\tau = \Id_{\{0,\ldots,q\}}$, and such that
\[
(\eta \circ \tau)_{I} = \eta_{I} \circ \tau_{\FZ_{\bullet}(\eta)(I)} \; .
\]
\forget{
\? {\tt Jos\'e, could you check what follows now (until the end of this
section...)} \?
}

\begin{Def} \label{5G}
Define $(\MHM_F Z)^{\FZ}$ as the full Abelian sub-category of
$(\MHM_F Z)^{\FZ}_\infty$ consisting of objects
\[
\left( (\BM_{I})_I , (\tau_{I})_{\tau,I} \right)
\]
satisfying the following property:
\[
\tau_{I}:\BM_{\FZ_{\bullet}(\tau )(I)}\longrightarrow
\BM_{I}
\]
is an isomorphism for any increasing
$\tau :\{0,\ldots,q\} \to \{0,\ldots,p\}$, and
for any $I\in \FZ_{p}$ such that the subsets of $\Sigma$ underlying
the $q+1$-tuple $\FZ_{\bullet}(\tau )(I)$ and the $p+1$-tuple $I$
are the same.
\end{Def}

\begin{Rem} \label{5Gb}
By definition, the components of an object of the category $(\MHM_F Z)^{\FZ}$
represent a finite number of isomorphism classes of Hodge modules on $Z$.
\end{Rem}

Observe that the theory of mixed Hodge modules over general
simplicial schemes is not well established because for general morphisms,
inverse images of
mixed Hodge modules are only defined in the derived category.
However, in our situation, there is no problem since all the
morphisms of the simplicial scheme $Z\times \FZ_{\bullet}$ are
given by the identity on $Z$. \\

Observe also that we can define $(\MHM_F Z)^{\FZ}$ for locally finite
$Z$ and infinite coverings $\FZ$, or even for any simplicial set
$\FZ_{\bullet}$ not
necessarily associated to a \v{C}ech covering.
An object $(\BM_{I})_I$ of $(\MHM_{F} Z)^{\FZ}$ defines a
co-simplicial object, denoted $\BM_{\bullet}$, of the category $\MHM_{F}
Z$ if $\FZ$ is finite,
and a co-simplicial object of the category $\Pro(\MHM_{F} Z)$ if $\FZ$ is
infinite: put
\[
\BM_{p}=\prod_{I\in \FZ_{p}}\BM_{I} \; ,
\]
with the induced morphisms. \\

We go back to the hypothesis of Definition \ref{5G}. Thus, $Z$ is of
finite type, and $\FZ$ is finite.
The following observation will be used repeatedly:

\begin{Prop} \label{3aaA}
Let $f: (\BM_I)^\bullet_I \to (\BN_I)^\bullet_I$ be a morphism
in the ca\-te\-gory
$D^b ((\MHM_F Z)^{\FZ})$. Then
$f$ is an isomorphism if and only if
$f_I: \BM^\bullet_I \to \BN^\bullet_I$
is an isomorphism in $D^b (\MHM_F Z)$, for all
$I \in \FZ_p$, and all $p \ge 0$.
\end{Prop}

Next, we need to define functors
\[
\vcenter{\xymatrix@R-10pt{
D^b (\MHM_F Z) \ar@/^2pc/[r]^-{S_\bullet} &
D^b ((\MHM_F Z)^{\FZ}) \ar@/^2pc/[l]^-{\Tot}
\\}}
\]
The functor $\Tot$ is induced by the exact functor
\[
(\MHM_F Z)^{\FZ} \longto C^+ (\MHM_F Z) \; ,
\]
denoted by the same symbol, that sends a co-simplicial mixed Hodge module
$(\BM_I)_I$ to the normalized cochain complex of mixed Hodge modules
associated to $\BM_\bullet$. Note that by definition of the category
$(\MHM_F Z)^{\FZ}$, the resulting functor
\[
\Tot: D^b ((\MHM_F Z)^{\FZ}) \longto D^+ (\MHM_F Z)
\]
factorizes through $D^b (\MHM_F Z)$.
\forget{
\? {\tt really} \?
} \\

The construction of the functor $S_\bullet$ depends on
the covering $\FZ$, not only on the index set $\Sigma$.
For any element $I=(\sigma _{0}, \dots , \sigma _{p})$ of $\FZ_{p}$,
set
\[
Z_I := \bigcap_{k=0}^{p}Z_{\sigma _{k}} \; ,
\]
with its
reduced scheme structure. We shall write $i_I$ for the closed immersion
of $Z_I$ into $Z$. For any increasing map $\tau $, if
$J=\FZ_{\bullet}(\tau )(I)$ then $Z_{I}$ is a closed subset of $Z_{J}$.
The basic idea for the construction of the functor
\[
S_\bullet: D^b (\MHM_F Z) \longto D^b ((\MHM_F Z)^{\FZ})
\]
is to associate to a
complex $\BM^\bullet$ of Hodge modules on $Z$ the class of a certain
complex of Hodge modules over
$Z\times \FZ_{\bullet}$ that restricts to
\[
(i_I)_* \, i_I^* \, \BM^\bullet
\]
on the component $Z\times \{I\}$. In order to do this rigorously,
we recall the definition of $(i_I)_* \, i_I^*$
in Saito's formalism \cite[(4.4.1)]{Sa}: \emph{choose} an
open affine
covering of the complement $j_I: U_I \into Z$ of $Z_I$,
and use the \v{C}ech complex associated to that covering
to define the functor
$(j_I)_! \, j_I^*$ on the level of complexes, together with a
transformation $(j_I)_! \, j_I^* \to \Id$. The functor
\[
(i_I)_* \, i_I^* : C^b (\MHM_F Z) \longto C^b (\MHM_F Z)
\]
maps a complex $\BM^\bullet$ to the simple complex associated to
\[
(j_I)_! \, j_I^* \BM^\bullet \longto \BM^\bullet \; .
\]
This construction descends to the level of derived categories,
and the induced functor
\[
(i_I)_* \, i_I^* : D^b (\MHM_F Z) \longto D^b (\MHM_F Z)
\]
does not depend on the choice of the affine covering of $U_I$.\\

In our situation, we can choose the affine coverings for the different
closed sub-schemes in such a
way that, for every
inclusion $Z_{I}\subset Z_{J}$, any open affine subset
occurring in the covering of $U_{J}$ is
contained in an open affine subset occurring in the covering
of $U_{I}$. This choice induces a
compatible set of morphisms
$(i_J)_* \, i_J^* \, \BM^\bullet \rightarrow (i_I)_* \,
i_I^* \, \BM^\bullet $ at the level of complexes.
Putting $S_{I}(\BM^{\bullet})=(i_I)_* \, i_I^* \, \BM^\bullet $,
we thus obtain a functor
\[
S_\bullet := ((i_I)_* \, i_I^*)_I:
                    C^b (\MHM_F Z) \longto C^b ((\MHM_F Z)^{\FZ}) \; .
\]
This construction descends to the level of derived categories.
The induced functor is independent of the choices.

\begin{Rem} \label{5H}
For later use, it will be important to observe that the above construction
also defines a filtered version of the functor $S_\bullet$:
\[
S_\bullet: DF^b (\MHM_F Z) \longto DF^b ((\MHM_F Z)^{\FZ}) \; ,
\]
where $DF^b$ denote the \emph{filtered bounded derived categories}
used in \cite[Section~3]{BBD} and \cite[Appendix~A]{B}. Thus, the term
``bounded'' refers to boundedness of the complexes as well as finiteness
of the filtrations.
\end{Rem}

\begin{Prop} \label{5I}
There is a canonical isomorphism of functors
\[
\Id \cong \Tot \circ S_\bullet: D^b (\MHM_F Z) \longto D^b (\MHM_F Z) \; .
\]
\end{Prop}

\begin{Proof}
By construction, the functors $(i_I)_* \, i_I^*$ come with natural
transformations $\Id \to (i_I)_* \, i_I^*$, which
induce a natural transformation $\Id \to \Tot \circ S_\bullet$. That it is
an isomorphism can be checked after application of the forgetful
functor to the bounded derived category of perverse sheaves on $Z(\BC)$.
By \cite[Main~Thm.~1.3]{B}, this latter category can be identified
with a full sub-category of the derived category
of Abelian sheaves on $Z(\BC)$. Thus, our claim
follows from the fact that the \v{C}ech complex of any
sheaf $\FF$ is a resolution of $\FF$.
\end{Proof}

We need to discuss functoriality of our constructions.
Let $p: Z \to Y$ be a morphism of reduced schemes,
which are separated and of finite type
over $\BC$, and assume given finite coverings
$\{ Z_\sigma \}_{\sigma \in \Sigma }$ and
$\{ Y_\sigma \}_{\sigma \in \Sigma }$
of $Z$ and $Y$, respectively.
Since the index set $\Sigma$ is the same for the two coverings,
we shall write $(\MHM_F Z)^{\FZ}$ and $(\MHM_F Y)^{\FZ}$
for the respective categories defined in \ref{5G}.
We have the following co-simplicial version of
direct images under $p$:

\begin{Prop} \label{5Q}
(a) There is a canonical functor
\[
p^{\FZ}_*:
D^b ((\MHM_F Z)^{\FZ}) \longto
D^b ((\MHM_F Y)^{\FZ}) \; .
\]
(b) Let $q \ge 0$, and $J \in \FZ_q$.
There is a natural commutative diagram
\[
\vcenter{\xymatrix@R-10pt{
D^b ((\MHM_F Z)^{\FZ})
               \ar[rr]^-{(\BM_I)^\bullet_I \mapsto \BM^\bullet_J}
               \ar[d]_{p^{\FZ}_*} &&
                            D^b (\MHM_F Z) \ar[d]^{p_*} \\
D^b ((\MHM_F Y)^{\FZ})
               \ar[rr]^-{(\BN_I)^\bullet_I \mapsto \BN^\bullet_J} &&
                            D^b (\MHM_F Y)
\\}}
\]
(c) There is a natural commutative diagram
\[
\vcenter{\xymatrix@R-10pt{
D^b ((\MHM_F Z)^{\FZ})
                 \ar[r]^-{\Tot} \ar[d]_{p^{\FZ}_*} &
                            D^b (\MHM_F Z) \ar[d]^{p_*} \\
D^b ((\MHM_F Y)^{\FZ})
                 \ar[r]^-{\Tot} &
                            D^b (\MHM_F Y)
\\}}
\]
\end{Prop}

\begin{Proof}
In order to define $p^{\FZ}_*$, we recall
part of the definition
of $p_*$ in Saito's formalism \cite[Thm.~4.3]{Sa}:
if $k: V \into Z$ is the immersion of an open affine subset,
then $(p \circ k)_*$ is the total left derived functor of the functor
$\CH^0 (p \circ k)_*$ (the definition of $\CH^0 (p \circ k)_*$
will not be recalled, since
it will not be needed). \emph{Choose}
a finite open affine covering
$\FV = \{ V_1,\ldots,V_r \}$ of $Z$.
Call a Hodge module $\BL$ on $Z$ \emph{$p_*$-acyclic with
respect to $\FV$} if the restriction of $\BL$ to
any intersection of the $V_l$ is
$(p \circ k)_*$-acyclic, where $k$ denotes the open immersion
of that intersection into $Z$.
We then have (see the proof of \cite[Thm.~4.3]{Sa}, which
in turn is based on \cite[Section~3]{B}):
\begin{enumerate}
\item[(1)] for any Hodge module $\BM$ on $Z$, there is an
epimorphism $\BL \onto \BM$, whose source is $p_*$-acyclic with
respect to $\FV$.
\item[(2)] for any system consisting of
Hodge modules $\BM_n$ representing a finite number of isomorphism
classes, the epimorphism in (1)
can be chosen functorially with respect to all
morphisms between the $\BM_n$. Indeed, this can be seen from
\cite[proof of Lemma~3.3]{B}; e.g., if $Z$
is quasi-projective, the $\BL_n$ can be chosen as $j_! j^* \BM_n$,
for the open immersion $j$ of some suitable affine open subset $U$ of $Z$
(the same for all $n$).
\end{enumerate}
Given a bounded complex $\BM^\bullet$
of Hodge modules, we use (1) and (2) to construct
a complex $\BL^\bullet$, all of whose components
are $p_*$-acyclic with respect to $\FV$,
and a morphism $\varphi: \BL^\bullet \to \BM^\bullet$, which becomes
an isomorphism in $D^-(\MHM_F Z)$.
Observe that $\BL^\bullet$
can be chosen to be bounded because of the finite cohomological
dimension of the $(p \circ k)_* \, k^*$.
(For later use, it will be important to note that furthermore,
the morphism of complexes
$\varphi$ can be chosen to be epimorphic in all degrees.)
Replace $\BL^\bullet$
by the \v{C}ech complex
$C(\BL^\bullet)^\bullet$ associated to $\FV$
(note that since the $k$ are affine, the $k_* \, k^*$ are exact).
We get an actual complex $p_* C(\BL^\bullet)^\bullet$, whose class in
$D^b (\MHM_F Y)$ does not depend on any of the choices. By definition,
this class is $p_* \BM^\bullet$.

Thanks to (2), and to Remark~\ref{5Gb}, the above can be imitated
on the level of complexes of simplicial objects.
This is the functor $p^{\FZ}_*$, and it satisfies properties
(b) and (c).
\end{Proof}

\begin{Cor} \label{5R}
There is a natural commutative diagram
\[
\vcenter{\xymatrix@R-10pt{
D^b (\MHM_F Z) \ar[r]^-{S_\bullet} \ar[d]_{p_*} &
             D^b ((\MHM_F Z)^{\FZ})
                                      \ar[d]^{p^{\FZ}_*} \\
D^b (\MHM_F Y) \ar@{<-}[r]^-{\Tot} &
             D^b ((\MHM_F Y)^{\FZ})
\\}}
\]
\end{Cor}

\begin{Proof}
This follows from Propositions~\ref{5Q}~(c) and \ref{5I}.
\end{Proof}


\bigskip
%
%

\section{Degeneration in relative torus embeddings}
\label{3a}

The aim of this short section is to study
the degeneration
of local systems, and of variations of Hodge structure in
a relative torus embedding.
We are going to use a number of concepts related to torus embeddings
as explained in \cite[Chap.~I]{KKMS} or \cite[5.1--5]{P1}.  \\

We shall consider the
following situation: $B$ is a scheme over $\BC$, and
$T$ a complex torus with cocharacter group $Y$. Fix a
smooth rational polyhedral decomposition $\Fs$ of $Y_\BR$.
Consider the (constant) torus $T_B$ over $B$.
We get a (partial) compactification $T_\Fs$ of $T_B$ relative to $B$,
which is naturally endowed
with a stratification indexed by the cones in $\Fs$.
More generally, this is true for any $T$-torsor $X$ over $B$.
Fix a cone $\theta \in \Fs$, and consider the
diagram
\[
\vcenter{\xymatrix@R-10pt{
X \ar@{^{ (}->}[r]^{\jmath} &
X_\Fs \ar@{<-^{ )}}[r]^-{\imath_\theta} &
X_\theta \ar@{<-^{ )}}[r]^-{\jmath_\theta} &
X_\theta^\circ \ar@/^2pc/[ll]^{\imath_\theta^\circ}
\\}}
\]
Here, $X_\Fs$ denotes the (partial)
compactification, $X_\theta^\circ$ the
stratum associated to $\theta$, and $X_\theta$ its closure.
One refers to $X_\Fs$ as the \emph{relative torus embedding} over $B$
associated to $\Fs$.
The stratum $X_\theta^\circ$ is itself a torsor under a complex torus
$T_\theta^\circ$, and there is a canonical isomorphism
\[
Y / (\langle \theta \rangle_\BR \cap Y) \isoto Y_\theta \; ,
\]
where $\langle \theta \rangle_\BR \subset Y_\BR$ denotes the subspace
generated by $\theta$, and $Y_\theta$ the
cocharacter group of $T_\theta^\circ$. Recall that the cocharacter group of a torus
is canonically identified with its fundamental group. The above isomorphism between
$Y / (\langle \theta \rangle_\BR \cap Y)$ and $Y_\theta$ is induced from a
projection from $Y$ to $Y_\theta$, which corresponds to the canonical projection
from $T$ to $T_\theta^\circ$.\\

Now let $\FF^\bullet$ be a complex of
$F$-linear local systems
on $X(\BC)$. Denote by $a$, resp.\ $a_\theta$, resp.\ $a_\theta^\circ$
the structure morphisms to $B$ from
$X$, resp.\ from $X_\theta$, resp.\ from $X_\theta^\circ$.
We have:

\begin{Prop} \label{5A}
(a) The adjunction morphism
\[
\imath_\theta^* \, \jmath_* \, \FF^\bullet \longto
{\jmath_\theta}_* \, (\imath_\theta^\circ)^* \, \jmath_* \, \FF^\bullet
\]
is an isomorphism in the derived category of
Abelian sheaves on $X_\theta (\BC)$. \\[0,2cm]
(b) Adjunction induces an isomorphism
\[
a_* \, \FF^\bullet \isoto
(a_\theta)_* \, \imath_\theta^* \, \jmath_* \, \FF^\bullet
\]
in the derived category of Abelian sheaves on $B (\BC)$.
\end{Prop}

\begin{Proof}
Since $a$ is locally a projection, we may assume that $B$ is a point,
and that $X = T$.
Claim~(a) can be shown after taking inverse images for all the strata
in the natural stratification of $T_\theta$. These correspond
to cones $\phi$ in $\Fs$ containing $\theta$ as a face. We have to show that
\[
(\imath_\phi^\circ)^* \, \jmath_* \, \FF^\bullet \longto
(\imath_\phi^\circ)^* {\jmath_\theta}_* \,
(\imath_\theta^\circ)^* \, \jmath_* \, \FF^\bullet
\]
is an isomorphism.
For this, we may assume (by passing to a subset of $\Fs$) that
$\phi$ is the unique open cone in $\Fs$, or equivalently, that
$T_\phi = T_\phi^\circ$ is the unique closed stratum of $T_\Fs$.
Then the structure morphism of $T_\Fs$ factors over the projection
to $T_\phi$, identifying $T_\Fs$ with a relative torus embedding over
$T_\phi$. As before, we may therefore assume that $T_\phi$ is a
point. Locally around this point, we can choose coordinates
$t_1, \ldots, t_n$,
and assume that the torus embedding
$T_\Fs$ equals $\BA^n$, with the canonical action of
$T = \Gnm \subset \BA^n$, that the stratification is the one
induced by the coordinates, that $T_\phi = \{ 0 \}$, and that the
intermediate stratum $T_\theta$ is defined by the vanishing of
the first $k$ of the $t_i$. We therefore have
\[
T_\theta^\circ = \{ (t_1, \ldots, t_n) \; , \; t_1 = \ldots = t_k = 0 \, , \,
t_{k+1} \cdots t_n \ne 0 \} \; .
\]
Now recall that the complexes of sheaves $(\imath_\theta^\circ)^* \, \jmath_*$,
$(\imath_\phi^\circ)^* \, {\jmath_\theta}_*$, and
$(\imath_\phi^\circ)^* \, \jmath_*$ can be
computed from direct limits over the (analytic)
neighbourhoods of $T_\theta$ in $T_\Fs$,
of $T_\phi$ in $T_\theta$, and of $T_\phi$ in $T_\Fs$, respectively. In each of these
direct systems, we can find a co-final system of neighbourhoods,
all of whose members are homotopically equivalent to each other.
If we evaluate on complexes of local systems,
we see that the direct limits over these co-final systems are constant.

Now denote by $U_i$ the image of the positively oriented generator of
the fundamental group of $\Gm (\BC)$
under the embedding of $\Gm$ into $\Gnm$ via the $i$-th coordinate.
Denote by $\Loc_F$ the category of $F$-linear local systems,
and by $\Sh$ the category of Abelian sheaves.
Identify local systems on $T (\BC)$ and $T_\theta^\circ (\BC)$
with representations of $Y = \langle U_1, \ldots, U_n \rangle_\BZ$
and of $Y_\theta = Y / (\langle \theta \rangle_\BR \cap Y) = \langle U_{k+1}, \ldots, U_n \rangle_\BZ$,
respectively.
For an abstract group $H$, denote by $R\Gamma (H, \argdot)$ the derived functor of
the $H$-invariants. From the above discussion, we see:
\begin{enumerate}
\item[(1)]
There is a commutative diagram of functors
\[
\vcenter{\xymatrix@R-10pt{
        D^{+}( \Loc_F T(\BC) )
                                   \ar[rr]^-{R\Gamma(\langle \theta \rangle_\BR \cap Y, \argdot)} \ar[d] &&
        D^{+}( \Loc_F T_\theta^\circ (\BC) ) \ar[d] \\
        D^{+}( \Sh T(\BC) )
                                   \ar[rr]^-{(\imath_\theta^\circ)^* \, \jmath_*} &&
        D^{+}( \Sh T_\theta^\circ (\BC) )
\\}}
\]
\item[(2)]
There is a commutative diagram of functors
\[
\vcenter{\xymatrix@R-10pt{
        D^{+}( \Loc_F T_\theta^\circ (\BC) )
                                   \ar[rr]^-{R\Gamma(Y_\theta,\argdot)} \ar[d] &&
        D^{+}( \Loc_F T_\phi (\BC) ) \ar[d] \\
        D^{+}( \Sh T_\theta^\circ (\BC) )
                                   \ar[rr]^-{(\imath_\phi^\circ)^* \, {\jmath_\theta}_*} &&
        D^{+}( \Sh T_\phi (\BC) )
\\}}
\]
\item[(3)]
There is a commutative diagram of functors
\[
\vcenter{\xymatrix@R-10pt{
        D^{+}( \Loc_F T (\BC) )
                                   \ar[rr]^-{R\Gamma(Y,\argdot)} \ar[d] &&
        D^{+}( \Loc_F T_\phi (\BC) ) \ar[d] \\
        D^{+}( \Sh T (\BC) )
                                   \ar[rr]^-{(\imath_\phi^\circ)^* \, \jmath_*} &&
        D^{+}( \Sh T_\phi (\BC) )
\\}}
\]
\item[(4)]
Under the identifications of diagrams (1)--(3),
the natural transformation \[
(\imath_\phi^\circ)^* \, \jmath_* \longto
(\imath_\phi^\circ)^* {\jmath_\theta}_* \circ
(\imath_\theta^\circ)^* \, \jmath_*
\]
of functors from
$D^{+}( \Sh T (\BC) )$ to $D^{+}( \Sh T_\phi (\BC) )$
restricts to the canonical isomorphism
\[
R\Gamma(Y,\argdot) \isoto
R\Gamma(Y_\theta,\argdot) \circ
R\Gamma(\langle \theta \rangle_\BR \cap Y, \argdot)
\]
of functors from $D^{+}( \Loc_F T (\BC) )$ to $D^{+}( \Loc_F T_\phi (\BC) )$.
\end{enumerate}
This shows claim~(a).
Claim~(b) follows from (a), and the fact that adjunction
induces an isomorphism
\[
a_* \, \FF^\bullet \isoto
(a_\theta^\circ)_* \, (\imath_\theta^\circ)^* \, \jmath_* \, \FF^\bullet \; .
\]
This in turn results from the fact that both sides are computed by $R\Gamma(Y,\argdot)$.
\end{Proof}

\begin{Cor} \label{5B}
Under the hypotheses of \ref{5A}, assume that $B$ is
smooth, and of finite type over $\BC$. Let
$\BM^\bullet$ be a complex of admissible variations of $F$-Hodge structure
on $X$. Adjunction induces an isomorphism
\[
a_* \, \BM^\bullet \isoto
(a_\theta)_* \, \imath_\theta^* \, \jmath_* \, \BM^\bullet
\]
in $D^b (\MHM_F B)$.
\end{Cor}

\begin{Proof}
Isomorphisms in the category $D^b (\MHM_F B)$ can be recognized after
application of the forgetful functor $rat$ to the bounded derived category
$D^b (\Perv_F B(\BC))$ of perverse sheaves on $B(\BC)$.
This is a formal consequence of exact- and faithfulness of $rat$ on the
level of Abelian categories
\[
rat: \MHM_F B \longto \Perv_F B(\BC)
\]
\cite[p.~222]{Sa}.
By \cite[Main~Thm.~1.3]{B}, the category $D^b (\Perv_F B(\BC))$ can be identified
with a full sub-category of the derived category
of Abelian sheaves on $B(\BC)$. So the claim follows from part (b) of
Proposition~\ref{5A}.
\end{Proof}


\bigskip
%
%

\section{Specialization of local systems,
and of variations of Hodge structure}
\label{3c}

In order to prove the part of Theorem~\ref{2D} concerning the
comparison of weight filtrations,
it will be necessary to recall the explicit description of
the \emph{nearby cycle functor}, as well as fundamental results on
\emph{nilpotent orbits}. The main result of this section is Theorem~\ref{3cY}.
It will be used in the form of Corollary~\ref{3cZ}, in the proof of
Proposition~\ref{5M}.
Because of the rather technical nature of the material, we chose to present
the main result first, and then recall the theory needed for its proof. \\

Throughout this section, assume that $X$ is a smooth analytic space,
and that $j: U \into X$ is the open immersion of a dense analytic subset,
such that the complement $Z'$ of $U$ is a divisor with normal crossings.
Let $Z$ be a closed analytic subspace of $Z'$, which is still a divisor in $X$.
Thus, locally on $X$, the set underlying $Z$ is the union of components
of $Z'$. \\

Recall Verdier's construction of the \emph{specialization functor} $Sp_Z$
in the ana\-lytic
context \cite[Section~9]{V2}. It preserves perversity \cite[(SP7)]{V2},
and can thus be seen as an exact functor
\[
\Perv_F X \longto \Perv_F N_{Z/X} \; ,
\]
where $N_{Z/X}$ denotes the (analytic) normal cone of $Z$ in $X$.
Recall from \cite[(SP1)]{V2}
that the image of $Sp_Z$ is contained in the category of \emph{monodromical}
perverse sheaves on $N_{Z/X}$ \cite[p.~356]{V2}. Perverse sheaves of the
form $Sp_Z \BW$ are thus equipped with a canonical \emph{monodromy
automorphism} $T$, and hence also with a nilpotent \emph{monodromy endomorphism},
namely, the logarithm of the unipotent part of $T$.\\

When $Z$ is a principal divisor, defined by a holomorphic function $g$,
then we also have the nearby cycle functor $\psi_g$,
which respects perversity up to a shift by $[-1]$
\cite[Section~3, Claim~4)]{V3}, and is exact. Let us write
\[
\psi_g^p := \psi_g [-1]: \Perv_F X \longto \Perv_F Z \; .
\]
By the very definition of $\psi_g^p$, perverse sheaves in its image
are again cano\-ni\-cally
equipped with a monodromy automorphism and a monodromy endomorphism.
By \cite[(SP6)]{V2},
$\psi_g^p$ and its monodromy automorphism can be recovered from $Sp_Z$.\\

We shall
be particularly interested in the composition of $Sp_Z$, resp.\ of $\psi_g^p$
with the functor
\[
j_*: \Perv_F U \longto \Perv_F X \; ,
\]
which respects perversity, since $j$ is affine. \\

For an object $\BE$ of an Abelian category $\mathcal{A}$,
equipped with a nilpotent endomorphism $N$,
recall the notion of \emph{monodromy weight filtration of $N$} on $\BE$
\cite[(1.6.1)]{D3}.
If $\BE$ is equipped with a finite ascending filtration $W_\bullet$,
then one defines the \emph{monodromy weight filtration of $N$ relative to
$W_\bullet$} on $\BE$ \cite[(1.6.13)]{D3}. (Caution! This latter filtration does not
always exist.)   \\

Set $D := \{ x \in \BC \tei |x| < 1 \}$, and
$D^* := D - \{ 0 \}$.
Fix a point $z$ of $Z'$. If locally around $z$, the divisor $Z'$ is the union of $m$
smooth components, then the fundamental group of $V \cap U$, for
small neighbourhoods $V$ of $z$ in $X$ isomorphic to $D^n$, is
free Abelian of rank $m$, and independent of $V$.
Let us refer to this group as the \emph{local monodromy group around
$Z'$ at $z$}.
Call a $\BZ$-base $\FT = (T_1, \ldots, T_m)$ of the local mo\-no\-dromy group
\emph{adapted to local coordinates} if
there is an isomorphism of $D^n$ onto a neighbourhood $V$
of $z$ identifying $(D^*)^m \times D^{n-m}$ with $V \cap U$,
and the canonical $\BZ$-base of $\pi_1 ((D^*)^m \times D^{n-m})$
with $(T_1, \ldots, T_m)$.  \\

We say that a local system $\BV$
on $U$ has \emph{unipotent local monodromy around $Z$}
if for any point $z$ of $Z \subset Z'$, the action of the
local monodromy group around $Z'$ at $z$ is unipotent.
Since the local monodromy group is Abelian,
its elements act as automorphisms of local systems on the restriction
of $\BV$ to $V \cap U$, for $V$ as above. Let us write
$\BV_{\FU^\circ(z)}$ for the direct limit of these restrictions.

\begin{Def} \label{3cB}
Let $\BV$ be a local system on $U$ with unipotent local mo\-no\-dromy around $Z$,
and $W_\bullet$ a finite filtration of $\BV$ by local systems.  \\[0,2cm]
(a)~We say that
the pair $(\BV, W_\bullet)$ \emph{satisfies condition $(MON)$ relative
to the divisor $Z$} if for any point $z$ of $Z$
there exists a $\BZ$-base $\FT$ of the local monodromy group around $Z'$ at $z$, adapted to local coordinates,
such that the following holds:
for any non-empty subsystem $\FT' = (T_{i_1}, \ldots, T_{i_r})$ of $\FT$,
denote by $W_\bullet^{\FT'}$
the monodromy weight filtration
of $N_{i_1} + \ldots + N_{i_r}$ on $\BV_{\FU^\circ(z)}$, where $N_j$ denotes the
logarithm of the image of $T_j$ under the action.
Then $W_\bullet^{\FT'}$ coincides with the filtration induced by $W_\bullet$
on $\BV_{\FU^\circ(z)}$.
(In particular, $W_\bullet^{\FT'}$ is independent of $\FT' \ne \emptyset$.)
\\[0.2cm]
(b)~Let $k$ be an integer. We say that the triple $(\BV, k, W_\bullet)$
\emph{satisfies condition $(MON)$ relative to the divisor $Z$} if the pair
$(\BV, W[k]_\bullet)$ satisfies condition $(MON)$ relative to the divisor $Z$.
\end{Def}

Here as usual, the shifted filtration $W[k]_\bullet$ is defined by
$W[k]_n \BV := W_{k+n} \BV$.
Condition $(MON)$ will be studied in a situation
where the local system in question underlies a variation
of Hodge structure, pure of weight $k$. In order to analyze its specialization
along $Z$, it is natural to consider the shift by $-k$ of the monodromy weight
filtration, rather than the monodromy weight filtration itself, whence
part (b) of the above definition.
Let us remark that condition $(MON)$
is very restrictive, but occurs naturally in the context
of toroidal compactifications of
Shimura varieties, as we shall see in Section~\ref{5}.
The main result
of this section reads as follows:

\begin{Thm} \label{3cY}
Assume that $Z$ is a reduced principal divisor defined by a function $g$.
Let $(\BV, k, W_\bullet)$ be a triple satisfying condition $(MON)$
relative to $Z$. Assume that $\BV$ underlies a variation of Hodge structure,
pure of weight $k$. Then the following filtrations on the perverse sheaf
$\psi_g^p \, j_* \BV$ on $Z$ coincide:
\begin{enumerate}
\item[(a)] the monodromy weight filtration (of the canonical monodromy
endomorphism), shifted by $-k$,
\item[(b)] the monodromy weight filtration relative to
$\psi_g^p \, j_* W_\bullet$.
\end{enumerate}
\end{Thm}

\begin{Cor} \label{3cZ}
Assume that $Z$ is a reduced divisor.
Let $(\BV, k, W_\bullet)$ be a triple satisfying condition $(MON)$
relative to $Z$. Assume that $\BV$ underlies a variation of Hodge structure,
pure of weight $k$. Then the following filtrations on the perverse sheaf
$Sp_Z \, j_* \BV$ on $N_{Z/X}$ coincide:
\begin{enumerate}
\item[(a)] the monodromy weight filtration, shifted by $-k$,
\item[(b)] the monodromy weight filtration relative to $Sp_Z \, j_* W_\bullet$.
\end{enumerate}
\end{Cor}

By what was said before, the theorem follows from its corollary. But we shall
prove the results in the above logical order:\\

\begin{Proofof}{Corollary~\ref{3cZ}, assuming Theorem~\ref{3cY}}
By \cite[(SP0)]{V2}, the question is local, so we can assume that
$Z$ is a principal divisor defined by $g$. By \cite[Section~3, Claim~1)]{V3}
and \cite[2nd proposition of Section~4]{V3},
the weight filtrations on $Sp_Z \, j_* \BV$
are uniquely determined by the weight filtrations on both $\psi_g^p \, j_* \BV$
and $\phi_g^p \, j_* \BV$, where $\phi_g$  denotes
the \emph{vanishing cycle functor}, and $\phi_g^p := \phi_g [-1]$.
Our claim follows
from \ref{3cY}, and from the fact that there is a
canonical isomorphism
\[
\var: \phi_g^p \, j_* \BV \isoto \psi_g^p \, j_* \BV \; ,
\]
which is compatible with the action of the $N_i$
\cite[Section~3, Claim~5)]{V3}.
\end{Proofof}

For the proof of Theorem~\ref{3cY}, we shall use two main ingredients:\\
(A) First (see Proposition~\ref{3cF}), the explicit description,
due to Galligo--Granger--Maisonobe
\cite{GGM}, of the full sub-category $(\Perv_\BC X)_{nc}$
of $\Perv_\BC X$
of perverse sheaves \emph{of normal crossing type},
when $X$ is a product
of unit disks, $U$ the corresponding product of punctured disks,
and $g$ a product of coordinates. Actually, we shall restrict ourselves to
the unipotent objects in this category.
We shall follow the presentation
of \cite[3.1]{Sa}.
Next (see Proposition~\ref{3cG}), using the above,
the explicit description, due to Saito
\cite[Thm.~3.3]{Sa}, of the functor $\psi_g^p$, together with its
monodromy endomorphism.\\
\noindent (B) The theory of nilpotent orbits, in particular,
the comparison of monodromy
weight filtrations
in nilpotent orbits of several variables, due to Cat\-ta\-ni--Kaplan
(see Theorem~\ref{3cI}). \\
\noindent We can conclude, thanks to a result of Kashiwara's
(see Proposition~\ref{3cJ}), which can be interpreted as permanence
of nilpotent orbits under $\psi_g^p$. \\

\begin{Rem} \label{3cC}
Theorem~\ref{3cY} and Corollary~\ref{3cZ} continue to hold in a larger
generality: in the definition of condition $(MON)$, we can allow
quasi-unipo\-tent local mo\-no\-dromy; in the hypotheses of \ref{3cY} and \ref{3cZ},
the divisor $Z$ need not be reduced.
We chose to add the conditions which ensure that the mo\-no\-dromy of
$\psi_g^p \, j_* \BV$, resp.\ of $Sp_Z \, j_* \BV$ is unipotent. First, this
covers the situation we shall be considering in Section~\ref{5}. Second,
restriction to unipotent objects in the explicit description of
$(\Perv_\BC X)_{nc}$ simplifies considerably the presentation of the material.
\end{Rem}

Until the actual proof of Theorem~\ref{3cY},
we shall study the situation $X = D^n$ and $U = (D^*)^n$, for
some $n \ge 1$. Put $Z_i := \{ x_i = 0 \}$, and
\[
Z_I := \bigcap_{i \in I} Z_i \; ,
\]
for $I \subset \{ 1, \ldots, n \}$.

\begin{Def} \label{3cD}
(a) Let $(\Perv_\BC X)_{nc}$
be the category
of perverse sheaves on $X$, whose characteristic
varieties are contained in the
union of the conormal bundles of $Z_I$. \\[0.2cm]
(b) Denote by $(\Perv_\BC X)_{ncu}$
the full sub-category of unipotent perverse shea\-ves, i.e., the objects $\BW$
of $(\Perv_\BC X)_{nc}$ satisfying the following: the canonical monodromy
automorphisms on $\psi_{x_i}^p (\BW)$ and $\phi_{x_i}^p (\BW)$
are unipotent, for all $i$.
\end{Def}

\begin{Def} \label{3cE}
Define the category $\Pnu$ as follows: objects are finite-dimensional
$\BC$-vector spaces $E_I$ indexed by $I \subset \{ 1, \ldots, n \}$, together
with morphisms
\begin{eqnarray*}
\can_i: E_I \longto E_{I \cup \{ i \} } & \text{for} & i \not\in I \; , \\
\var_i: E_I \longto E_{I - \{ i \} } & \text{for} & i \in I \; , \\
N_i: E_I \longto E_I & \text{for} & i \in \{ 1, \ldots, n \} \; ,
\end{eqnarray*}
subject to the following conditions: $N_i = \can_i \circ \var_i$ and
$N_i = \var_i \circ \can_i$ whenever these compositions make sense;
all $N_i$ are nilpotent;
furthermore, $A_i$ and $B_j$ commute for $i \ne j$ and $A,B \in \{ \can, \var,
N \}$ such that the composition makes sense. Morphisms in $\Pnu$ are the
morphisms of vector spaces compatible with the $\can_i$, $\var_i$, and $N_i$.
\end{Def}

We then have the following:

\begin{Prop}[Galligo--Granger--Maisonobe] \label{3cF}
(a) There is a na\-tu\-ral equivalence of categories
\[
\Psi^n: (\Perv_\BC X)_{ncu} \isoto \Pnu \; .
\]
It is defined by associating to $\BW \in (\Perv_\BC X)_{ncu}$ the data
$(E_I)_I$, where
\[
E_I:= \Psi_{x_1,I} \circ \Psi_{x_2,I} \circ \ldots \circ \Psi_{x_n,I} (\BW) \; ,
\]
with
\[
\Psi_{x_i,I}:=
\begin{cases}
\psi_{x_i}^p \; , \ & \text{ if } i \not\in I \; ,\\
\phi_{x_i}^p \; , \ & \text{ if } i \in I \; .
\end{cases}
\]
The morphisms $\can_i$, $\var_i$ and $N_i$ are the ones naturally associated to
$\psi_{x_i}^p$ and $\phi_{x_i}^p$. \\[0.2cm]
(b) There is a natural quasi-inverse $(\Psi^n)^{-1}$ of $\Psi^n$.
\end{Prop}

Let us illustrate the effect of the functor $\Psi^n$ for $n=1$, and for a perverse
sheaf on $X = D$ of the form $j_* \, \FF$, for a unipotent
local system $\FF$ on $U = D^*$. We identify $\FF$ with a vector space $H$, together
with a unipotent automorphism $T$.
We then have $E_\emptyset = E_{\{1\}} = H$,
the morphism $\var = \var_1: E_{\{1\}} \to E_\emptyset$ is the identity on $H$, and $\can = \can_1: E_\emptyset \to E_{\{1\}}$ is
the logarithm of $T$.  \\

\begin{Proofof}{Proposition~\ref{3cF}}
Our claim is in fact a particular case of \cite[Thm.~IV.3]{GGM}. There, the
condition on (1)~unipotency of the perverse sheaves (see \ref{3cD}~(b))
is dropped,
(2)~nilpotency of the endomorphisms $N_i$ (see \ref{3cE}) is replaced
by ``the sums $\id + N_i$ are invertible''.
By splitting the objects into generalized eigenspaces with respect to
the commuting operators $\id + N_i$, one gets the description of \cite[3.1]{Sa}.
Our situation corresponds to the (multiple) eigenvalues $(1,1,\ldots,1)$.
In the description of \cite[3.1]{Sa}, this means that the components $E_I^\nu$
are trivial whenever $v \in (\BC / \BZ)^n$ is unequal to zero.
\end{Proofof}

Now fix $m \in \{ 1, \ldots, n \}$, set $g := \prod_{i=1}^m x_i$, and consider
the principal reduced divisor $Z = \cup_{i=1}^m Z_i$ defined by $g$. Using the
fact that the direct image of a closed embedding is fully faithful, we
may view $\psi_g^p$ as a functor from $\Perv_\BC X$ to itself. We then have:

\begin{Prop}[Saito] \label{3cG}
For any subset $I$ of $\{ 1, \ldots, n \}$, denote by $I_\le$ the intersection
$I \cap \{ 1, \ldots, m \}$, and by $I_>$ the complement of $I_\le$
in $I$. \\[0.2cm]
(a) The functor $\psi_g^p$ respects the sub-category $(\Perv_\BC X)_{ncu}$
of $\Perv_\BC X$.\\[0.2cm]
(b) The composition
\[
\Psi^n \circ \psi_g^p \circ (\Psi^n)^{-1}: \Pnu \longto \Pnu
\]
is given as follows: let
\[
\BE = (E_I, \can_i, \var_i, N_i)_{I,i}
\]
be an object of $\Pnu$. Then
\[
\Psi^n \circ \psi_g^p \circ (\Psi^n)^{-1} (\BE) =
(\tE_I, \tcan_i, \tvar_i, \tN_i)_{I,i} \; ,
\]
with
\[
\tE_I := \coker \left( \prod_{i \in I_\le} (N_i - N) \; \tei \;
                                                   E_{I_>}[N] \right)\; ,
\]
where we define $E_{I_>}[N]$ as the tensor product of $E_{I_>}$ and the
polynomial ring $\BC[N]$ in one variable $N$. The variable acts on $\BC[N]$
by multiplication.
The actions of $N_i$ and of $N$ on $E_{I_>}[N]$
are the ones induced by the tensor product structure.
In particular, the endomorphisms $\prod_{i \in I_\le} (N_i - N)$ of
$E_{I_>}[N]$ are injective, so we may identify their cokernels with their
mapping cones. The morphisms $\tcan_i$, $\tvar_i$, and $\tN_i$
are given as morphisms of complexes concentrated in two degrees:
\[
\tcan_i :=
\begin{cases}
(\id , N_i - N) \; , \ & \text{ if } i \not\in I \; , \; 1 \le i \le m \; , \\
(\can_i , \can_i ) \; , \ & \text{ if } i \not\in I \; , \; i > m \; ,
\end{cases}
\]
\[
\tvar_i :=
\begin{cases}
(N_i - N , \id) \; , \ & \text{ if } i \in I_\le \; , \\
(\var_i , \var_i ) \; , \ & \text{ if } i \in I_> \; ,
\end{cases}
\]
\[
\tN_i :=
\begin{cases}
(N_i - N , N_i - N) \; , \ & \text{ if } i \in I_\le \; , \\
(N_i , N_i ) \; , \ & \text{ if } i \in I_> \; .
\end{cases}
\]
(c) For $\BE \in \Pnu$, the canonical monodromy endomorphism on
\[
\Psi^n \circ \psi_g^p \circ (\Psi^n)^{-1} (\BE) =
(\tE_I, \tcan_i, \tvar_i, \tN_i)_{I,i} \; ,
\]
in the description
of (b), is given by the endomorphism $\tN := (N , N)$
on all components $\tE_I$.
\end{Prop}

\begin{Proof}
This is part of the information provided by \cite[Thm.~3.3]{Sa}. There,
the specialization $Sp_Z$ is described in terms of the categories $\Pnu$
and $\BP (n+1)_u$. By the last line of \cite[Thm.~3.3]{Sa}, in order
to read off $\psi_g^p$ from the given description, one has to restrict to
the components ``$0 \not\in I$''. Note that Saito admits
quasi-unipotent objects. Thanks to reducedness of $Z$, the
condition ``$E_I^\nu = 0$
whenever $v$ is unequal to zero'' is respected by $Sp_Z$,
hence by $\psi_g^p$. This shows parts~(a) and (b) of our claim.
As for (c), observe that by \cite[Thm.~3.3]{Sa},
$\tN$ occurs in the explicit description of $Sp_Z$.
More precisely, it is the restriction to
the components ``$0 \not\in I$'' of the collection of the $n+1$st
nilpotent endomorphisms $M_0$ of the components of
\[
\Psi^{n+1} \circ Sp_Z \circ (\Psi^n)^{-1} (\BE) =
(F_I, \can_i, \var_i, M_i)_{I,0 \le i \le n} \; .
\]
Our claim follows thus from \cite[(SP6)]{V2}.
\end{Proof}

In order to prove Theorem~\ref{3cY}, one is thus naturally led to study
mo\-no\-dromy weight filtrations on objects of the form
\[
\tE_I := \coker \left( \prod_{i \in I_\le} (N_i - N) \; \tei \;
                                                   E_{I_>}[N] \right)\; ,
\]
for certain objects $\BE = (E_I, \can_i, \var_i, N_i)_{I,i}$ of $\Pnu$.
The result we want to use requires an additional structure on $\BE$.
Recall the notion of nilpotent orbits of a weight $k \in \BZ$
and dimension $n \in \BN$
(e.g.\ \cite[(3.1)]{CK}; cmp.\ also \cite[4.1]{K}).
For such objects $\BH$, we shall use the notation
\[
\BH = ((H,F^\bullet,W_\bullet), N_i \, (1 \le i \le n), S) \; .
\]
As for the nature of the components of $\BH$,
note in particular
that $H$ is a finite-dimensional $\BC$-vector space, with finite descending,
resp.\ ascending filtrations $F^\bullet$ and $W_\bullet$, $S$ is a sesquilinear
form on $H$, the $N_i$ are mutually commuting nilpotent endomorphisms,
and $W_\bullet$ is the monodromy weight filtration of the sum
$\sum_{i=1}^n N_i$, shifted by $-k$. \\

The main motivation for this concept stems from Schmid's Nilpotent Orbit
Theorem, which we shall use in the following form:

\begin{Thm}[Schmid] \label{3cH}
Let $\BV$ be a local system on $U$, which underlies a variation of Hodge
structure, pure of weight $k$. Write
\[
\Psi^n (j_* \BV) =: \BE = (E_I, \can_i, \var_i, N_i)_{I,i} \: .
\]
Then for any subset $I$ of $\{ 1, \ldots , n \}$, the data
\[
(E_I, N_i \, (1 \le i \le n))
\]
underly a nilpotent orbit of weight $k$.
\end{Thm}

\begin{Proof}
$\BV$ is given by a vector space $H$, together with commuting mo\-no\-dromy
automorphisms $T_i$, $1 \le i \le n$. We then have $E_\emptyset = H$, and
the $N_i$ are the logarithms of the unipotent parts of the $T_i$.
Now apply \cite[Thm.~(4.12)]{Sch} to show the claim for $I=\emptyset$.
But all the other components of $\Psi^n (j_* \BV)$ are
isomorpic to $E_\emptyset$ \cite[Section~3, Claim~5)]{V3}.
\end{Proof}

On nilpotent orbits, comparison of
monodromy weight filtrations is possible thanks
to the following result:

\begin{Thm}[Cattani--Kaplan] \label{3cI}
Let
\[
\BH = ((H,F^\bullet,W_\bullet), N_i \, (1 \le i \le n), S)
\]
be a nilpotent
orbit of weight $k$ and dimension $n$, and $I_1$ and $I_2$ two disjoint
subsets of $\{ 1, \ldots , n \}$.
Denote by $W_\bullet^{I_1}$ the monodromy weight filtration of
$\sum_{i \in I_1} N_i$, shifted by $-k$.
Then the following filtrations on $H$ coincide:
\begin{enumerate}
\item[(a)] the monodromy weight filtration of
$\sum_{i \in I_1 \cup I_2} N_i$, shifted by $-k$,
\item[(b)] the monodromy weight filtration of $\sum_{i \in I_2} N_i$
relative to $W_\bullet^{I_1}$.
\end{enumerate}
\end{Thm}

\begin{Proof} This is the content of
\cite[Thm.~(3.3)]{CK}. Note that the original statement of loc.\ cit.\ is
misprinted; the correct version can be found in \cite[Prop.~(4.72)]{CKS}.
\end{Proof}

\begin{Prop}[Kashiwara] \label{3cJ}
Let
\[
\BH = ((H,F^\bullet,W_\bullet), N_i \, (1 \le i \le n), S)
\]
be a nilpotent
orbit of weight $k$ and dimension $n$, and
$\emptyset \ne I \subset \{ 1, \ldots, n \}$. Set
\[
\tH_I := \coker \left( \prod_{i \in I} (N_i - N) \; \tei \;
                                                   H[N] \right)\; .
\]
(a) The vector space
$\tH_I$ underlies in a natural way a nilpotent orbit of weight $k+1- | I |$
and dimension $n+1$
\[
\tBH = ((\tH_I,F^\bullet,M_\bullet), N, N_i \, (1 \le i \le n), \tS) \; .
\]
In particular, $M_\bullet$ is the monodromy weight filtration of the sum
$N + \sum_{i=1}^n N_i$, shifted by $-(k+1-  | I |)$. \\[0.2cm]
(b) The filtration
$M_\bullet$ coincides with the monodromy weight filtration of
$N$, shifted by $-(k+1- | I |)$.
\end{Prop}

\begin{Proof}
Part~(a) is contained in \cite[Prop.~3.19]{Sa}. In order to see that
(b) holds, one has to look at Kashiwara's proof of loc.\ cit.
\cite[(A.3.1), A.4]{Sa}.
\end{Proof}

Combining the two preceding results, we get:

\begin{Cor} \label{3cK}
Keep the assumptions of Proposition~\ref{3cJ}, and denote by $W_\bullet$ the
filtration on the vector space $\tH_I$ induced by the filtration $W_\bullet$
on $H$. (Note that the functor $H \mapsto \tH_I$ is exact.) Then
the following filtrations on $\tH_I$ coincide:
\begin{enumerate}
\item[(a)] the monodromy weight filtration of
$N$, shifted by $-k$,
\item[(b)] the monodromy weight filtration of $N$ relative to $W_\bullet$.
\end{enumerate}
\end{Cor}

Finally, we can show the main result of this section:\\

\begin{Proofof}{Theorem~\ref{3cY}}
Since the question is local, we may assume that we are in the situation
discussed in \ref{3cD}--\ref{3cK}. Propositions~\ref{3cF} and \ref{3cG}
tell us that we need to compare monodromy weight filtrations on
the
\[
\tE_I := \coker \left( \prod_{i \in I_\le} (N_i - N) \; \tei \;
                                                   E_{I_>}[N] \right)\; ,
\]
for $\Psi^n (j_* \BV) =: \BE = (E_I, \can_i, \var_i, N_i)_{I,i}$.
By Theorem~\ref{3cH}, the $E_{I_>}$ underly nilpotent orbits of
weight $k$. We omit the $N_i$ with $i > m$, and consider $E_{I_>}$
as nilpotent orbit of dimension $m$. Now apply Corollary~\ref{3cK}, with
$n$ replaced by $m$, and $I$ replaced by $I_\le = \{ i_1 , \ldots , i_r \}$.
Thanks to condition $(MON)$, the filtration $W_\bullet$
in \ref{3cK} coincides with the one induced by the filtration $W_\bullet$
of $\BV$.
\end{Proofof}


\bigskip
%
%

\section{Strata in toroidal compactifications}
\label{4}

In order to prepare the proof of Theorem~\ref{2D}, to be given in
Section~\ref{5}, we need to discuss the geometry of
\emph{toroidal compactifications} $\MKS$ of $M^K$.
We keep the notations and hypotheses of Section~\ref{2}. In particular,
the subgroup $K \subset G(\BA_f)$ is neat, and $(G,\FH)$ satisfies $(+)$.
Choose a \emph{$K$-admissible cone decomposition} $\FS$ satisfying the
conditions of \cite[(3.9)]{P2}. In particular, all cones occurring in
$\FS$ are \emph{smooth}, and the decomposition is \emph{complete}.\\

Let us denote by $\MKS := M^K (G,\FH,\FS)$ the toroidal compactification
associated to $\FS$. It is a smooth projective scheme over $\BC$,
which in a natural way contains $M^K$ as an open sub-scheme.
The complement
is a union of smooth divisors with normal crossings. The identity on
$M^K$ extends uniquely to a surjective morphism
\[
p = p_\FS: \MKS \longonto (M^K)^* \; .
\]
The inverse images under $p$ of the strata described in Section~\ref{1}
form a stratification of $\MKS$.
We follow \cite[(3.10)]{P2} for the description of these
inverse images: as usual,
fix a proper boundary component $(P_1,\FX_1)$ of $(G,\FH)$, and an element
$g \in G(\BA_f)$. To the given data, the following are
canonically associated:
\begin{enumerate}
\item[(i)] an Abelian scheme $A \to \Mp$, and an $A$-torsor $B \to \Mp$,
\item[(ii)] a torus $T$, and a $T$-torsor $X \to B$,
\item[(iii)] a rational partial polyhedral decomposition of $Y_*(T)_\BR$
($Y_*(T) :=$ the cocharacter group of $T$), again denoted by $\FS$, and a
non-empty subset $\FT \subset \FS$,
\item[(iv)] an action of $H_Q$ on $B$, $X$, and $T$.
\end{enumerate}
These objects satisfy the following properties:
\begin{enumerate}
\item[(A)] the $H_Q$-action is equivariant with respect to the group and
torsor structures and stabilizes $\FS$ and $\FT$,
\item[(B)] the subgroup $P_1(\BQ)$ of $H_Q$ acts trivially on $B$, $X$,
and $T$,
\item[(C)] the group $\Delta_1 = H_Q/P_1(\BQ)$ acts freely on $\FT$,
\item[(D)] the pair $(\FS,\FT)$ satisfies conditions~\cite[(2.3.1--3)]{P2}
(see below).
\end{enumerate}
Consider the relative torus embedding $X \into X (\FS)$.
Condition~\cite[(2.3.1)]{P2} is equivalent to
saying that $\FT$ defines a closed sub-scheme $Z$ of $X (\FS)$.
In order to state the other two conditions, define
\[
D := \bigcup_{\sigma \in \FT} \sigma^\circ \; ,
\]
where for each cone $\sigma$ we denote by $\sigma^\circ$ the topological
interior of $\sigma$ inside the linear subspace of $Y_*(T)_\BR$ generated by
$\sigma$. The subset $D$ of $Y_*(T)_\BR$ is endowed with the induced
topology. Condition~\cite[(2.3.2)]{P2} says that every point of $D$ admits
a neighbourhood $U$ such that $U \cap \sigma \ne \emptyset$ for only a
finite number of $\sigma \in \FT$. Condition~\cite[(2.3.3)]{P2} states that
$D$ is contractible.\\

By (C),
the induced action of $\Delta_1$ on $Z$ is free and proper in the
sense of \cite[(1.7)]{P2}. The geometric quotient
$\MoS$ exists and is canonically
isomorphic to the inverse image of $\Mo$ under $p$.
Furthermore, the analytic space $\MoS (\BC)$ is the quotient of $Z (\BC)$
by $\Delta_1$ in the analytic category.
We summarize the situation
by the following diagram:
\[
\vcenter{\xymatrix@R-10pt{
M^K \ar@{^{ (}->}[r]^-{j_\FS} \ar@{=}[d] &
\MKS \ar@{<-^{ )}}[r]^-{i_\FS} \ar@{->>}[d]^p &
\MoS = \Delta_1 \backslash Z \ar@{<<-}[r]^-{\tilde{q}} \ar@{->>}[d]^p &
Z \ar@{->>}[d]^{\tilde{p}} \\
M^K \ar@{^{ (}->}[r]^-{j} &
(M^K)^* \ar@{<-^{ )}}[r]^-{i} &
\Mo = \Delta \backslash \Mp \ar@{<<-}[r]^-q &
\Mp \\}}
\]
The left and the middle square are Cartesian, and the maps $p$ are proper.\\

It will be necessary to consider a refinement of the
stratification of $\MKS$. The induced stratification of $Z$
is the natural one given by $\FT$.
For any cone $\sigma \in \FT$, denote by
\[
\isc: \Msc \longinto Z
\]
the immersion of the corresponding stratum into $Z$, and by
\[
\is: \Ms \longinto Z
\]
the immersion of its closure.
In the same way, we shall write
\[
\istc: \Mstc \longinto \MoS \longinto \MKS
\]
and
\[
\ist: \Mst \longinto \MoS \longinto \MKS
\]
for the respective immersions into $\MoS$, or into $\MKS$.
These immersions
are indexed by the quotient $\tilde{\FT} := \Delta_1 \backslash \FT$.
Note that $\Mst$ is closed in $\MoS$, but in
general not in $\MKS$. \\

In order to describe the situation on the level of the underlying
analytic spaces, let us connect
the present notation (which is that of \cite[(3.10)]{P2})
to the one of \cite{P1}.
Consider the factorization of $\pi: (P_1,\FX_1) \longto (G_1,\FH_1)$
corresponding to the weight filtration of the unipotent radical $W_1$:
\[
\vcenter{\xymatrix@R-10pt{
(P_1,\FX_1) \ar@{->>}[r]^-{\pi_t} \ar@/_2pc/[rr]_\pi&
(P'_1,\FX'_1) := (P_1,\FX_1)/U_1 \ar@{->>}[r]^-{\pi_a} &
(G_1,\FH_1)
\\}}
\]
where $U_1$ denotes the weight $-2$ part of $W_1$ \cite[Def.~2.1~(v)]{P1}.
On the level of Shimura varieties, the picture looks as follows:
\[
\vcenter{\xymatrix@R-10pt{
\MKo = \MKo (P_1,\FX_1) \ar@{->>}[r]^-{\pi_t} \ar@/_2pc/[rr]_\pi &
M^{\pi_t (K_1)} := M^{\pi_t (K_1)} (P'_1,\FX'_1) \ar@{->>}[r]^-{\pi_a} &
\Mp
\\}}
\]
By \cite[3.12--3.22~(a)]{P1}, $\pi_a$ is in a natural way a torsor under an
Abelian scheme, while $\pi_t$ is a torsor under a torus.
In fact, we have (i)~$B = M^{\pi_t (K_1)}$, and (ii)~$X = \MKo$.
Furthermore, the action~(iv) of $H_Q$ on $B$, $X$ and $T$ is induced by the
natural action of $H_Q$ on the Shimura data
involved in the above factorization
of $\pi$. Since $P_1(\BQ)$ acts trivially on the associated Shimura
varieties, this explains property~(B). \\

The map $\tilde{p}: Z \to \Mp$ thus factors through $\pi_a$. \cite[6.13]{P1}
contains the definition of a $K_1$-admissible smooth
cone decomposition $\FS_1^0$
canonically associated to $(P_1,\FX_1)$ and $g$. It is \emph{concentrated in
the unipotent fibre} \cite[6.5~(d)]{P1}, and thus defines a smooth torus
embedding $j_1: \MKo \into \MKoS$ over $M^{\pi_t (K_1)}$.
In fact, we have
$X(\FS) = \MKoS$. Furthermore \cite[6.13]{P1}, there is a
closed analytic subset $\dU := \dU(P_1,\FX_1,g)$
of $\MKoS (\BC)$ canonically associated to our data.
The proof of \cite[Prop.~6.21]{P1} shows that $\dU = Z (\BC)$.
In fact, the projection $\tilde{q}: Z \to \Delta_1 \backslash Z$
corresponds to the quotient map
\[
\dU \longonto \Delta_1 \backslash \dU
\]
of \cite[7.3]{P1}.

\begin{Prop} \label{4B}
The morphism $\tilde{q}: Z \longto \MoS$ induces an isomorphism
\[
\Ms \isoto \Mst
\]
for any $\sigma \in \FT$. In particular, it
induces an isomorphism on every irreducible component of Z.
\end{Prop}

\begin{Proof}
By \cite[Cor.~7.17~(a)]{P1}, the morphism
\[
\tilde{q}: \Ms \longonto \Mst
\]
identifies $\Mst$ with the quotient of $\Ms$ by a certain
subgroup $\Stab_{\Delta_1} ([\sigma])$ of $\Delta_1$. (The hypotheses of
loc.~cit.\ are satisfied because they are implied by the conditions
of \cite[(3.9)]{P2}, which we assume throughout.)
By \cite[Lemma~1.7]{W2}, condition $(+)$ and neatness of $K$ imply that
\[
\Stab_{\Delta_1} ([\sigma]) = 1 \; .
\]
\end{Proof}

Consider the diagram
\[
\vcenter{\xymatrix@R-10pt{
Z \ar@{^{ (}->}[rr]^-{i_1} \ar@{->>}[d]_{\tilde{q}} &&
X(\FS) = \MKoS \\
\MoS \ar@{^{ (}->}[r]^-{i_\FS} &
\MKS & \\}}
\]
By \cite[6.13]{P1}, there is an open neighbourhood
$\FU := \bar{\FU}(P_1,\FX_1,g)$
of $Z (\BC)$ in $\MKoS (\BC)$, and a natural extension of
$\tilde{q}$ to $\FU$. It will equally be denoted by $\tilde{q}$:
\[
\vcenter{\xymatrix@R-10pt{
Z (\BC) \ar@{^{ (}->}[r]^-{i_1} \ar@{->>}[d]_{\tilde{q}} &
\FU \ar@{^{ (}->}[r] \ar@{->>}[d]_{\tilde{q}} &
\MKoS (\BC)
\\
\MoS (\BC) \ar@{^{ (}->}[r]^-{i_\FS} &
\MKS (\BC) & \\}}
\]
Furthermore, the open subset $\FU$ is stable under $\Delta_1$.
By \cite[Prop.~1.9]{W2}, the map $\tilde{q}$ (which in loc.~cit.\
was denoted by $f$) is open, and we have the equality
\[
\tilde{q}^{-1} (M^K (\BC)) = \FU \cap \MKo (\BC) \; .
\]
Furthermore \cite[Thm.~1.11~(i)]{W2},
$\tilde{q}$ is locally biholomorphic near $Z$. It thus
induces an isomorphism
between the quotient of the formal analytic
completion of $\MKoS (\BC)$ along $Z (\BC)$
by the free action of
$\Delta_1$, and the formal analytic
completion of $\MKS (\BC)$ along $\MoS (\BC)$.
According to \cite[p.~224]{P2}, we have:

\begin{Prop} \label{4D}
This isomorphism is algebraic in the following sense: \\[0.2cm]
\noindent (a) The action of $\Delta_1$ on the formal
completion $\FF = \FF_{Z / \MKoS}$
of $\MKoS$ along $Z$ is free and proper. The geometric quotient
$\Delta_1 \backslash \FF$ exists.
Furthermore, the analytic space
$(\Delta_1 \backslash \FF) (\BC)$ is the quotient of
$\FF (\BC)$ by $\Delta_1$ in the analytic category. \\[0.2cm]
\noindent (b) $\tilde{q}$ induces an isomorphism
between $\Delta_1 \backslash \FF$ and the formal
completion of $\MKS$ along $\MoS$.
\end{Prop}

\forget{\? {\tt It would be more satisfactory to give a proof of this.
The problem is that \cite[Lemma~9.38]{P1} cannot be applied directly
since $\MoS$ is not proper...} \? \\
}
In fact, $\tilde{q}$ descends to the reflex field of our Shimura
varieties. In Section~\ref{5}, the following
consequence of Proposition~\ref{4D} will be needed:

\begin{Cor} \label{4E}
The map $\tilde{q}$ induces an isomorphism
\[
\Delta_1 \backslash N_{Z / \MKoS}
\isoto N_{\MoS / \MKS}
\]
between the quotient of the normal cone of $Z$ in $\MKoS$
by the free and proper
action of $\Delta_1$, and the normal cone of $\MoS$ in $\MKS$.
\end{Cor}


\bigskip
%
%

\section{Proof of the main result}
\label{5}

Recall the situation
considered in Section~\ref{4}:
\[
\vcenter{\xymatrix@R-10pt{
M^K \ar@{^{ (}->}[r]^-{j_\FS} \ar@{=}[d] &
\MKS \ar@{<-^{ )}}[r]^-{i_\FS} \ar@{->>}[d]^p &
\MoS = \Delta_1 \backslash Z \ar@{<<-}[r]^-{\tilde{q}} \ar@{->>}[d]^p &
Z \ar@{->>}[d]^{\tilde{p}} \\
M^K \ar@{^{ (}->}[r]^-{j} &
(M^K)^* \ar@{<-^{ )}}[r]^-{i} &
\Mo = \Delta \backslash \Mp \ar@{<<-}[r]^-q &
\Mp \\}}
\]
Proper base change \cite[(4.4.3)]{Sa} yields the following:

\begin{Prop} \label{5K}
There is a canonical isomorphism of functors
\[
i^* j_* \cong p_* \, i_\FS^* \, {j_\FS}_*:
                           D^b (\MHM_F M^K) \longto D^b (\MHM_F \Mo) \; .
\]
\end{Prop}

We are thus led to study the inverse image
\[
i_\FS^*: D^b (\MHM_F \MKS) \longto D^b (\MHM_F \MoS) \; .
\]
According to Corollary~\ref{4E}, the normal cone $N_{\MoS / \MKS}$ is
canonically isomorphic to
the quotient of the normal cone $N_{Z / \MKoS}$
by the free and proper action of $\Delta_1$. Using Corollary~\ref{5E},
we make the following identifications:
\[
D^b (\MHM_F \MoS) = D^b (\Delta_1 \text{-} \MHM_F Z) \; ,
\]
\[
D^b (\MHM_F N_{\MoS / \MKS}) =
D^b (\Delta_1 \text{-} \MHM_F N_{Z / \MKoS}) \; .
\]
Since the action of $\Delta_1$ on $N_{Z / \MKoS}$ respects
the natural inclusion of $Z$, we can think of the inverse image
\[
i_0^*: D^b (\MHM_F N_{\MoS / \MKS}) \longto D^b (\MHM_F \MoS)
\]
as the $\Delta_1$-equivariant inverse image
\[
i_0^*: D^b (\Delta_1 \text{-} \MHM_F N_{Z / \MKoS}) \longto
D^b (\Delta_1 \text{-} \MHM_F Z) \; .
\]
Recall the specialization functor
\[
Sp_{\MoS}: D^b (\MHM_F \MKS) \longto D^b (\MHM_F N_{\MoS / \MKS}) \; .
\]
According to \cite[2.30]{Sa}, we have:

\begin{Prop} \label{5L}
There is a canonical isomorphism of functors
\[
i_\FS^* \cong i_0^* \, Sp_{\MoS} : D^b (\MHM_F \MKS) \longto D^b (\MHM_F {\MoS}) \;.
\]
\end{Prop}

We summarize the situation by the following commutative diagram:
\[
\vcenter{\xymatrix@R-10pt{
& D^b (\Delta_1 \text{-} \MHM_F N_{Z / \MKoS})
                          \ar@{=}[d] \ar@/^10pc/[ddd]^{i_0^*}   \\
D^b (\MHM_F \MKS) \ar[r]^-{Sp_{\MoS}} \ar[dr]_{i_\FS^*} &
  D^b (\MHM_F N_{\MoS / \MKS}) \ar[d]^{i_0^*}      \\
& D^b (\MHM_F \MoS) \ar@{=}[d]    \\
& D^b (\Delta_1 \text{-} \MHM_F Z)
\\}}
\]
Recall the open immersion $j_1: \MKo \into \MKoS$
introduced in Section~\ref{4}.
It is a smooth relative torus embedding, hence in particular affine.
This allows to define the exact functor
\[
{j_1}_*: \MHM_F \MKo \longto \MHM_F \MKoS
\]
even though $\MKoS$ is only locally
of finite type: cover $\MKoS$ by open affines,
use exactness of the direct image of the restriction of $j_1$ to each
such affine \cite[4.2.11]{Sa}, and glue. The same technique allows to define
the specialization functor
\[
Sp_Z: \MHM_F \MKoS \longto \MHM_F N_{Z / \MKoS} \; .
\]
Because of the functorial behaviour of ${j_1}_*$ and $Sp_Z$,
these functors admit $\Delta_1$-equivariant versions. Since they are
exact, they induce functors on the level of bounded derived categories.

\begin{Prop} \label{5M}
There is a natural commutative diagram
\[
\vcenter{\xymatrix@R-10pt{
D^b (\Rep_F G) \ar[r]^-{\mu_K} \ar[d]_{\Res^G_Q} &
                                  D^b (\MHM_F M^K) \ar[dddd]^{{j_\FS}_*} \\
D^b (\Rep_F Q) \ar[d] & \\
D^b (\Rep_F P_1 , H_Q) \ar[d]_{\mu_{K_1}} & \\
D^b (\Delta_1 \text{-} \MHM_F \MKo) \ar[d]_{{j_1}_*} & \\
D^b (\Delta_1 \text{-} \MHM_F \MKoS) \ar[d]_{Sp_Z} &
                                  D^b (\MHM_F \MKS) \ar[d]^{Sp_{\MoS}} \\
D^b (\Delta_1 \text{-} \MHM_F N_{Z / \MKoS}) \ar@{=}[r] &
                                  D^b (\MHM_F N_{\MoS / \MKS})
\\}}
\]
\end{Prop}

\begin{Rem}
This result implies a comparison isomorphism on the
level of singular cohomology, which is already known.
In fact, it can be seen to be equivalent to \cite[Prop.~(5.6.12)]{HZ1}.
\forget{
The composition of
$Sp_{\MoS} \, {j_\FS}_* \circ \mu_K$ and $i_\FS^*$ is calculated
in \cite[Thm.~4.3.7~(i)]{HZ2} \? {\tt in fact -- it seems that our
main result is implied by this...} \? ; indeed, as one sees from the proof
of loc.\ cit., one of the vital ingredients is an analytic version of the
above comparison statement \cite[proof of Prop.~4.3.10]{HZ2}.
However, the proof of this statement is incomplete: the
uniqueness property of \cite[Prop.~2.11]{Sa} applies to Hodge modules of
the form $j'_* \BW$, where $j'$ denotes the open immersion of the smooth
part of $N_{\MoS / \MKS}$ into $N_{\MoS / \MKS}$. But the objects in
question are in general not of this form, unless $\MoS$ consists of a single
stratum. Indeed, this phenomenon necessitated the whole of our
Section~\ref{3c}.}
\end{Rem}

\begin{Proofof}{Proposition~\ref{5M}}
Since all the functors in the diagram are exact
on the level of Abelian categories, it suffices to show the
result for objects $\BV$ of $\Rep_F G$.
Recall from Corollary~\ref{4E}
that the isomorphism between $\Delta_1 \backslash N_{Z / \MKoS}$ and
$N_{\MoS / \MKS}$ is induced by the analytic map
\[
\tilde{q}: \FU \longto \MKS (\BC) \; ,
\]
which is a local analytic isomorphism near $Z$.
The restriction to the pre-image of $M^K (\BC)$ of $\tilde{q}$ looks as follows
(see \cite[6.10]{P1}, or the proof of \cite[Prop.~2.1]{W2}): we have
\[
\tilde{q}^{-1} (M^K (\BC)) = P_1 (\BQ) \backslash
                                (\FX^+ \times P_1 (\BA_f) / K_1 ) \; ,
\]
for a certain complex manifold $\FX^+$, which is open in both $\FH$
and $\FX_1$. On $\tilde{q}^{-1} (M^K (\BC))$,
the map $\tilde{q}$ is given by
\[
[(x,p_1)] \longmapsto [(x,p_1 g)] \in G (\BQ) \backslash
                                (\FH \times G (\BA_f) / K ) = M^K (\BC) \; .
\]
It follows that the local system
$\tilde{q}^{-1} \circ \mu_{K,\topp} (\BV)$ is the restriction to
$\tilde{q}^{-1} (M^K (\BC))$ of the local system
$\mu_{K_1,\topp} (\Res^G_{P_1} \BV)$, and
that the natural action of $\Delta_1$ on
$\tilde{q}^{-1} \circ \mu_{K,\topp} (\BV)$ corresponds to the
action of $H_Q$ on $\Res^G_Q \BV$. Since the
topological version of specialization can be computed locally
\cite[(SP0)]{V2},
we thus obtain the desired comparison
result on the level of perverse sheaves.
It remains to show that this isomorphism,
call it $\alpha$, respects the weight and Hodge
filtrations.

Denote by $\BV_G$ and $\BV_Q$ the two $\Delta_1$-equivariant
variations on the open subset $\tilde{q}^{-1} (M^K (\BC))$ of
$M^{K_1} (\BC)$ obtained by restricting $\mu_K (\BV)$ and
$\mu_{K_1} (\Res^G_Q \BV)$, respectively.
By \cite[Prop.~4.12]{P1}, the Hodge filtrations on $\BV_G$ and $\BV_Q$
coincide.
By the proof of \cite[Thm.~3.27]{Sa}, the Hodge filtrations of the mixed
Hodge modules ${j_1}_* \BV_G$ and ${j_1}_* \BV_Q$
depend only on the Hodge filtrations
of $\BV_G$ and $\BV_Q$ respectively. Therefore, they coincide as well.
By definition of the functor $Sp_Z$ (see in particular \cite[2.30]{Sa}
and \cite[2.3]{Sa}),
the Hodge filtrations of $Sp_Z {j_1}_* \BV_G$ and $Sp_Z {j_1}_* \BV_Q$
depend only on the Hodge filtrations of
${j_1}_* \BV_G$ and ${j_1}_* \BV_Q$, respectively.
They are therefore respected by $\alpha$.

It remains to compare the weight filtrations of
$Sp_Z {j_1}_* \BV_G$ and $Sp_Z {j_1}_* \BV_Q$.
Recall the \emph{barycentric subdivision}
$\FS'$ of $\FS$ (e.g., \cite[5.24]{P1}). By the proof of \cite[Prop.~9.20]{P1},
the cone decomposition $\FS'$ still satisfies the
conditions of \cite[(3.9)]{P2}. The refinement induces a projective
and surjective morphism
\[
M^K (G,\FH,\FS') \longto \MKS = M^K (G,\FH,\FS) \; ,
\]
and the pre-image $Z'$ of $Z$ is a divisor (with normal crossings).
Now recall the definition of $Sp_Z$
via the nearby cycle functor \cite[2.30]{Sa}. Apply projective cohomological
base change for the latter \cite[Thm.~2.14]{Sa}, and the fact that in our
situation, the cohomology objects are trivial in degree non-zero. This shows
that without loss of generality, we may assume that $Z$
is a divisor with normal crossings.

Because of the semi-simplicity of $\Rep_F G$, we may also assume that $\BV$ is
pure of weight $k$ (say). Via $\alpha$, we view the local system
underlying $\BV_G$
as being equipped with the filtration $W_\bullet$ (coming from $\BV_Q$).
By definition of the weight filtration on $Sp_Z {j_1}_* \BV_G$
(see \cite[2.3]{Sa}), it remains to show that
the following coincide:
\begin{enumerate}
\item[(a)] the monodromy weight filtration on $Sp_Z {j_1}_* \BV_G$,
shifted by $-k$,
\item[(b)] the monodromy weight filtration on $Sp_Z {j_1}_* \BV_G$
relative to $Sp_Z \, {j_1}_* W_\bullet$.
\end{enumerate}
By \cite[Prop.~1.3]{W2}, the triple $(\BV_G, k, W_\bullet)$ satisfies condition
$(MON)$ relative to the divisor $Z$. Our claim follows
thus from Corollary~\ref{3cZ}.
\end{Proofof}

\begin{Rem}
We use the opportunity to point out a minor error in \cite{W2}. The proof
of loc.\ cit., Prop.~1.3 relies on loc.\ cit., Lemma~1.2, which is not
correctly stated: the claim ``$\imm(\iota x) \in U_1(\BR)(-1)$''
should be replaced by
``$\imm(\iota x) - \imm(x) \in U_1(\BR)(-1)$''.
As a consequence, the proof of loc.\ cit., Prop.~1.3 (but not its statement)
has to be slightly modified: in line~12 of page~328,
replace ``maps $u_0$ to $\pm \frac{1}{2 \pi i} T$'' by
``maps $u_0$ to $\pm \frac{1}{2 \pi i} T \mod U$'', where $U$ denotes the
weight $-2$ part of the unipotent radical of the group $P$. Since $U$ acts
trivially on the weight-graded parts of any representation $\BV$ of $P$,
the rest of the proof remains unchanged.
\end{Rem}

By the preceding results,
we have to compute the composition of the following three functors:
(I)~the functor
\[
Sp_Z \, {j_1}_* \, \mu_{K_1}: D^b (\Rep_F P_1 , H_Q) \longto
                                    D^b (\Delta_1 \text{-} \MHM_F N_{Z / \MKoS})
\]
(see Proposition~\ref{5M}); (II)~the functor
\[
i_0^*: D^b (\Delta_1 \text{-} \MHM_F N_{Z / \MKoS}) \longto
       D^b (\Delta_1 \text{-} \MHM_F Z) \; ,
\]
whose target is equal to $D^b (\MHM_F \MoS)$; (III)~the functor
\[
p_*: D^b (\MHM_F \MoS) \longto D^b (\MHM_F \Mo) \; .
\]
This computation
is complicated by the fact that $p_*$ is neither
left nor right exact --- remember that we are working in a (derived) category
of objects which behave like perverse sheaves, hence there are no exactness
properties for Grothendieck's functors associated to arbitrary morphisms. This is why
we construct a certain factorization of $p_*$,
which will represent it as the composition of a left exact and
a right exact functor (Proposition~\ref{5Ra}).
In order to do so, consider the diagram
\[
\vcenter{\xymatrix@R-10pt{
\MoS = \Delta_1 \backslash Z \ar@{<<-}[r]^-{\tilde{q}} \ar@{->>}[d]^p &
Z \ar@{->>}[d]^{\tilde{p}} \\
\Mo = \Delta \backslash \Mp \ar@{<<-}[r]^-q &
\Mp \\}}
\]
It is \emph{not} Cartesian. However, setting $\bar{Z}:= \HC \backslash Z$
(remember that thanks to Corollary~\ref{1E}, we consider $\HC$ as a
subgroup of $\Delta_1$ in a natural way), we get a natural
factorization
of the morphism $\tilde{q}$, which fits into the diagram
\[
\vcenter{\xymatrix@R-10pt{
\MoS = \Delta_1 \backslash Z \ar@{<<-}[r]^-q \ar@{->>}[d]^p &
\bar{Z}= \HC \backslash Z \ar@{<<-}[r]^-{\bar{q}} \ar@{->>}[d]^{p} &
Z \ar@/_2pc/[ll]^-{\tilde{q}} \ar@{->>}[d]^{\tilde{p}} \\
\Mo = \Delta \backslash \Mp \ar@{<<-}[r]^-q &
\Mp \ar@{=}[r] &
\Mp \\}}
\]
Observe that the left half of this diagram \emph{is} Cartesian,
and that the morphisms $q$ are finite Galois coverings, with Galois
group $\Delta$. We identify $D^b (\MHM_F \MoS)$ with
$D^b (\Delta \text{-} \MHM_F \bar{Z})$, and the functor
\[
p_*: D^b (\MHM_F \MoS) \longto D^b (\MHM_F \Mo)
\]
with its $\Delta$-equivariant version
\[
p_*: D^b (\Delta \text{-} \MHM_F \bar{Z}) \longto
D^b (\Delta \text{-} \MHM_F \Mp) \; .
\]
In order to study this last functor,
recall the closed covering of $Z$ by the closures $\Ms$ of the
toric strata (Section~\ref{4}). It induces a (finite!) closed covering
of $\bar{Z}$ by the $\Msb$, for $\bar{\sigma} \in \bar{\FT} := \HC
\backslash \FT$. By
Proposition~\ref{4B}, the morphism $\bar{q}$ induces isomorphisms
\[
\Ms \isoto \Msb \; ,
\]
for any $\sigma \in \FT$. Observe moreover that, for any $\bar{I}\in
\bar{\FT}_{\bullet}$  the intersection $\MIb$ is either one stratum
$\Msb$ or the empty set.

\begin{Def} \label{5O}
Let $M$ denote one of the varieties $\bar{Z}$ or $\Mp$. Define the
Abelian category
$(\MHM_F M)^{\bar{\FT}}$ applying Definition \ref{5G} to the closed covering
$\{\Msb\}_{\bar{\sigma}\in \bar{\FT}}$ in the case of $\bar{Z}$ and to
the trivial closed covering
$\{M_{\bar{\sigma }}\}_{\bar{\sigma}\in \bar{\FT}}$ with
$M_{\bar{\sigma }}=\Mp$ for all
$\bar{\sigma}\in \bar{\FT}$, in the case of $\Mp$.
\end{Def}

Now remember the action of the finite group $\Delta$ on our geometric
situation. This group acts on the spaces $\bar{Z}$ and $\Mp$ and on
the set of indexes $\bar{\FT}$, hence on the simplicial schemes
$\bar{Z}\times \bar{\FT}_{\bullet}$ and $\Mp\times
\bar{\FT}_{\bullet}$. Therefore, as in
\ref{3S}~(c), we can define the categories $\Delta \text{-}(\MHM_F
\bar{Z})^{\bar{\FT}}$ and $\Delta \text{-}(\MHM_F
\Mp)^{\bar{\FT}}$. For instance the former is the category of mixed
Hodge modules $\BM$ over the simplicial scheme $\bar{Z}\times
\bar{\FT}_{\bullet}$ together with isomorphisms, for $\gamma \in
\Delta $,
\[
\rho _{\gamma }:\gamma ^{\ast} \BM \isoto \BM,
\]
that satisfy the cocycle condition.\\

Since $\Delta $
respects the stratification of $\bar{Z}$ indexed by $\bar{\FT}$, we
can define the equivariant version of the functors $S_{\bullet}$ and
$\Tot$.  We leave it to the reader to check that the
$\Delta$-equivariant versions of \ref{5Q} and \ref{5R} hold (in the proof
of the analogue of \ref{5Q}, choose a finite open affine covering
$\FV$ closed under the action of the finite group $\Delta$, and
observe that the open subset $U$ occurring in point (2) can be
replaced by the intersection of all its translates under $\Delta$).
In particular, we have:

\begin{Prop} \label{5Ra}
(a) There is a canonical functor
\[
p^{\bar{\FT}}_*:
D^b (\Delta \text{-}(\MHM_F \bar{Z})^{\bar{\FT}}) \longto
D^b (\Delta \text{-}(\MHM_F \Mp)^{\bar{\FT}}) \; .
\]
(b) There is a natural commutative diagram
\[
\vcenter{\xymatrix@R-10pt{
D^b (\MHM_F \MoS) \ar@{=}[d] \ar@/_7pc/[ddd]_{p_*} & \\
D^b (\Delta \text{-}\MHM_F \bar{Z}) \ar[r]^-{S_\bullet} \ar[d]_{p_*} &
             D^b (\Delta \text{-}(\MHM_F \bar{Z})^{\bar{\FT}})
                                      \ar[d]^{p^{\bar{\FT}}_*} \\
D^b (\Delta \text{-}\MHM_F \Mp) \ar@{<-}[r]^-{\Tot} \ar@{=}[d] &
             D^b (\Delta \text{-}(\MHM_F \Mp)^{\bar{\FT}}) \\
D^b (\MHM_F \Mo) &
\\}}
\]
\end{Prop}

Using $\Delta =\Delta _{1}/\HC$ and $\bar{Z}=\HC\backslash Z$, and
a slight generalization of Proposition \ref{5D} we make
the identification
\[
D^b (\Delta \text{-}\MHM_F \bar{Z}) =
D^b (\Delta_1 \text{-}\MHM_F Z) \; .
\]
Since the group $\Delta _{1}$ acts on $\Mp$ (by its
quotient $\Delta $) and on the set $\FT$, it acts also on $\Mp\times
\FT_{\bullet}$. We define the category
$\Delta_1 \text{-} (\MHM_F \Mp)^{\FT}$, in the same way as we have defined
$\Delta \text{-} (\MHM_F \Mp)^{\bar{\FT}}$ but using
the infinite version of Definition~\ref{5G}.
We have the following variant
of Proposition~\ref{5D}:

\begin{Prop} \label{5Rb}
The inverse image
\[
\Delta \text{-}(\MHM_F \Mp)^{\bar{\FT}} \longto
\Delta_1 \text{-} (\MHM_F \Mp)^{\FT}
\]
is an equivalence of categories, which possesses a canonical
pseudo-inverse.
\end{Prop}

\begin{Proof}
The group $\Delta_1$ acts freely and properly on the simplicial scheme
$\Mp \times \FT_{\bullet}$.
Hence so does the subgroup $\HC$ of $\Delta_1$.
The quotient by $\HC$ of $\Mp \times \FT_{\bullet}$ equals
$\Mp \times \bar{\FT}_{\bullet}$.
The action of $\Delta = \Delta_1 / \HC$
on this quotient is free and proper, and
\[
\Delta \backslash (\Mp \times \bar{\FT}_{\bullet}) =
\Delta_1 \backslash (\Mp \times \FT_{\bullet}) \; .
\]
\end{Proof}

We now start to evaluate our functors.
Consider the composition
\[
\nu:= p^{\bar{\FT}}_* \, S_\bullet \, i_0^* \, Sp_Z \, {j_1}_* \, \mu_{K_1}:
D^b (\Rep_F P_1 , H_Q) \longto
D^b (\Delta \text{-}(\MHM_F \Mp)^{\bar{\FT}}) \;
\]
where we use the identification
\[
D^b (\Delta \text{-}\MHM_F \bar{Z}) =
D^b (\Delta_1 \text{-}\MHM_F Z) \; .
\]
before applying the functor $S_{\bullet}$.
We have a variant of the canonical construction
\[
\mupT : (\Rep_F G_1, \HQ) \longto \Delta_1 \text{-} (\MHM_F \Mp)^{\FT} \; ,
\]
which associates to a representation $\BV_1$ the mixed Hodge module,
whose  component over $\Mp \times \{I\}$ is
\[
\BM_{I} =
\begin{cases}
\mup (\BV_1) \; , \ &\text{ if }\MI \not = \emptyset \; ,\\
0 \; , \ &\text{ if }\MI=\emptyset \; .
\end{cases}
\]
For any increasing map $\tau $ and $I\in \FT_{\bullet}$ with
$J=\FT_{\bullet}(\tau )(I)$, we put $\tau _{I}=\Id$ if $\MI$ is not
empty and zero otherwise. For $\gamma \in \Delta _{1}$ we let
the isomorphisms $\rho_\gamma$ be given by the action of
$\HQ$.\\

The functor $\mupT$ is exact.
As before, we let $c$ denote the codimension of
$\Mo$ in $(M^K)^*$, which is the same as the relative dimension of
the morphism $\pi: \MKo \to \Mp$.

\begin{Prop} \label{5S}
There is a natural commutative diagram
\[
\vcenter{\xymatrix@R-10pt{
D^b (\Rep_F P_1 , H_Q) \ar[rr]^-{\nu [-c]} \ar[d]_{R\Gamma(W_1,\argdot)} &&
         D^b (\Delta \text{-}(\MHM_F \Mp)^{\bar{\FT}}) \ar@{=}[d]   \\
D^b (\Rep_F G_1, \HQ)  \ar[rr]^-{\mupT}  &&
         D^b (\Delta_1 \text{-} (\MHM_F \Mp)^{\FT})
\\}}
\]
\end{Prop}

\begin{Proof}
Let us first determine the cohomology functors
\[
\CH^r \nu: (\Rep_F P_1 , H_Q) \longto \Delta_1 \text{-} (\MHM_F \Mp)^{\FT} \; .
\]
Let $\BV_1$ be in $(\Rep_F P_1 , H_Q)$, and $I \in \FT$.
By Proposition~\ref{5Q}~(b), the component $(\CH^r \nu (\BV_1))_{I}$
is given by
\[
\CH^r (\tilde{p} \circ i_{I})_* \, i_{I}^* \, Sp_Z \, {j_1}_* \,
                                                     \mu_{K_1} (\BV_1) \; .
\]
Now remember the factorization
\[
\vcenter{\xymatrix@R-10pt{
\MKo \ar@{->>}[r]^-{\pi_t} \ar@/_2pc/[rr]_\pi &
M^{\pi_t (K_1)} \ar@{->>}[r]^-{\pi_a} &
\Mp
\\}}
\]
of the morphism $\pi: \MKo \to \Mp$. As explained in Section~\ref{4},
the map $\tilde{p}: Z \to \Mp$ factors through $\pi_a$, and
identifies $Z$ with a closed union of strata in a
torus embedding of $\MKo$ over $M^{\pi_t (K_1)}$.
The above object involves the direct images $(\tilde{p} \circ i_{I})_*$
of individual strata of this torus embedding. Note that the corresponding
direct image of the generic stratum $\MKo$ equals $\pi_{a *} \circ \pi_{t *} = \pi_*$.
Using Corollary~\ref{5B} and the compatibility of $Sp_Z$ with
$i_{I *} \circ i_{I}^*$
\cite[2.30]{Sa}, we see that
\[
\CH^r (\tilde{p} \circ i_{I})_* \, i_{I}^* \, Sp_Z \, {j_1}_* \,
                                                     \mu_{K_1} (\BV_1) = \CH^r \pi_* \, \mu_{K_1} (\BV_1) \;
\]
when $\MI$ is not empty and zero otherwise.
Furthermore, this identification is compatible with the simplicial
structure:
assume that $\FT_{\bullet}(\tau )(I)=J$ for some increasing map $\tau
$ and that $\MI$ is not empty.
Then the morphisms
\[
(\tau )_{I}: (\CH^r \nu (\BV_1))_{J} \longto
(\CH^r \nu (\BV_1))_{I}
\]
correspond to the identity on $\CH^r \pi_* \, \mu_{K_1} (\BV_1)$.
By \cite[Thm.~2.3]{W1}, there is a natural isomorphism
\[
\mup (H^{r+c} (W_1, \BV_1)) \isoto \CH^r \pi_* \, \mu_{K_1} (\BV_1) \; .
\]
In fact, this is the \emph{canonical} isomorphism given by
the universal property of the cohomological derived functor
\cite[II.2.1.4]{V1}, and by the fact that
the functors on the right hand side are \emph{effa\c{c}able} for $r > -c$.
Since this isomorphism is compatible with automorphisms of Shimura data, we see
that the natural actions of $\Delta_1$ on both sides are compatible.
This proves the claim after passage to the cohomology objects.
We see in particular that the functor $\CH^{-c} \nu$ is left exact, and that
its total right derived functor is equal to
\[
\mupT \circ R\Gamma(W_1,\argdot) \; .
\]
Let us assume for a moment that $\nu$ admits a natural \emph{$f$-lifting}
in the sense of \cite[Def.~A.1~(c)]{B} to the
filtered bounded derived categories (see below).
By \cite[A.7]{B}, this $f$-lifting induces a natural transformation
\[
\eta: \mupT \circ R\Gamma(W_1,\argdot) \longto \nu[-c]
\]
of triangulated functors. Furthermore, the $H^n \eta$ are the natural
transformations corresponding to the universal property of the
cohomological derived functor.
Since we already know that these are isomorphisms, we get the desired
conclusion.

It remains to construct the natural extension
\[
\nu: DF^b (\Rep_F P_1 , H_Q) \longto
DF^b (\Delta \text{-}(\MHM_F \Mp)^{\bar{\FT}})
\]
of $\nu$ satisfying the conditions of \cite[Def.~A.1~(c)]{B}. Since
$\nu$ is a composition of functors, we need to define such an extension
for each of them. For the exact functors $\mu_{K_1}$, ${j_1}_*$, and
$Sp_Z$, there is no problem. For $S_\bullet$, we have Remark~\ref{5H}.
It remains to consider $i_0^*$ and $p^{\bar{\FT}}_*$.
For the construction of
\[
i_0^*: DF^b (\MHM_F N_{\MoS / \MKS}) \longto DF^b (\MHM_F \MoS) \; ,
\]
observe first that by the construction of ${i_0}_* \, i_0^*$
recalled earlier, this latter functor admits a filtered version:
\[
{i_0}_* \, i_0^*: DF^b (\MHM_F N_{\MoS / \MKS}) \longto
DF^b (\MHM_F N_{\MoS / \MKS}) \; .
\]
In fact, the image of ${i_0}_* \, i_0^*$ is contained in
$DF^b_{\MoS} (\MHM_F N_{\MoS / \MKS})$, the full triangulated sub-category of
$DF^b (\MHM_F N_{\MoS / \MKS})$ of filtered complexes $(\BM, F^\bullet \BM)$
with support in $\MoS$, i.e., for which the cohomology objects of all
$F^r \BM$ are supported in $\MoS$. It remains to show that the functor
\[
{i_0}_*: DF^b (\MHM_F \MoS) \longto DF^b_{\MoS} (\MHM_F N_{\MoS / \MKS})
\]
(which exists since the unfiltered ${i_0}_*$ is exact) is an equivalence
of categories. For this, we need to check ($\alpha$)~full faithfulness and
($\beta$)~essential surjectivity.
For ($\alpha$), let $\BM$ and $\BN$ be two objects
of $DF^b (\MHM_F \MoS)$. In order to show that
\[
{i_0}_*: \Hom_{\MoS} (\BM,\BN) \longto
                   \Hom_{N_{\MoS / \MKS}} ({i_0}_* \BM,{i_0}_* \BN)
\]
is an isomorphism, we may, using the exact triangles associated to
the (finite!) filtrations of both
$\BM$ and $\BN$, suppose that these are concentrated in single degrees,
say $m$ and $n$. The same is then true for the filtrations of
${i_0}_* \BM$ and ${i_0}_* \BN$.
By \cite[Def.~A.1~(a)~(iii)]{B}, there are no non-trivial morphisms if $m > n$.
Furthermore, loc.\ cit.\ allows to reduce the case $m \le n$ to the case
$m=n$. But then the morphisms can be calculated in the unfiltered derived
categories \cite[Def.~A.1~(c)]{B}, and the claim follows from
\cite[(4.2.10)]{Sa}. For ($\beta$), we use induction on the length of the
filtration of a given object $\BM$ in $DF^b_{\MoS} (\MHM_F N_{\MoS / \MKS})$.
If the filtration is concentrated in a single degree, use \cite[Def.~A.1~(c)]{B}
and \cite[(4.2.10)]{Sa}. If not, then $\BM$ is a cone of a morphism
$\BM'' \to \BM'[1]$ in $DF^b_{\MoS} (\MHM_F N_{\MoS / \MKS})$ of two
objects in the image of ${i_0}_*$. By ($\alpha$), this morphism comes
from a morphism $f$ in $DF^b (\MHM_F \MoS)$. Thus there is an isomorphism
between $\BM$ and the image under ${i_0}_*$ of a cone of $f$.

For the construction of
\[
p^{\bar{\FT}}_*:
DF^b (\Delta \text{-}(\MHM_F \bar{Z})^{\bar{\FT}}) \longto
DF^b (\Delta \text{-}(\MHM_F \Mp)^{\bar{\FT}}) \; ,
\]
observe first that the functor
$(\BM, F^\bullet \BM) \mapsto (\BM, \BM / F^\bullet \BM)$
identifies the filtered (derived) category of complexes of Hodge modules
with the co-filtered (derived) category. This latter point of view will be
better adapted to our needs. Next, fix a finite $\Delta$-equivariant
open affine covering
$\FV = \{ V_1,\ldots,V_r \}$ of $\bar{Z}$ as in the proof of
Proposition~\ref{5Ra}.
Now imitate the proof of Proposition~\ref{5Q},
using the following observation (see the proof of \cite[Thm.~4.3]{Sa}): assume given a diagram
of bounded complexes of $\Delta$-equivariant Hodge modules on $\bar{Z}$
\[
\vcenter{\xymatrix@R-10pt{
            & \BN_0^\bullet \ar[d]^{\varphi_0} \\
\BM_1^\bullet \ar[r]^-{f_M} & \BM_0^\bullet
\\}}
\]
where
\begin{enumerate}
\item[(1)] the morphism $f_M: \BM_1^\bullet \to \BM_0^\bullet$ is epimorphic
in all degrees,
\item[(2)] the components of $\BN_0^\bullet$ are
$p_*$-acyclic with respect to $\FV$,
\item[(3)] the morphism $\varphi_0: \BN_0^\bullet \to \BM_0^\bullet$
is epimorphic in all degrees, and becomes an isomorphism in
$D^b (\MHM_F \bar{Z})$.
\end{enumerate}
Then this diagram can be
completed in the following way:
\[
\vcenter{\xymatrix@R-10pt{
\BN_1^\bullet \ar[r]^-{f_N} \ar[d]_{\varphi_1} &
                                      \BN_0^\bullet \ar[d]^{\varphi_0} \\
\BM_1^\bullet \ar[r]^-{f_M} & \BM_0^\bullet
\\}}
\]
where
\begin{enumerate}
\item[(4)] the morphism $f_N$ is epimorphic in all degrees,
\item[(5)] the components of $\BN_1^\bullet$ are
$p_*$-acyclic with respect to $\FV$,
\item[(6)] the morphism $\varphi_1$
is epimorphic in all degrees, and becomes an isomorphism in
$D^b (\MHM_F \bar{Z})$.
\end{enumerate}

\end{Proof}

We can now complete the proof of our main result:\\

\begin{Proofof}{Theorem~\ref{2D}}
By Propositions~\ref{5K}, \ref{5L}, \ref{5M}, \ref{5Ra}, and
\ref{5S}, all that remains to be proved is that there is
a natural commutative diagram
\[
\vcenter{\xymatrix@R-10pt{
D^b (\Rep_F G_1, \HQ)
\ar[r]^-{\mupT} \ar[d]_{R\Gamma(\HC,\argdot)} &
         D^b (\Delta_1 \text{-} (\MHM_F \Mp)^{\FT}) \ar[d]^{\Tot}    \\
D^b (\Rep_F G_1, \HQ/\HC) \ar[r]^-{\mup} &
         D^b (\Delta \text{-} (\MHM_F \Mp))
\\}}
\]
Recall that we identify $\Delta_1 \text{-} (\MHM_F \Mp)^{\FT}$
and $\Delta \text{-} (\MHM_F \Mp)^{\bar{\FT}}$,
as well as $\Delta \text{-} (\MHM_F \Mp)$
and $\MHM_F \Mo$, and that the functor $\Tot$ is formed with respect to the
stratification $\bar{\FT}$ (not with respect to $\FT$).

\forget{
In order to prove the claim, recall that
by Proposition~\ref{3R}, we have
\[
\mup \circ R\Gamma(\HC,\argdot) = R\Gamma(\HC,\argdot) \circ \mup \; ,
\]
where
\[
R\Gamma(\HC,\argdot): D^b(\Delta_1 \text{-} (\MHM_F \Mp)) \longto
D^b(\Delta \text{-} (\MHM_F \Mp))
\]
is defined as in Variant~\ref{3S}~(c).}

This is where the conditions \cite[(2.3.1--3)]{P2} listed in the beginning of
Section~\ref{4} enter. For any $\sigma \in \FT$ we denote by $\starT(\sigma )$
the union of $\tau ^{\circ}$ for all $\tau \in \FT$ such that $\sigma $ is
a face of $\tau $. Then the $\starT(\sigma )$ form an open covering of
the set $D$. Moreover, all these open sets are contractible and the
intersection of a finite number of them is also contractible
\cite[Lemma~(2.4.1)]{P2}. We denote
by $C_{\bullet}(\{\starT(\sigma )\},\BZ)$ the \v{C}ech chain complex
associated to this covering (i.e. the dual of the usual \v{C}ech
cochain complex). Since the set $D$ is contractible, the natural
augmentation
\begin{displaymath}
  C_{\bullet}(\{\starT(\sigma )\},\BZ)\longrightarrow \BZ
\end{displaymath}
is a resolution. The group $\Delta _{1}$ acts freely and
properly on the set $\FT$. Therefore, the Abelian groups
$C_{p}(\{\starT(\sigma
)\},\BZ)$ have a natural structure of free $\BZ \Delta _{1}$-modules.
Moreover, since the combinatorics of the open covering
$\{\starT(\sigma )\}$ of $D$ agrees with that of the closed covering
$\{\Ms\}$ of $Z$, by the definition of $\mupT$, the
composition of functors $\Tot\circ \mupT$ agrees with the composition
of functors
\begin{displaymath}
  \mup\circ (\argdot)^{\HC}\circ \Hom(C_{\bullet}(\{\starT(\sigma
  )\},\BZ),\argdot),
\end{displaymath}
which by Proposition \ref{3Ta} agrees with
\begin{displaymath}
  \mup\circ R\Gamma (\HC,\argdot).
\end{displaymath}
\end{Proofof}

\begin{Rem}
(a) Our proof uses a \emph{choice} of toroidal compactification.
However, as can
be seen by passing to simultaneous refinements of two cone decompositions,
the isomorphism of Theorem~\ref{2D} does not depend on this choice.
We leave the details of the proof to the reader. \\[0.2cm]
(b) We also leave it to the reader to formulate and prove results like
\cite[Prop.~(4.8.5)]{P2} on the behaviour of the isomorphism of \ref{2D}
under change of the subgroup $K \subset G(\BA_f)$, and of
the element $g \in G(\BA_f)$.
\end{Rem}


\bigskip
%
%

\end{document}